\documentclass[times, final, 10pt]{elsarticle}

\usepackage[english]{babel}
\usepackage[utf8]{inputenc}
\usepackage{macros}
\usepackage{macros1D}

\usepackage{mathtools}
\usepackage[ruled,vlined]{algorithm2e}
\usepackage{floatrow}
\usepackage{placeins}

\usepackage{placeins}

\usepackage[tmargin=1.0in,bmargin=1.0in, lmargin=1.0in, rmargin=1.0in]{geometry}

\hypersetup{
    colorlinks=true,      		
    linkcolor=blue,       		
    citecolor=blue,       		
    filecolor=black,      		
    urlcolor=red,       		
    bookmarks=false,
    pdffitwindow=true,
    pdfpagelayout=SinglePage
}

\usepackage{xcolor}
\definecolor{legendBlue}{rgb}{0.04705882352, 0.02745098039, 0.5294117647}
\definecolor{legendViolet}{rgb}{0.415686275, 0., 0.658823529}
\definecolor{legendMagenta}{rgb}{0.694117647, 0.164705882, 0.564705882}
\definecolor{legendRed}{rgb}{0.882352941, 0.02745098039, 0.384313725}
\definecolor{legendOrange}{rgb}{0.988235294, 0.650980392, 0.211764706}
\definecolor{legendYellow}{rgb}{0.941176471, 0.976470588, 0.129411765}

\usepackage{tikz}
\usetikzlibrary{arrows.meta}
\usetikzlibrary{calc}
\usetikzlibrary{decorations.markings}

\usepackage{mathpazo} 

\usepackage[displaymath, mathlines,running]{lineno}

\stackMath
\newcommand\reallywidereversedhat[1]{%
\savestack{\tmpbox}{\stretchto{%
  \scaleto{%
    \scalerel*[\widthof{\ensuremath{#1}}]{\kern-.6pt\bigvee\kern-.6pt}%
    {\rule[-\textheight/2]{1ex}{\textheight}}
  }{\textheight}%
}{0.5ex}}%
\stackon[1pt]{#1}{\tmpbox}%
}

\newcommand{\setofindices}{\nabla}
\newcommand{\physicaltree}[1]{R(#1)}
\newcommand{\leaves}[1]{L(#1)}
\newcommand{\physicalleaves}[1]{S(#1)}
\newcommand{\multiresolutiontransform}{\mathcal{M_R}}
\newcommand{\myrightleftarrows}[1]{\mathrel{\substack{\xrightarrow{#1} \\[-.9ex] \xleftarrow{#1}}}}
\newcommand{\collided}{\star}

\usepackage{soul}

\usepackage{calc}
\newsavebox\CBox
\newcommand\hcancel[2][0.5pt]{%
  \ifmmode\sbox\CBox{$#2$}\else\sbox\CBox{#2}\fi%
  \makebox[0pt][l]{\usebox\CBox}%
  \rule[0.5\ht\CBox-#1/2]{\wd\CBox}{#1}}

\begin{document}

\begin{frontmatter}

\title{Multidimensional fully adaptive lattice Boltzmann methods with error control based on multiresolution analysis}
\address[CMAP]{CMAP, CNRS, \'Ecole polytechnique, Institut Polytechnique de Paris, 91128 Palaiseau Cedex, France.}
\address[ORSAY]{Institut de Mathématique d'Orsay,  Université Paris-Saclay, 91405 Orsay Cedex, France.}
\author[CMAP]{Thomas Bellotti}
\ead{thomas.bellotti@polytechnique.edu}
\author[CMAP]{Lo\"{\i}c Gouarin}
\ead{loic.gouarin@polytechnique.edu}
\author[ORSAY]{Benjamin Graille}
\ead{benjamin.graille@universite-paris-saclay.fr}
\author[CMAP]{Marc Massot}
\ead{marc.massot@polytechnique.edu}





\begin{abstract}
Lattice-Boltzmann methods are known for their simplicity, efficiency and ease of parallelization, usually relying on uniform Cartesian meshes with a strong bond between spatial and temporal discretization.
This fact complicates the crucial issue of reducing the computational cost and the memory impact by automatically coarsening the grid where a fine mesh is unnecessary, still ensuring the overall quality of the numerical solution through error control. This work provides a possible answer to this interesting question, by connecting, for the first time, the field of lattice-Boltzmann Methods (LBM) to the adaptive multiresolution (MR) approach based on wavelets.
To this end, we employ a MR multi-scale transform to adapt the mesh as the solution evolves in time according to its local regularity. 
The collision phase is not affected due to its inherent local nature and because we do not modify the speed of the sound, contrarily to most of the LBM/Adaptive Mesh Refinement (AMR) strategies proposed in the literature, thus preserving the original structure of
any LBM scheme. Besides, an original use of the MR allows the scheme to resolve the proper physics by efficiently controlling the accuracy of the transport phase.
We carefully test our method to conclude on its adaptability to a wide family of existing lattice Boltzmann schemes, treating both hyperbolic and parabolic systems of equations, thus being less problem-dependent than the AMR approaches, which have a hard time guaranteeing an effective control on the error. 
The ability of the method to yield a very efficient compression rate and thus a computational cost reduction for solutions involving localized structures with loss of regularity is also shown, while guaranteeing a precise control on the approximation error introduced by the spatial adaptation of the grid. The numerical strategy is implemented on a specific open-source platform called \texttt{SAMURAI} with a dedicated data-structure relying on set algebra.
\end{abstract}

\end{frontmatter}



\section{Introduction}


Steep fronts and shocks are omnipresent structures in natural phenomena and can be modeled by very diverse mathematical models ranging from hyperbolic systems of conservation laws to hyperbolic-parabolic or parabolic systems. 
One finds them in a large variety of phenomena including both gas dynamics, when studying the air-flow around solid objects such as airplanes at relatively large speeds, and combustion and detonation waves in various devices or pollutant
transport in the atmosphere  and  even in astrophysics. 

In these physical phenomena, there are small areas where all the variation of the solution is concentrated -- such as shocks -- and a vast remaining part of the system where the solution varies smoothly or remains almost constant in large \emph{plateaux}.
Therefore, it is interesting to study how to exploit this property of the system in order to reduce both the computational cost of the resolution and to decrease the memory footprint.
This comes from the fact that in the areas where the solution varies gently, we can sample the solution coarsely without giving up on the overall quality and accuracy and we can expect to perform less computations.

Reducing the size of the spatial mesh can basically rely on three approaches.
The first one consists in creating a fixed adapted mesh based on some \emph{a priori} knowledge of the solution and which will not evolve with time. This simple approach is not recommended since it is highly problem-dependent, far from being optimal and can hardly provide an estimate on the additional error introduced by the adaptive mesh. There are some situations, namely in aeroacoustics, when the impact on the resolved physics is important and can be hard to recover from \cite{gendre2017,feng2020,horstmann2018hybrid}.
A second approach is the so-called Adaptive Mesh Refinement (AMR) approach \cite{berger1989}, where the computational mesh is locally refined or coarsened dynamically in time according to a criterion based on the value of some sensor(s), such as the magnitude of the gradient of the solution. Though this method provides a relatively simple time-adaptive approach, it is built on heuristic arguments and thus prevents, in the vast majority of the cases, from having a direct control on the error.
Finally, the last approach is based on the so-called adaptive multiresolution analysis (MR). The discrete solution is decomposed on a local wavelet basis and the coefficients provide a precise measure of the local regularity of the solution, thus supplying essential information on the necessity of refining or coarsening the mesh.
This approach guarantees a precise control on the additional error introduced by the mesh adaptation and proves to be fundamentally problem-independent contrarily to AMR. 
This comes at the price of a more complex implementation and optimization of the procedure and calls for specific data-structures \cite{muller2012}.

Multiresolution analysis, stemming from the pioneer works of Daubechies \cite{daubechies1988} and Mallat \cite{mallat1989}, has been successfully used in combination with the Finite Volume method by Harten \cite{harten1994, harten1995} for the solution of conservation laws with shocks, without exploiting its full capabilities in terms of cost reduction, but showing that one can recover an error control on the approximated solution. Cohen \emph{et al.} \cite{cohen2003} pushed the method forward by devising a fully adaptive strategy which reduces both the computational cost and the memory demand of the method.
Such an approach has had a strong impact for a large spectrum of mathematical models and applications (see the non exhaustive list of examples 
 \cite{bramkamp2004, roussel2005, lamby2005, coquel2006, burger2008, coquel2010, duarte2011,duarte2013,duartecandel2013,duarte2014,duarte2015,nguessan2019} for hyperbolic, parabolic, mixed and coupled systems).
Interestingly, even if such comparisons are a difficult task and the results should be handled carefully, the technique has been compared with the AMR \cite{deiterding2016} where it proves to yield better compression rates. The major feature of the method is the error control but the difficulties concerning an efficient implementation are already pointed out.


In this work, we concentrate on the class of \lb methods, first introduced by McNamara and Zanetti \cite{mcnamara1988} and which are capable of solving many different types of Partial Differential Equations.
In all likelihood, the most well-known \lb scheme is the family of nine-velocities 2D scheme  (see for example \cite{lallemand2000}), used to approximate the solution of the incompressible Navier-Stokes system in the limit of low Mach numbers.
Overall, a large variety of such schemes exists in the literature and a general view on the possible approaches and applications is far from the scope of this work. 
Still, we want to emphasize the fact that the strategy developed in this paper can be \emph{a priori} applied to any \lb scheme.

The inherent structure of the \lb methods comes from a discretization of a mesoscopic equation, with a high level of similarity to 
the Boltzmann equation, with a small family of discrete velocities which are multiple of the so-called lattice velocity, linked to the space and time discretizations.
Therefore, the advection of the particles can be discretized exactly and the collision phase is kept explicit and local by seeing it as a linear relaxation towards the equilibrium distributions, coherently with the so-called BGK approximation of the collision operator as for the Boltzmann equation.
The major strengths of the \lb methods are their {simplicity and} computational performance, since they are fully explicit and they lead to efficient parallelization and implementation on GPUs \cite{schonherr2011, obrecht2013, lin2013}.

Until today, the available mesh adaptation strategies for \lb methods were only based either on a fixed mesh which does not evolve in time \cite{filippova1998, lin2000, kandhai2000, yu2002, dupuis2003} or on the AMR approach \cite{crouse2003, eitel2013, fakhari2014, fakhari2015, fakhari2016}, with the previously highlighted strengths and weaknesses.
Due to the special relation between space and time discretization on the uniform lattice, one should pay special attention to the way of performing the time advancement. 
Besides the way of constructing the mesh, we can identify two main trends:
\begin{itemize}
    \item Methods using local time steps for each level of refinement, thus needing spatial and temporal interpolations and an adaptation of the collision phase. We can cite those acting on fixed grids \cite{filippova1998, lin2000, kandhai2000, dupuis2003} with nested patches and \cite{yu2002} with patches communicating only through edges. Rohde \emph{et al.} \cite{rohde2006} employ the previous approach utilizing volumes to simplify the issue of conservation, which is a feature we retain for our strategy. The same method is combined with AMR in \cite{eitel2013} with sensors based on the vorticity and on the matching of a free-flow total energy and by \cite{crouse2003} with a criterion on the divergence of the velocity field.
    
    The advantage of this procedure is that it minimizes the number of time steps but the shortcomings are the need of an adaptation of the collision phase, with possibly singular parameters, and the need of interpolation which calls for a massive storage of the solution at previous time steps.
    
    \item Methods using a global time step given by the finest space step and no need to adapt the collision. This is the strategy we shall employ in the work. We can cite \cite{fakhari2014, fakhari2015}, who combine this with the AMR where the sensor is again based on the vorticity or use it \cite{fakhari2016} to simulate diphasic flows with sensor on the gradient of the phase-field.
    In these works, the authors employ a Lax-Wendroff scheme for the advection phase where the adaptation to the local level of refinement is done by modifying the CFL number.
    The collision remains untouched and performed locally.
    This method is simpler to implement, more flexible and needs less storage than the previous approach.
    
\end{itemize}

In a companion paper \cite{bellotti2021}, we have made a first theoretical step in order to bridge the gap, in the one-dimensional context, between multiresolution analysis and \lb schemes, bringing in a proof of concept and the potential ability to address a large variety of schemes. However, several crucial questions remain. First, the extension to multi-dimension is far from obvious and is required in order to check the ability to capture properly the resolution of classical test-cases in the \lb community. Second, the multidimensional extension has to rely on an efficient implementation of the method; dynamic mesh handling, even if especially important to reduce the memory footprint, should remain marginal in terms of computational cost. Since this is strongly problem-dependent as demonstrated in \cite{descombes2017}, we will focus on algorithms rather than getting into the field of tailored data-structure and implementation optimization.


In this paper, we introduce a strategy to adapt pre-existing multidimensional \lb methods to be used on non-uniform time-evolving spatial meshes built using the multiresolution analysis.
We keep the procedure as general as possible, ensuring a straightforward adaptation of the pre-existing \lb strategies to address a wide range of applications involving steep and localized fronts.
The developed technique proves, for the first time in the \lb literature, to provide a control on the additional error introduced by the mesh adaptation in time and proves to be successful in coping both with hyperbolic and parabolic equations in several dimensions. 
Successful means here that we reach very interesting compression rates, while guaranteeing a resolution with a small and controlled error with respect to the simulation on the finest mesh. 

Let us insist on the fact that we plainly rely on the original version of any \lb scheme and, combining efficiently with the MR approach, we are able to predict the proper physics with error control for both parabolic and hyperbolic problems. 
The algorithms are well-tailored and designed in such a way as to avoid the recursive nature of the multiresolution transform.
Since the efficient implementation of multiresolution is always a key issue, we rely on a dedicated data-structure and an open-source code \texttt{SAMURAI}.\footnote{\url{https://github.com/hpc-maths/samurai}.} 
Even if the data-structure has been created for the sake of efficient parallelization, this issue is kept out of the scope of the present paper as it is a scientific subject {\sl  per se} and has to be properly investigated within the framework of \lb schemes (see \cite{descombes2017} and \cite{schonherr2011}).


The paper is organized as follows: Sections \ref{sec:SpaceAndTimeDiscretization}, \ref{sec:LatticeBoltzmannOnVolumes} and \ref{sec:MultiresolutionAnalysis} introduce the basic bricks of the hybrid method, namely the concept of nested lattices, the \lb methods on volumes and the multiresolution analysis on Cartesian grids.
Then, these ingredients are put together in the core Section \ref{sec:AdaptiveLatticeBoltzmann}, which elucidates how to construct the mesh in time and how to adapt the existing methods.
The following Section \ref{sec:ImplementationAndOptimization} provides key issues on the implementation of the method and some preliminary paths towards optimization. 
The new strategy is systematically tested in Section \ref{sec:VerificationAndExperiments} on a vectorial four-velocities scheme for hyperbolic conservation laws, on a nine-velocities scheme for the incompressible Navier-Stokes system of equations and finally on a six-velocities scheme to simulate the transport equation in three dimensions.
Here, we show that we actually control the additional error and that the memory needed by the method is effectively reduced.
We conclude in Section \ref{sec:Conclusions}.

\section{Space and time discretization}\label{sec:SpaceAndTimeDiscretization}

Consider a multidimensional bounded domain $\domain \subset \reals^{\spatialdimension}$, where the dimension $\spatialdimension = 1, 2, 3$ is fixed from now on.
For the sake of presentation, we consider hyper-cubic domains of the form $\domain = [0, 1]^{\spatialdimension}$.

\subsection{Nested lattices}


\begin{figure}[h]
    \begin{center}
        \includegraphics[width=0.7\textwidth]{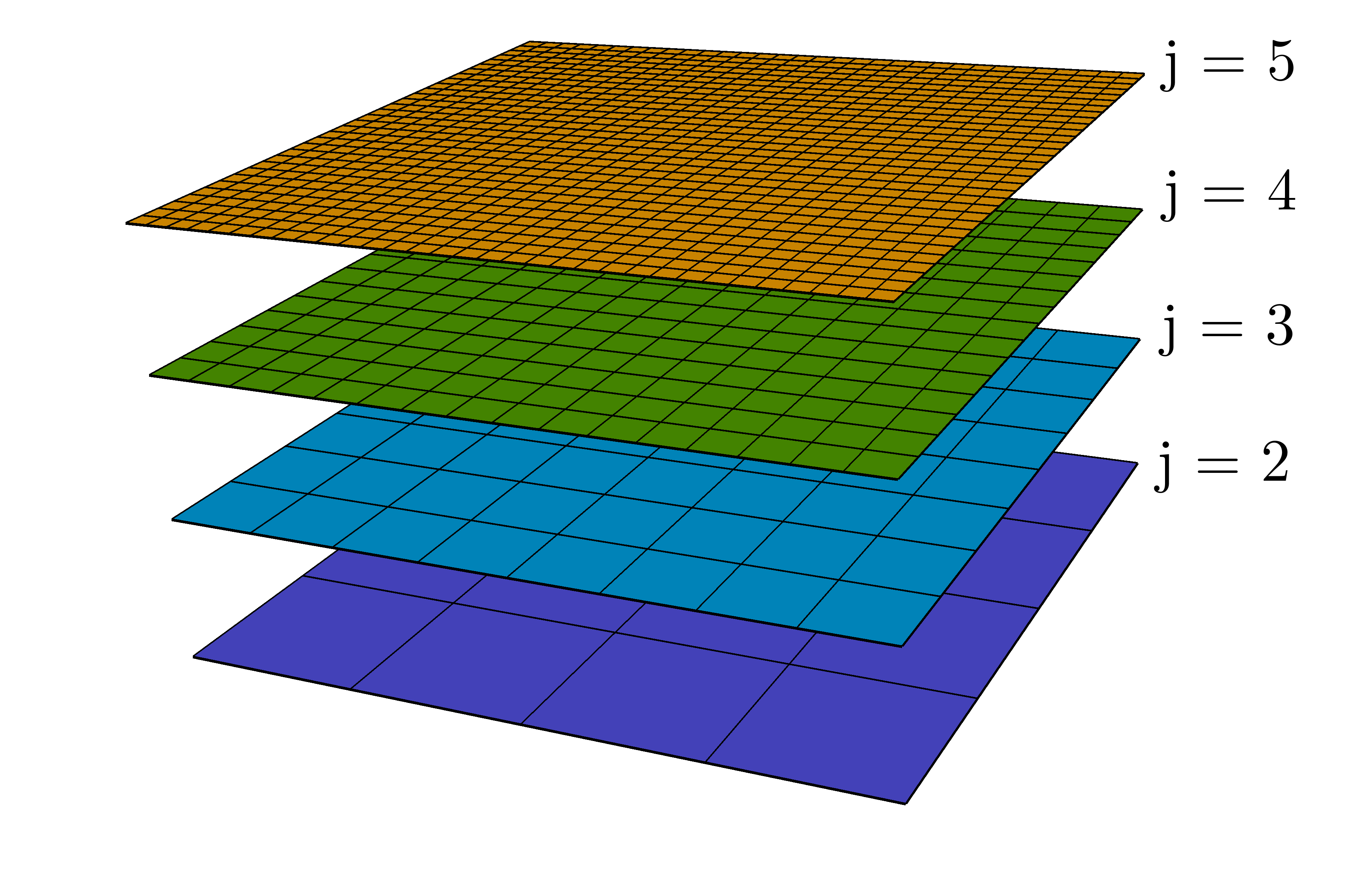}
    \end{center}\caption{\label{fig:nestedlattices}Example of nested dyadic grids for $\spatialdimension = 2$. }
\end{figure}

As done by \cite{cohen2003} and \cite{duarte2011}, the spatial discretization of the domain $\Omega$ is carried out by considering a maximal level of resolution $\maxlevel \in \naturals$, as well as $\levelsnumber + 1$ nested uni-variate dyadic grids with $\levelsnumber \in \naturals$.
Each lattice is identified by its level $\levelletter = \minlevel, \dots, \maxlevel$ with $\minlevel \definitionequality \maxlevel - \levelsnumber$, going from the coarsest to the finest discretization.
The lattices are given by
\begin{linenomath}\begin{equation*}
    \lattice{\levelletter} \definitionequality \subscript{\adaptiveroundbrackets{\subscript{\cellletter}{\levelletter, \vectorial{\indexletter}}}}{\vectorial{\indexletter}}, \qquad \text{with} \quad \subscript{\cellletter}{\levelletter, \vectorial{\indexletter}} \definitionequality \prod_{\dimensionindex = 1}^{\spatialdimension} \adaptivesquarebrackets{ \powertwo{-\levelletter}\subscript{\indexletter}{\dimensionindex},  \powertwo{-\levelletter}\adaptiveroundbrackets{\subscript{\indexletter}{\dimensionindex}+1}},
\end{equation*}\end{linenomath}
with the set of admissible indices $\vectorial{\indexletter} \in \{0, \dots, \numbercells{\levelletter} - 1\}^d$, setting $\numbercells{\levelletter} = \powertwo{\levelletter}$.
These shall be the sets in which $\levelletter$ and $\vectorial{\indexletter}$ are allowed to vary, unless otherwise stated.
As shown in Figure \ref{fig:nestedlattices}, these lattices form a sequence of nested grids where the union of the children cells renders the parent cell
\begin{equation}\label{eq:NestedChildren}
    \bigcup_{\vectorial{\delta} \in \Sigma} C_{\levelletter + 1, 2\vectorial{\indexletter} + \vectorial{\delta}} = C_{\levelletter, \vectorial{\indexletter}}, \qquad \text{with} \qquad \spanningsubcells \definitionequality \power{\adaptivecurlybrackets{0, 1}}{\spatialdimension}, \quad \levelletter = \minlevel, \dots, \maxlevel-1, \quad \text{and} \quad \vectorial{k} \in \{0, \dots, \numbercells{\levelletter} - 1\}^d.
\end{equation}
We stress the fact that the nested structure can be seen in terms of trees (a binary tree for $\spatialdimension = 1$, a quadtree for $\spatialdimension = 2$ and an octree for $\spatialdimension = 3$), even if we shall avoid to explicitly construct it in the sequel because of the intrinsic recursive nature of such a structure.
Thanks to the spatial discretization we have chosen, each cell at the finest level $\maxlevel$ has a edge length $\Delta \spacevariable \definitionequality \ratioobl{1}{\numbercells{\maxlevel}}$. Thus the length of the edge of the cells at level $\levelletter$ and the $\spatialdimension$-dimensional Lebesgue measures of the cells are
\begin{linenomath}\begin{equation*}
    \subscript{\Delta \spacevariable}{\levelletter} := \powertwo{\maxlevel - \levelletter} \Delta x, \qquad
     \adaptiveabs{\subscript{\cellletter}{\levelletter, \vectorial{\indexletter}}}_{\spatialdimension} = \adaptiveroundbrackets{\subscript{\Delta \spacevariable}{\levelletter}}^{\spatialdimension} =  \powertwo{\spatialdimension \adaptiveroundbrackets{\maxlevel - \levelletter}} \power{\adaptiveroundbrackets{\Delta \spacevariable}}{\spatialdimension},
\end{equation*}\end{linenomath}
for $\levelletter = \minlevel, \dots, \maxlevel$ and $\vectorial{k} \in \{0, \dots, \numbercells{\levelletter} - 1\}^d$.

\subsection{Acoustic scaling}

Once the question of the spatial discretization of $\domain$ is successfully handled, we have to specify how do we intend to discretize the time variable.
We have already mentioned that in the standard lattice-Boltzmann methods on uniform grids, the space step $\Delta \spacevariable$ and the time step $\Delta \timevariable$ are indeed linked. This is achieved by considering the so-called lattice velocity $\latticevelocity > 0$ and taking the acoustic scaling\footnote{A parabolic scaling of type $\Delta \timevariable \sim \Delta \spacevariable^2$ is also possible but not covered here.}
\begin{equation}\label{eq:acousticscaling}
    \latticevelocity = \ratio{\Delta \spacevariable}{\Delta \timevariable}.
\end{equation}

So far, we have seen that we have many spatial scales spanned by $\levelletter$, coexisting on the same domain $\domain$, at our disposal.
The approaches stemming from that of Filippova and H\"{a}nel \cite{filippova1998} consider as many local time-steps as $L + 1$ to preserve the same lattice velocity everywhere. This strategies have the following drawbacks: they oblige one to modify the collision phase of the scheme to recover the right viscosity and they require to synchronize the different levels in time with time interpolations and rescalings the variables.
To avoid these complications, we only consider one time scale $\Delta \timevariable$ dictated by the finest spatial scale $\Delta \spacevariable$ through \eqref{eq:acousticscaling}, as also done by Fakhari \emph{et al.} \cite{fakhari2014, fakhari2015, fakhari2016}.
However, we shall emphasize the differences with their strategy in the sequel.
A unique time-step across levels has notable advantages, in particular the fact that the relaxation parameters do not need to be adjusted, avoiding to deal with singular parameters (infinite relaxation times) and the MRT (Multiple Relaxation Times, \cite{dhumieres2002}) case is very naturally handled.
Moreover, we do not have to perform temporal interpolations across the grids and we do not need to store the solution at many different sub time-steps.

\section{Lattice Boltzmann schemes on volumes}\label{sec:LatticeBoltzmannOnVolumes}

Being given spatial and temporal discretizations, we now derive the lattice Boltzmann schemes we employ on volumes.
Indeed, it is most common to interpret the solutions of a \lb scheme in terms of point values rather than averages on volumes \cite{dubois2008lattice, caetano2019}. In this paper, the latter interpretation is utilized  to keep the conservation properties unchanged when adapting the mesh and to be coherent with the most common version of the multiresolution procedure.

In this Section -- since we do not deal with boundary conditions (treated in Section \ref{sec:BoundaryConditions}) -- we temporarily forget to work on a bounded domain $\domain$ by allowing $\vectorial{\indexletter} \in \power{\relatives}{\spatialdimension}$.

\subsection{The origin from the discrete Boltzmann equation}

Consider the discrete Boltzmann equation without external forcing with $q \in \naturals$ discrete velocities spanned by the indices $\populationindex = 0, \dots, \velocitiesnumber - 1$ and collision kernel under the MRT approximation
\begin{linenomath}\begin{equation}\label{eq:BoltzmannEquation}
    \ratio{\partial \superscript{f}{\populationindex}}{\partial \timevariable} + \superscript{\vectorial{\xi}}{\populationindex} \cdot \nabla \superscript{f}{\populationindex} = - \sum_{l = 0}^{\velocitiesnumber - 1} \subscript{\omega}{hl} \adaptiveroundbrackets{\superscript{f}{l} - \superscript{f}{l, \text{eq}}}, \qquad \populationindex = 0, \dots, \velocitiesnumber - 1.
\end{equation}\end{linenomath}
where $\superscript{\vectorial{\xi}}{\populationindex} = \lambda \superscript{\vectorial{\eta}}{\populationindex} $ is the microscopic velocity with corresponding particle distribution function $f^{\populationindex} = f^{\populationindex}(\timevariable, \vectorial{\spacevariable})$, with $\vectorial{\eta}^{\populationindex} \in \mathbb{Z}^d$ being the dimensionless ``logical'' velocity. $ \subscript{\omega}{hl}$ are the entries of the matrix of relaxation frequencies denoted $\operatorial{\omega}$ and $\superscript{f}{l, \text{eq}} = \superscript{f}{l, \text{eq}}(f^0, \dots, f^{q-1})$ are the equilibrium distributions which in general are non linear functions of the conserved moments and thus of the distribution functions.\footnote{The choice $\subscript{\omega}{hl} = \ratioobl{\subscript{\delta}{hl}}{\tau}$ with relaxation time $\tau > 0$ gives the well-known BGK model with single relaxation time.}

The velocity space has already been discretized, thus we still have to handle time and space, as fully described in Appendix 1. The main points are the time integration along the characteristics of $\superscript{\vectorial{\xi}}{\populationindex}$ with a trapezoidal rule to approximate the collision term and the averaging on each cell with commutation of the integral and the collision operator.
Considering a cell  $\subscript{\cellletter}{\maxlevel, \vectorial{\indexletter}}$ for any $\vectorial{k} \in \mathbb{Z}^{\spatialdimension}$ at the finest level of discretization, at the end of the process, the fully discrete scheme reads
\begin{equation}\label{eq:LBMReferenceSchemeNotSplit}
    \subscript{\overline{{f}}}{\maxlevel, \vectorial{\indexletter} + \superscript{\vectorial{\eta}}{\populationindex}}^{\populationindex} \adaptiveroundbrackets{\timevariable + \Delta \timevariable} = \subscript{\overline{{f}}}{\maxlevel, \vectorial{\indexletter}}^{\populationindex} \adaptiveroundbrackets{\timevariable}  - \Delta \timevariable \operatorial{\omega}  (\operatorial{I} + \Delta \timevariable \operatorial{\omega}/2)^{-1} \adaptiveroundbrackets{\overline{\vectorial{f}}_{\maxlevel, \vectorial{\indexletter}}(\timevariable) - \superscript{\vectorial{f}}{\text{eq}}\adaptiveroundbrackets{\superscript{\overline{f}}{0}_{\maxlevel, \vectorial{\indexletter}}(\timevariable), \dots, \superscript{\overline{f}}{q-1}_{\maxlevel, \vectorial{\indexletter}}(\timevariable)} } |_{\populationindex},
\end{equation}
where the bars indicate averages on the considered cell.

\begin{remark}[Not a way of approximating the Boltzmann equation]
 We emphasize the fact that even if \eqref{eq:LBMReferenceSchemeNotSplit} is one of the possible space-time discretizations of the discrete Boltzmann equation  \eqref{eq:BoltzmannEquation}, the link between  \eqref{eq:BoltzmannEquation} and the continuous Boltzmann equation, where the velocities vary in a large continuous set, is in general far from being clear. It is therefore misleading to consider that \eqref{eq:LBMReferenceSchemeNotSplit} is a convergent discretization of a mesoscopic level of description describing a specific physics such as the continuous Boltzmann equation, because $\velocitiesnumber$ remains finite whereas $\Delta x, \Delta t \to 0$. Our way of proceeding is just a manner of proposing a scheme, which can be seen as an algorithmic recipe for solving a macroscopic level of physical modeling.
\end{remark}

\subsection{Collision and transport}

The fully discretized \lb scheme \eqref{eq:LBMReferenceSchemeNotSplit} is made up of a collision and a stream part.
For the sake of clarity and for future use, since they undergo different treatments, we separate them more explicitly.
The link with the operator splitting and the fractional step methods has been studied by Dellar \cite{dellar2013}.
The sub-steps are given by:
\begin{itemize}
    \item \textbf{Collision}, local to each cell and possibly involving the evaluation of non-linear equilibrium functions. This reads
    \begin{linenomath}\begin{equation*}
         \subscript{\overline{\vectorial{f}}}{\maxlevel, \vectorial{\indexletter}}^{\collided} \adaptiveroundbrackets{\timevariable} = \subscript{\overline{\vectorial{f}}}{\maxlevel, \vectorial{\indexletter}} \adaptiveroundbrackets{\timevariable}  - \Delta \timevariable \operatorial{\omega}  (\operatorial{I} + \Delta \timevariable \operatorial{\omega}/2)^{-1} \adaptiveroundbrackets{\overline{\vectorial{f}}_{\maxlevel, \vectorial{\indexletter}}(\timevariable) - \superscript{\vectorial{f}}{\text{eq}}\adaptiveroundbrackets{\superscript{\overline{f}}{0}_{\maxlevel, \vectorial{\indexletter}}(\timevariable), \dots, \superscript{\overline{f}}{q-1}_{\maxlevel, \vectorial{\indexletter}}(\timevariable)}}.
    \end{equation*}\end{linenomath}
    The sign $\collided$ shall denote any post-collisional state.
    To diagonalize the collision matrix using a change of basis $\operatorial{\changeofbasisletter} \in \text{GL}_{\velocitiesnumber}(\mathbb{R})$, it is more common to express it in the space of moments, yielding
    \begin{linenomath}\begin{align}
        \subscript{\overline{\vectorial{m}}}{\maxlevel, \vectorial{\indexletter}}(t) &= \operatorial{\changeofbasisletter} \subscript{\overline{\vectorial{f}}}{\maxlevel, \vectorial{\indexletter}} (t), \nonumber \\ 
        \subscript{\overline{\vectorial{f}}}{\maxlevel, \vectorial{\indexletter}}^{\collided}\adaptiveroundbrackets{\timevariable} &= \power{\operatorial{\changeofbasisletter}}{-1}  \adaptiveroundbrackets{(\operatorial{I} - \operatorial{S})\subscript{\overline{\vectorial{m}}}{\maxlevel, \vectorial{\indexletter}}(t) + \operatorial{S} \superscript{\vectorial{m}}{\text{eq}} \adaptiveroundbrackets{\subscript{\superscript{\overline{m}}{0}}{\maxlevel, \vectorial{\indexletter}}(t), \dots}},\label{eq:CollisionOnFinest}
    \end{align}\end{linenomath}
    where $\operatorial{I}$ is the $q\times q$ identity, $\operatorial{S} = \text{diag} (s_0, \dots, s_{\velocitiesnumber - 1})$ is the matrix of relaxation coefficients such that $s_{\populationindex} \in [0, 2]$ for $\populationindex = 0, \dots, \velocitiesnumber - 1$ and $\vectorial{m}^{\text{eq}}$ are the moments at equilibrium depending on the conserved moments.
    Indeed, in every \lb scheme, there is a certain number of moments which are conserved through the collision phase: their relaxation parameter is conventionally set to zero.
    \item \textbf{Stream}, non-local but (intrinsically) linear
    \begin{equation}\label{eq:StreamFinest}
        \subscript{\superscript{\overline{f}}{\populationindex}}{\maxlevel, \vectorial{\indexletter} + \superscript{\vectorial{\eta}}{\populationindex}} \adaptiveroundbrackets{\timevariable + \Delta \timevariable} = \subscript{\superscript{\overline{f}}{\populationindex, \collided}}{\maxlevel, \vectorial{\indexletter}} \adaptiveroundbrackets{\timevariable }, \qquad \text{or} \qquad  \subscript{\superscript{\overline{f}}{\populationindex}}{\maxlevel, \vectorial{\indexletter}} \adaptiveroundbrackets{\timevariable + \Delta \timevariable} = \subscript{\superscript{\overline{f}}{\populationindex, \collided}}{\maxlevel, \vectorial{\indexletter} - \superscript{\vectorial{\eta}}{\populationindex}} \adaptiveroundbrackets{\timevariable }.
    \end{equation}
    These equation emphasize, in a very compact and elegant fashion, the upwind nature of the stream phase with respect to the direction of each microscopic velocity. From now and for future use on, given $\populationindex = 0, \dots, \velocitiesnumber - 1$, we indicate $\superscript{\populationindex}{\dag}$ the indices such that $\superscript{\vectorial{\eta}}{\superscript{\populationindex}{\dag}} = -\superscript{\vectorial{\eta}}{\populationindex}$.
\end{itemize}

We drop the bar to indicate mean values, assuming that this is understood in the next pages.

\begin{remark}[Vectorial schemes]\label{rem:VectorialSchemes}
    To simulate systems of hyperbolic equations, in particular the Euler system (see Section \ref{sec:EulerTest}), we shall utilize the notion of vectorial scheme, so called in the works of Graille \cite{graille2014} and Dubois \cite{dubois2014}. Vectorial schemes are a subclass of lattice Boltzmann schemes introduced originally in \cite{bartoloni1993lbe}, \cite{shan1997simulation} and \cite{he1998novel} to simulate the Navier-Stokes equations coupled with thermal phenomena. Even if such schemes bear some similarities with discrete velocity kinetic schemes \cite{bouchut1999construction} and \cite{aregba2000discrete}, they are fundamentally different in both the streaming and collision resolution and have been developed in parallel but different communities.
\end{remark}

\section{Multiresolution analysis}\label{sec:MultiresolutionAnalysis}

In the previous Sections, we have presented the volumetric \lb methods on a uniform grid and the framework of nested dyadic grids, without saying how to locally adapt the latter.
We answer this question in the present Section, providing the core concepts of multiresolution analysis on Cartesian meshes. Since this procedure is still static with respect to the time variable $\timevariable$, we omit the time for the sake of clarity.

\subsection{Projection and prediction operators}

        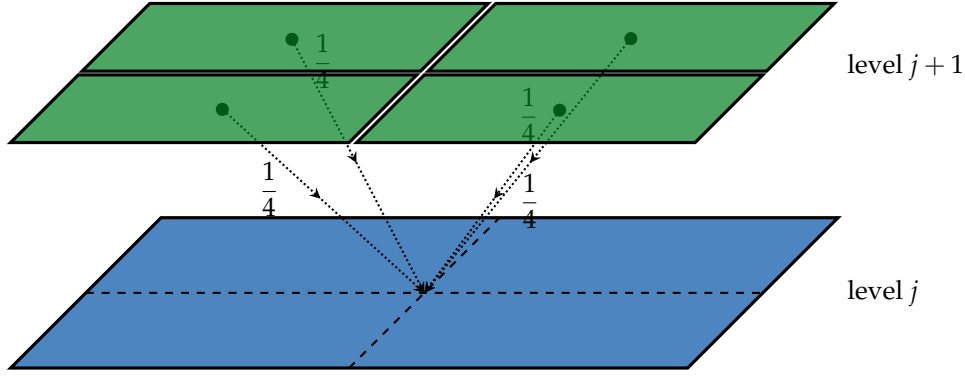
\begin{figure}[h]
        \begin{center}
        \begin{tikzpicture}[x=0.1cm, y=0.1cm]
    \newcount\levdis
    \levdis = 30
    
    \newcount\tol
    \tol = 1
    \definecolor{bluej}{rgb}{0., 0.31764705882, 0.65490196078}

	\coordinate (P1) at (0, 0);
	\coordinate (P2) at (90, 0);
	\coordinate (P3) at (110, 20);
	\coordinate (P4) at (20, 20);
	\fill[fill = bluej, opacity=0.7] (P1) -- (P2) -- (P3) -- (P4) -- cycle; 
    \draw[very thick] (P1) -- (P2) -- (P3) -- (P4) -- cycle; 

	\draw[dashed, thick] (45, 0) -- (65, 20);
	\draw[dashed,thick] (10, 10) -- (100, 10);
	\draw[] (110, 10) node [right] {level $\levelletter$};
	\draw[] (110, 10 + \levdis) node [right] {level $\levelletter + 1$};

    \coordinate (PC) at (55, 10);

    \begin{scope}[decoration={markings, mark=at position 0.5 with {\arrow{stealth}}}] 
    \coordinate (QC) at (27.5, 5 + \levdis);
    \draw[Circle-stealth, thick, densely dotted, postaction = {decorate}] (QC) -- (PC) node [near start, below] {$\ratio{1}{4}$};
    \coordinate (TC) at (27.5 + 45 + \tol, 5 + \levdis);
    \draw[Circle-stealth, thick, densely dotted, postaction = {decorate}] (TC) -- (PC) node [near start, above] {$\ratio{1}{4}$};
    \coordinate (VC) at (27.5 + 45 + \tol + 10 - 0.5*\tol, 5 + \levdis + 10 - 0.5*\tol);
    \draw[Circle-stealth, thick, densely dotted, postaction = {decorate}] (VC) -- (PC) node [midway, below] {$\ratio{1}{4}$};
    \coordinate (UC) at (27.5 + + 10 - 0.5*\tol, 5 + \levdis + 10 - 0.5*\tol);
    \draw[Circle-stealth, thick, densely dotted, postaction = {decorate}] (UC) -- (PC) node [near start, above] {$\ratio{1}{4}$};
    \end{scope}
    
    \definecolor{greenjp1}{rgb}{0., 0.49803921568, 0.12549019607}
    
    \coordinate (Q1) at (0, 0 + \levdis);
	\coordinate (Q2) at (45 , 0 + \levdis);
	\coordinate (Q3) at (55 - \tol, 10 - \tol + \levdis);
	\coordinate (Q4) at (10 - \tol, 10 - \tol + \levdis);
	\fill[fill = greenjp1, fill opacity=0.7] (Q1) -- (Q2) -- (Q3) -- (Q4) -- cycle; 
    \draw[very thick] (Q1) -- (Q2) -- (Q3) -- (Q4) -- cycle;

    \coordinate (T1) at (0 + 45 + \tol, 0 + \levdis);
	\coordinate (T2) at (45 + 45 + \tol, 0 + \levdis);
	\coordinate (T3) at (55 - \tol + 45 + \tol, 10 - \tol + \levdis);
	\coordinate (T4) at (10 - \tol + 45 + \tol, 10 - \tol + \levdis);
	\fill[fill = greenjp1, fill opacity=0.7] (T1) -- (T2) -- (T3) -- (T4) -- cycle; 
    \draw[very thick] (T1) -- (T2) -- (T3) -- (T4) -- cycle;

    \coordinate (V1) at (0 + 45 + \tol + 10 - 0.5*\tol, 0 + \levdis + 10 - 0.5*\tol);
	\coordinate (V2) at (45 + 45 + \tol + 10 - 0.5*\tol, 0 + \levdis + 10 - 0.5*\tol);
	\coordinate (V3) at (55 - \tol + 45 + \tol + 10 - 0.5*\tol, 10 - \tol + \levdis + 10 - 0.5*\tol);
	\coordinate (V4) at (10 - \tol + 45 + \tol + 10 - 0.5*\tol, 10 - \tol + \levdis + 10 - 0.5*\tol);
	\fill[fill = greenjp1, fill opacity=0.7] (V1) -- (V2) -- (V3) -- (V4) -- cycle; 
    \draw[very thick] (V1) -- (V2) -- (V3) -- (V4) -- cycle;

    \coordinate (U1) at (0 + 10 - 0.5*\tol, 0 + \levdis + 10 - 0.5*\tol);
	\coordinate (U2) at (45 + 10 - 0.5*\tol , 0 + \levdis + 10 - 0.5*\tol);
	\coordinate (U3) at (55 - \tol + 10 - 0.5*\tol, 10 - \tol + \levdis + 10 - 0.5*\tol);
	\coordinate (U4) at (10 - \tol + 10 - 0.5*\tol, 10 - \tol + \levdis + 10 - 0.5*\tol);
	\fill[fill = greenjp1, fill opacity=0.7] (U1) -- (U2) -- (U3) -- (U4) -- cycle; 
    \draw[very thick] (U1) -- (U2) -- (U3) -- (U4) -- cycle; 

\end{tikzpicture}
        \end{center}
    \caption[]{\label{fig:projectionOperator} Illustration of the action of the projection operator in the context of $d = 2$. The cell average on the cell at level $\levelletter$ is reconstructed by taking the average of the values on its four children at level $\levelletter + 1$.}
        \end{figure}

Multiresolution analysis is made possible by two basic operations, called ``projection'' and ``prediction'', allowing to pass from a fine to a coarse representation and \emph{vice-versa} through one-level jumps.
We first introduce the projection operator, whose rationale is illustrated on Figure \ref{fig:projectionOperator} and which is used to transform information known on a fine level $\levelletter+1$ to one known on a coarse level $\levelletter$
\begin{definition}[Projection operator]
    Let $\levelletter = \minlevel, \dots, \maxlevel - 1$. The projection operator $\subscript{\operatorial{P}}{\triangledown}: \reals^{2d} \to \reals$, taking data on level $\levelletter+1$ and yielding data on level $\levelletter$ acts in the following way:
    \begin{linenomath}\begin{equation*}
        \superscript{\subscript{f}{\levelletter, \vectorial{\indexletter}}}{\populationindex} = \subscript{\operatorial{P}}{\triangledown} \adaptiveroundbrackets{(\superscript{\subscript{f}{\levelletter + 1, 2 \vectorial{\indexletter} + \vectorial{\delta}}}{\populationindex})_{\vectorial{\delta} \in \Sigma}} \definitionequality \ratio{1}{|\subscript{\cellletter}{\levelletter, \vectorial{\indexletter}}|_{\spatialdimension}} \sum_{\vectorial{\delta} \in \spanningsubcells} |\subscript{\cellletter}{\levelletter + 1, 2 \vectorial{\indexletter} + \vectorial{\delta}}|_{\spatialdimension} ~  \superscript{\subscript{f}{\levelletter + 1, 2 \vectorial{\indexletter} + \vectorial{\delta}}}{\populationindex} = \frac{1}{\powertwo{\spatialdimension}}  \sum_{\vectorial{\delta} \in \spanningsubcells} \superscript{\subscript{f}{\levelletter + 1, 2 \vectorial{\indexletter} + \vectorial{\delta}}}{\populationindex} ~ .
    \end{equation*}\end{linenomath}
\end{definition}
Given a cell, the projection operator is fully local and is the unique operator conserving the average because it takes the average of the values defined on children cells which are nested inside their parent according to \eqref{eq:NestedChildren}.

The prediction operator acts in the opposite direction, thus transforming information known on a coarse level $\levelletter$ to one on a finer level $\levelletter + 1$. This is the crucial ingredient of multiresolution and it is not uniquely defined, because there are many ways of reconstructing the lacking pieces of information ``destroyed'' by the projection operator.
In order to ensure the feasibility of the whole procedure and enforce mass conservation that one expects from a \lb method, some constraints inspired by \cite{cohen2003} have to be imposed
\begin{definition}[Prediction operator]\label{def:PredictionOperator}
    Let $\levelletter = \minlevel, \dots, \maxlevel - 1$. The prediction operator $\subscript{\operatorial{P}}{\vartriangle}: \reals^{1+w} \to \mathbb{R}^{2d}$, used to transform information known on a level $\levelletter$ to a guessed value on level $\levelletter + 1$ (denoted by a hat) fulfills the following properties:
    \begin{itemize}
        \item It is local, namely, the value predicted on a cell $C_{\levelletter + 1, 2 \vectorial{\indexletter} + \vectorial{\delta}}$ for any $\vectorial{\delta} \in \Sigma$ can only depend on the values of cells of level $\levelletter$ belonging to a finite stencil $\mathcal{R} \adaptiveroundbrackets{\levelletter, \vectorial{\indexletter}}$ of $w \in \naturals$ elements around its parent $C_{j, \vectorial{k}}$.
        \item It is consistent with the projection operator, so the operator
        \begin{linenomath}\begin{equation*}
            \subscript{\operatorial{P}}{\triangledown} \composed \subscript{\operatorial{P}}{\vartriangle}: \reals^{1+w} \to \reals,
        \end{equation*}\end{linenomath}
        leaves the value on the parent cell invariant, namely
        \begin{equation}\label{eq:ConsistencyPrediction}
            \adaptiveroundbrackets{\subscript{\operatorial{P}}{\triangledown} \composed \subscript{\operatorial{P}}{\vartriangle}} (f_{\levelletter, \vectorial{\indexletter}}^{\populationindex}, \dots, \dots) = f_{\levelletter, \vectorial{\indexletter}}^{\populationindex}.
        \end{equation}
    \end{itemize}
\end{definition}
The first point ensures that the operator is local around a given cell, whereas the second\footnote{In wavelet theory, this consistency property is interpreted by saying that the mother wavelet is oscillating, with vanishing mean.} guarantees that the discrete integrals are conserved and that the parent necessarily belongs to $\mathcal{R} \adaptiveroundbrackets{\levelletter, \vectorial{\indexletter}}$.

    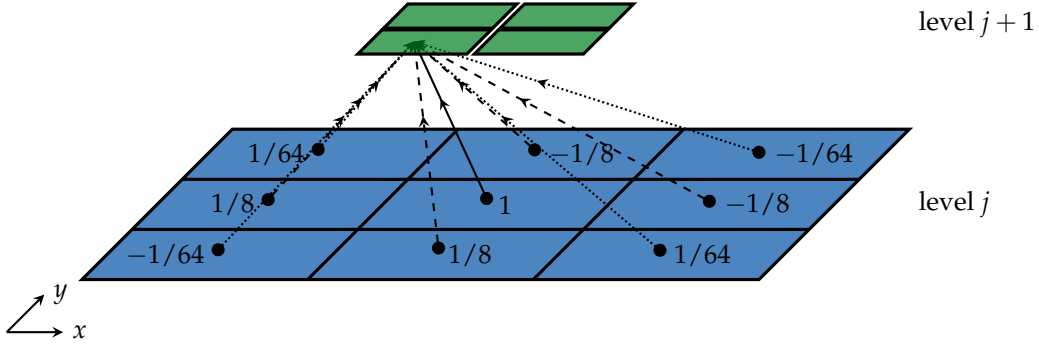
\begin{figure}
            \begin{center}
                   \begin{tikzpicture}[x=0.1cm, y=0.1cm]

    \newcount\levdis
    \levdis = 30
    
    \newcount\tol
    \tol = 1.0
    \definecolor{bluej}{rgb}{0., 0.31764705882, 0.65490196078}
    \definecolor{greenjp1}{rgb}{0., 0.49803921568, 0.12549019607}

    \draw[-stealth, thick] (-10, -7) -- (-2.5, -7); 
    \draw[-stealth, thick] (-10, -7) -- (-5, -2); 
    \coordinate (Px) at (-2.5, -7);
    \coordinate (Py) at  (-5, -2);
    \draw[] (Px) node [right]{$x$};
    \draw[] (Py) node [right]{$y$};

    \coordinate (P1) at (0, 0);
	\coordinate (P2) at (90, 0);
	\coordinate (P3) at (110, 20);
	\coordinate (P4) at (20, 20);
	\fill[fill = bluej, opacity=0.7] (P1) -- (P2) -- (P3) -- (P4) -- cycle; 
    \draw[very thick] (P1) -- (P2) -- (P3) -- (P4) -- cycle; 

    \draw[very thick] (30, 00) -- (50, 20); 
    \draw[very thick] (60, 00) -- (80, 20); 
    \draw[very thick] (6.666, 6.666) -- (6.666 + 90, 6.666);
    \draw[very thick] (2*6.666, 2*6.666) -- (2*6.666 + 90, 2*6.666); 
    
    \begin{scope}[decoration={markings, mark=at position 0.65 with {\arrow{stealth}}}] 
    \coordinate (Target) at (36.666 + 7.5, \levdis + 7.5 * 0.222);
    \coordinate (Central) at (45 + 7.5 * 1.222, 10);
    \draw[Circle-stealth, thick, postaction = {decorate}] (Central) -- (Target);
    \draw[] (Central) node [right]{$1$};
    
    \coordinate (E) at (45 + 7.5 * 1.222 + 30, 10);
    \draw[Circle-stealth, thick, dashed, postaction = {decorate}] (E) -- (Target);
    \draw[] (E) node [right]{$-\ratioobl{1}{8}$};
    \coordinate (W) at (45 + 7.5 * 1.222 - 30, 10);
    \draw[Circle-stealth, thick, dashed, postaction = {decorate}] (W) -- (Target);
    \draw[] (W) node [left]{$\ratioobl{1}{8}$};
    \coordinate (N) at (45 + 7.5 * 1.222 + 30 * 0.222, 16.666);
    \draw[Circle-stealth, thick, dashed, postaction = {decorate}] (N) -- (Target);
    \draw[] (N) node [right]{$-\ratioobl{1}{8}$};
    \coordinate (S) at (45 + 7.5 * 1.222 - 30 * 0.222, 3.333);
    \draw[] (S) node [right]{$\ratioobl{1}{8}$};
    \draw[Circle-stealth, thick, dashed, postaction = {decorate}] (S) -- (Target);
    
    \coordinate (SE) at (45 + 7.5 * 1.222 - 30 * 0.222 + 30, 3.333);
    \draw[] (SE) node [right]{$\ratioobl{1}{64}$};
    \draw[Circle-stealth, thick, densely dotted, postaction = {decorate}] (SE) -- (Target);
    \coordinate (SW) at (45 + 7.5 * 1.222 - 30 * 0.222 + 30 - 60, 3.333);
    \draw[] (SW) node [left]{$-\ratioobl{1}{64}$};
    \draw[Circle-stealth, thick, densely dotted, postaction = {decorate}] (SW) -- (Target);
    \coordinate (NW) at (45 + 7.5 * 1.222 + 30 * 0.222 - 30, 16.666);
    \draw[] (NW) node [left]{$\ratioobl{1}{64}$};
    \draw[Circle-stealth, thick, densely dotted, postaction = {decorate}] (NW) -- (Target);
    \coordinate (NE) at (45 + 7.5 * 1.222 + 30 * 0.222 + 30, 16.666);
    \draw[] (NE) node [right]{$-\ratioobl{1}{64}$};
    \draw[Circle-stealth, thick, densely dotted, postaction = {decorate}] (NE) -- (Target);

    \end{scope}
    \draw[] (110, 10) node [right] {level $\levelletter$};
	\draw[] (110, 4 + \levdis) node [right] {level $\levelletter + 1$};
    
    \coordinate (Q1) at (36.666, \levdis);
	\coordinate (Q2) at (36.666 + 15 - 0.5*\tol, \levdis);
	\coordinate (Q3) at (36.666 + 15 * 1.222 - 1.222*0.5*\tol, \levdis + 0.222*15 - 0.222 *0.5*\tol );
	\coordinate (Q4) at (36.666 + 15 * 0.222 - 0.222*0.5*\tol, \levdis + 0.222*15 - 0.222 *0.5*\tol );
	\fill[fill = greenjp1, opacity=0.7] (Q1) -- (Q2) -- (Q3) -- (Q4) -- cycle; 
    \draw[very thick] (Q1) -- (Q2) -- (Q3) -- (Q4) -- cycle; 

    \coordinate (R1) at (36.666 + 15 + 0.5*\tol, \levdis);
	\coordinate (R2) at (36.666 + 15 - 0.5*\tol + 15 + 0.5*\tol, \levdis);
	\coordinate (R3) at (36.666 + 15 * 1.222 - 1.222*0.5*\tol + 15 + 0.5*\tol, \levdis + 0.222*15 - 0.222 *0.5*\tol );
	\coordinate (R4) at (36.666 + 15 * 0.222 - 0.222*0.5*\tol + 15 + 0.5*\tol, \levdis + 0.222*15 - 0.222 *0.5*\tol );
	\fill[fill = greenjp1, opacity=0.7] (R1) -- (R2) -- (R3) -- (R4) -- cycle; 
    \draw[very thick] (R1) -- (R2) -- (R3) -- (R4) -- cycle;

    \coordinate (S1) at (36.666 + 0.222 * 15 + 0.222 * 0.5*\tol, \levdis + 0.222 * 15 + 0.222 * 0.5*\tol);
	\coordinate (S2) at (36.666 + 15 - 0.5*\tol + 0.222 * 15 + 0.222 * 0.5*\tol, \levdis + 0.222 * 15 + 0.222 * 0.5*\tol);
	\coordinate (S3) at (36.666 + 15 * 1.222 - 1.222*0.5*\tol + 0.222 * 15 + 0.222 * 0.5*\tol, \levdis + 0.222*15 - 0.222 *0.5*\tol + 0.222 * 15 + 0.222 * 0.5*\tol);
	\coordinate (S4) at (36.666 + 15 * 0.222 - 0.222*0.5*\tol + 0.222 * 15 + 0.222 * 0.5*\tol, \levdis + 0.222*15 - 0.222 *0.5*\tol + 0.222 * 15 + 0.222 * 0.5*\tol);
	\fill[fill = greenjp1, opacity=0.7] (S1) -- (S2) -- (S3) -- (S4) -- cycle; 
    \draw[very thick] (S1) -- (S2) -- (S3) -- (S4) -- cycle;

    \coordinate (V1) at (36.666 + 0.222 * 15 + 0.222 * 0.5*\tol + 15 + 0.5*\tol, \levdis + 0.222 * 15 + 0.222 * 0.5*\tol);
	\coordinate (V2) at (36.666 + 15 - 0.5*\tol + 0.222 * 15 + 0.222 * 0.5*\tol + 15 + 0.5*\tol, \levdis + 0.222 * 15 + 0.222 * 0.5*\tol);
	\coordinate (V3) at (36.666 + 15 * 1.222 - 1.222*0.5*\tol + 0.222 * 15 + 0.222 * 0.5*\tol + 15 + 0.5*\tol, \levdis + 0.222*15 - 0.222 *0.5*\tol + 0.222 * 15 + 0.222 * 0.5*\tol);
	\coordinate (V4) at (36.666 + 15 * 0.222 - 0.222*0.5*\tol + 0.222 * 15 + 0.222 * 0.5*\tol + 15 + 0.5*\tol, \levdis + 0.222*15 - 0.222 *0.5*\tol + 0.222 * 15 + 0.222 * 0.5*\tol);
	\fill[fill = greenjp1, opacity=0.7] (V1) -- (V2) -- (V3) -- (V4) -- cycle; 
    \draw[very thick] (V1) -- (V2) -- (V3) -- (V4) -- cycle; 
    
\end{tikzpicture}

                \end{center}
            \caption[]{\label{fig:predictionOperator} Illustration of the action of the prediction operator in the context of $d=2$ for $\predictionstencildepth = 1$. We indicate the value of the weight for each cell in the prediction stencil. In particular, we are predicting on the lower-left sibling and this yields a particular sign of the weights for this cell. The parent cell (weight $1$) is connected through a continuous arrow, whereas neighbors along the axis (weight $\pm \ratioobl{1}{8}$) with a dashed arrow and along the diagonals (weight $\pm \ratioobl{1}{64}$) with a dotted arrow.}
        \end{figure}

Though not compulsory, we consider a linear prediction operator. Using this kind of operator ensures the performance optimization  of the whole algorithm as described in Section \ref{sec:ImplementationAndOptimization} and guarantees a good generality of the procedure.
Let $\predictionstencildepth$ be a non-negative integer and consider the prediction operators generated by polynomial centered interpolations which are exact for polynomials up to degree $2 \predictionstencildepth$, thus being of order $\predictionorder \definitionequality 2 \predictionstencildepth + 1$ and their generalization by tensor product to the multidimensional context by Bihari and Harten \cite{bihari1997multiresolution}:
\begin{itemize}
    \item For $\spatialdimension = 1$, we utilize
    \begin{linenomath}\begin{equation*}
        \predicted{f}{\populationindex}{\levelletter + 1, 2\indexletter + \delta}  = \superscript{\subscript{f}{\levelletter, \indexletter}}{\populationindex} + \power{\adaptiveroundbrackets{-1}}{\delta} \subscript{\superscript{\predictionnoncenteredletter}{\predictionstencildepth}}{1} \adaptiveroundbrackets{\indexletter; \superscript{\subscript{\vectorial{f}}{\levelletter}}{\populationindex}}, \qquad \text{with} \quad \subscript{\superscript{\predictionnoncenteredletter}{\predictionstencildepth}}{1} \adaptiveroundbrackets{\indexletter; \superscript{\subscript{\vectorial{f}}{\levelletter}}{\populationindex}} \definitionequality \sum_{\alpha = 1}^{\predictionstencildepth} c_{\alpha} \adaptiveroundbrackets{\superscript{\subscript{f}{\levelletter, \indexletter + \alpha}}{\populationindex} - \superscript{\subscript{f}{\levelletter, \indexletter - \alpha}}{\populationindex}},
    \end{equation*}\end{linenomath}
        for $\delta = 0, 1$, where the first coefficients are given on Table \ref{tab:PredictionCoefficients} and recovered as presented in Appendix 2.
        \begin{table}
            \begin{center}\caption{\label{tab:PredictionCoefficients} Coefficients used in the polynomial centered prediction operators (taken from Duarte \cite{duarte2011}). The rationale and the procedure to recover these coefficient are presented in Appendix 2.}
                \begin{tabular}{| c | c c c |}
                    \hline
                    $\predictionstencildepth$ & $c_1$ & $c_2$ & $c_3$  \\ 
                    \hline
                    1 & $-\ratioobl{1}{8}$ & - & - \\
                    2 & $-\ratioobl{22}{128}$ & $\ratioobl{3}{128}$ & -\\
                    3 & $-\ratioobl{201}{1024}$ & $\ratioobl{11}{256}$ & $-\ratioobl{5}{1024}$ \\
                    \hline
                \end{tabular}
            \end{center}
        \end{table}

    \item For $\spatialdimension = 2$, we use the same idea than for $\spatialdimension = 1$ along the $x$ and $y$ axis completed by a tensor product term along the diagonals
    \begin{linenomath}\begin{align*}
         \predicted{f}{\populationindex}{\levelletter + 1, 2\vectorial{\indexletter} + \vectorial{\delta}}
          &= \superscript{\subscript{f}{\levelletter, \vectorial{\indexletter}}}{\populationindex} + \power{\adaptiveroundbrackets{-1}}{\subscript{\delta}{1}} \subscript{\superscript{\predictionnoncenteredletter}{\predictionstencildepth}}{1} \adaptiveroundbrackets{\subscript{\indexletter}{1}; \superscript{\subscript{\vectorial{f}}{\levelletter, \adaptiveroundbrackets{\cdot, \subscript{\indexletter}{2}}}}{\populationindex}} + \power{\adaptiveroundbrackets{-1}}{\subscript{\delta}{2}} \subscript{\superscript{\predictionnoncenteredletter}{\predictionstencildepth}}{1} \adaptiveroundbrackets{\subscript{\indexletter}{2}; \superscript{\subscript{\vectorial{f}}{\levelletter, \adaptiveroundbrackets{ \subscript{\indexletter}{1}, \cdot}}}{\populationindex}} - \power{\adaptiveroundbrackets{-1}}{\subscript{\delta}{1} + \subscript{\delta}{2}} \superscript{\subscript{\predictionnoncenteredletter}{2}}{\predictionstencildepth}\adaptiveroundbrackets{\vectorial{\indexletter};\subscript{\superscript{\vectorial{f}}{\populationindex}}{\levelletter}}, \\
         &\text{with} \qquad \\
         &\superscript{\subscript{\predictionnoncenteredletter}{2}}{\predictionstencildepth}\adaptiveroundbrackets{\vectorial{\indexletter};\subscript{\superscript{\vectorial{f}}{\populationindex}}{\levelletter}} \definitionequality \sum_{\alpha, \beta = 1}^{\predictionstencildepth} c_{\alpha} c_{\beta} \adaptiveroundbrackets{\superscript{\subscript{f}{\levelletter, \adaptiveroundbrackets{\subscript{\indexletter}{1} + \alpha, \subscript{\indexletter}{2} + \beta}}}{\populationindex} - \superscript{\subscript{f}{\levelletter, \adaptiveroundbrackets{\subscript{\indexletter}{1} - \alpha, \subscript{\indexletter}{2} + \beta}}}{\populationindex} - \superscript{\subscript{f}{\levelletter, \adaptiveroundbrackets{\subscript{\indexletter}{1} + \alpha, \subscript{\indexletter}{2} - \beta}}}{\populationindex} + \superscript{\subscript{f}{\levelletter, \adaptiveroundbrackets{\subscript{\indexletter}{1} - \alpha, \subscript{\indexletter}{2} - \beta}}}{\populationindex}}.
    \end{align*}\end{linenomath}
    for $\vectorial{\delta} \in \Sigma$.
    The way this formula acts for $\predictionstencildepth = 1$ is illustrated on Figure \ref{fig:predictionOperator}. This will be indeed the choice for all the examples we will show in pages to come.

    \item For $\spatialdimension = 3$, the strategy follows the same rationale than for the previous point
    \begin{linenomath}\begin{gather*}
        \predicted{f}{\populationindex}{\levelletter + 1, 2 \vectorial{\indexletter} + \vectorial{\delta}}
        = \subscript{\superscript{f}{\populationindex}}{\levelletter, \vectorial{\indexletter}} + \power{\adaptiveroundbrackets{-1}}{\subscript{\delta}{1}} \subscript{\superscript{\predictionnoncenteredletter}{\predictionstencildepth}}{1} \adaptiveroundbrackets{\subscript{\indexletter}{1}; \superscript{\subscript{\vectorial{f}}{\levelletter, \adaptiveroundbrackets{\cdot, \subscript{\indexletter}{2}, \subscript{\indexletter}{3}}}}{\populationindex}} + \power{\adaptiveroundbrackets{-1}}{\subscript{\delta}{2}} \subscript{\superscript{\predictionnoncenteredletter}{\predictionstencildepth}}{1} \adaptiveroundbrackets{\subscript{\indexletter}{2}; \superscript{\subscript{\vectorial{f}}{\levelletter, \adaptiveroundbrackets{ \subscript{\indexletter}{1}, \cdot, \subscript{\indexletter}{3}}}}{\populationindex}} + \power{\adaptiveroundbrackets{-1}}{\subscript{\delta}{3}} \subscript{\superscript{\predictionnoncenteredletter}{\predictionstencildepth}}{1} \adaptiveroundbrackets{\subscript{\indexletter}{3}; \superscript{\subscript{\vectorial{f}}{\levelletter, \adaptiveroundbrackets{ \subscript{\indexletter}{1}, \subscript{\indexletter}{2}, \cdot}}}{\populationindex}} \\
        - \power{\adaptiveroundbrackets{-1}}{\subscript{\delta}{1} + \subscript{\delta}{2}} \subscript{\superscript{\predictionnoncenteredletter}{\predictionstencildepth}}{2} \adaptiveroundbrackets{\adaptiveroundbrackets{\subscript{\indexletter}{1}, \subscript{\indexletter}{2}}; \superscript{\subscript{\vectorial{f}}{\levelletter, \adaptiveroundbrackets{\cdot, \cdot, \subscript{\indexletter}{3}}}}{\populationindex}} - \power{\adaptiveroundbrackets{-1}}{\subscript{\delta}{1} + \subscript{\delta}{3}} \subscript{\superscript{\predictionnoncenteredletter}{\predictionstencildepth}}{2} \adaptiveroundbrackets{\adaptiveroundbrackets{\subscript{\indexletter}{1}, \subscript{\indexletter}{3}}; \superscript{\subscript{\vectorial{f}}{\levelletter, \adaptiveroundbrackets{\cdot, \subscript{\indexletter}{2}, \cdot}}}{\populationindex}} - \power{\adaptiveroundbrackets{-1}}{\subscript{\delta}{2} + \subscript{\delta}{3}} \subscript{\superscript{\predictionnoncenteredletter}{\predictionstencildepth}}{2} \adaptiveroundbrackets{\adaptiveroundbrackets{\subscript{\indexletter}{2}, \subscript{\indexletter}{3}}; \superscript{\subscript{\vectorial{f}}{\levelletter, \adaptiveroundbrackets{ \subscript{\indexletter}{1}, \cdot, \cdot}}}{\populationindex}} \\
        + \power{\adaptiveroundbrackets{-1}}{\subscript{\delta}{1} + \subscript{\delta}{2} + \subscript{\delta}{3}} \subscript{\superscript{\predictionnoncenteredletter}{\predictionstencildepth}}{3} \adaptiveroundbrackets{\vectorial{\indexletter}; \superscript{\subscript{\vectorial{f}}{\levelletter}}{\populationindex}},
    \end{gather*}\end{linenomath}
    for $\vectorial{\delta} \in \Sigma$, where we have set
    \begin{linenomath}\begin{align*}
        \subscript{\superscript{\predictionnoncenteredletter}{\predictionstencildepth}}{3} \adaptiveroundbrackets{\vectorial{\indexletter}; \superscript{\subscript{\vectorial{f}}{\levelletter}}{\populationindex}} \definitionequality \sum_{\alpha, \beta, \tau = 1}^{\predictionstencildepth} c_{\alpha} c_{\beta} c_{\tau} &\bigg ( \superscript{\subscript{f}{\levelletter, \adaptiveroundbrackets{\subscript{\indexletter}{1} + \alpha, \subscript{\indexletter}{2} + \beta, \subscript{\indexletter}{3} + \tau}}}{\populationindex} - \superscript{\subscript{f}{\levelletter, \adaptiveroundbrackets{\subscript{\indexletter}{1} - \alpha, \subscript{\indexletter}{2} + \beta, \subscript{\indexletter}{3} + \tau}}}{\populationindex} - \superscript{\subscript{f}{\levelletter, \adaptiveroundbrackets{\subscript{\indexletter}{1} + \alpha, \subscript{\indexletter}{2} - \beta, \subscript{\indexletter}{3} + \tau}}}{\populationindex} \\
        & - \superscript{\subscript{f}{\levelletter, \adaptiveroundbrackets{\subscript{\indexletter}{1} + \alpha, \subscript{\indexletter}{2} + \beta, \subscript{\indexletter}{3} - \tau}}}{\populationindex} + \superscript{\subscript{f}{\levelletter, \adaptiveroundbrackets{\subscript{\indexletter}{1} - \alpha, \subscript{\indexletter}{2} - \beta, \subscript{\indexletter}{3} + \tau}}}{\populationindex} + \superscript{\subscript{f}{\levelletter, \adaptiveroundbrackets{\subscript{\indexletter}{1} - \alpha, \subscript{\indexletter}{2} + \beta, \subscript{\indexletter}{3} - \tau}}}{\populationindex} \\
        &+ \superscript{\subscript{f}{\levelletter, \adaptiveroundbrackets{\subscript{\indexletter}{1} + \alpha, \subscript{\indexletter}{2} - \beta, \subscript{\indexletter}{3} - \tau}}}{\populationindex} - \superscript{\subscript{f}{\levelletter, \adaptiveroundbrackets{\subscript{\indexletter}{1} - \alpha, \subscript{\indexletter}{2} - \beta, \subscript{\indexletter}{3} - \tau}}}{\populationindex} \bigg ).
    \end{align*}\end{linenomath}
\end{itemize}
These operators satisfy the requirements we have highlighted so far in Definition \ref{def:PredictionOperator}, namely the locality and the consistency with the projection operator.
It should be observed that for $\predictionstencildepth > 0$, these operators are not convex combinations of the data, so we cannot expect any maximum principle to hold.

\subsection{Details and local regularity}

For the prediction operator is essentially an interpolation from coarse data obtained by averaging through  the projection operator, one understands that the more a function is locally regular, the more the predicted values will be close to the actual values.
On the other hand, if the function or its derivatives have abrupt spatial changes, the prediction operator will be less suitable to correctly reconstruct the averages.
This is precisely what is encoded in the details, which are a metric to quantify the information loss introduced by the projection operator and which carry essential information on the regularity of the solution. Thus, they let us identify areas of the computation domain where the spatial resolution can be reduced without affecting the quality of the reconstruction.
\begin{definition}[Details]
    Let $\populationindex = 0, \dots, \velocitiesnumber - 1$, $\levelletter = \minlevel + 1, \dots, \maxlevel$ and $\vectorial{\indexletter} \in \{0, \dots, N_{\levelletter} - 1\}^{\spatialdimension}$. The detail $\subscript{\superscript{d}{\populationindex}}{\levelletter, \vectorial{\indexletter}}$ for the $h$-th field on the cell $\subscript{\cellletter}{\levelletter, \vectorial{\indexletter}}$ is the difference between the actual value of the mean and the predicted value on this cell, that is
    \begin{linenomath}\begin{equation*}
        \subscript{\superscript{d}{\populationindex}}{\levelletter, \vectorial{\indexletter}} \definitionequality \superscript{\subscript{f}{\levelletter, \vectorial{\indexletter}}}{\populationindex} - \predicted{f}{\populationindex}{\levelletter, \vectorial{\indexletter}}.
    \end{equation*}\end{linenomath}
\end{definition}

Without entering into technicalities (see \cite{bellotti2021} and references therein), if the function has $\nu \geq 0 $ bounded derivatives in a suitable neighborhood $\tilde{\Sigma}_{\levelletter, \vectorial{\indexletter}}$ of the cell $\subscript{\cellletter}{\levelletter, \vectorial{\indexletter}}$ depending on the choice of prediction operator (namely on $\predictionstencildepth$), one can show that
\begin{linenomath}\begin{equation}\label{eq:DetailsDecayEstimate}
   |\subscript{\superscript{d}{\populationindex}}{\levelletter, \vectorial{\indexletter}}| \leqwithconstants \powertwo{-j\min{(\nu, \mu)}} |\superscript{f}{\populationindex}|_{W_{\infty}^{\min(\nu, \mu)}(\tilde{\Sigma}_{\levelletter, \vectorial{\indexletter}})}, \quad \text{with} \quad |f|_{W_{\infty}^{n}(\Omega)} \definitionequality \max_{|\vectorial{\alpha}| \leq n} \lVert D^{\vectorial{\alpha}} f \rVert_{L^{\infty}(\Omega)},
\end{equation}\end{linenomath}
being an appropriate Sobolev semi-norm computed on the neighborhood $\tilde{\Sigma}_{\levelletter, \vectorial{\indexletter}}$.
The possibility of establishing such control comes from the fact that the weights $(c_{\alpha})_{\alpha = 1}^{\alpha = \predictionstencildepth}$ have been chosen (see Appendix 2) so that the prediction operator is exact for polynomial functions of degree up to $2\predictionstencildepth$.
This estimate means that the details decrease with larger level $\levelletter$ if the regularity $\nu$ is such that $\nu > 0$ (namely if the solution is more than bounded), according to the smoothness of the function. Moreover, close to a jump discontinuity ($\nu = 0$) -- the typical situation with Riemann problems -- details have constant magnitude throughout the levels.
Finally, we emphasize that details on cells having the same parent are redundant as a consequence of the consistency property \eqref{eq:ConsistencyPrediction}. Hence, the following linear constraint holds\footnote{For $\spatialdimension = 1$, only one detail on two is significant; for $\spatialdimension = 2$, only three out of four; for $\spatialdimension = 3$, only seven out of eight.} for any $\populationindex = 0, \dots, \velocitiesnumber - 1$
\begin{linenomath}\begin{equation}\label{eq:DetailRedundancy}
    \sum_{\vectorial{\delta} \in \Sigma} \superscript{\subscript{d}{\levelletter + 1, 2 \vectorial{\indexletter} + \vectorial{\delta}}}{\populationindex} = 0, \qquad \levelletter = \minlevel, \dots, \maxlevel-1, \quad \vectorial{k} \in \{0, \dots, N_{\levelletter} - 1 \}^{\spatialdimension}.
\end{equation}\end{linenomath}

This calls for the introduction of the set of significant details $\setofindices_{\levelletter} \subset \{(\levelletter, \vectorial{\indexletter}) ~ : ~ \vectorial{k} \in \{0, \dots, N_{\levelletter} \}^{\spatialdimension} \}$ for $\levelletter = \minlevel + 1, \dots, \maxlevel$ where for each $2^{\spatialdimension}$ cell of level $\levelletter$ sharing the same parent, we eliminate one of them in order to avoid the redundancy by \eqref{eq:DetailRedundancy}.
We also set $\setofindices_{\minlevel} \definitionequality \{(\minlevel, \vectorial{\indexletter}) ~ : ~ \vectorial{\indexletter} \in \{0, \dots, N_{\minlevel} - 1 \}^{\spatialdimension} \}$.
We also introduce the set $\setofindices \definitionequality \cup_{\levelletter = \minlevel}^{\maxlevel} \setofindices_{\levelletter}$.
In this way, we have a one-to-one correspondence between data discretized on the finest level $\maxlevel$ and the data at the coarsest level $\minlevel$ plus the details at each level $\levelletter = \minlevel + 1, \dots, \maxlevel$ upon removing the redundancy\footnote{The interested reader can check that both sides of the equation have the same number of elements.}
\begin{linenomath}\begin{equation*}
    \vectorial{f}_{\maxlevel}^{\populationindex} \qquad \underset{\multiresolutiontransform^{-1}}{\overset{\multiresolutiontransform}{\myrightleftarrows{\rule{1cm}{0cm}}}} \qquad \adaptiveroundbrackets{\vectorial{f}_{\minlevel}^{\populationindex}, \vectorial{d}_{\minlevel + 1}^{\populationindex}, \dots,  \vectorial{d}_{\maxlevel}^{\populationindex}},
\end{equation*}\end{linenomath}
where $\vectorial{f}_{\levelletter}^{\populationindex} = (f_{\levelletter, \vectorial{\indexletter}}^{\populationindex})_{\vectorial{\indexletter} \in \{0, \dots, N_{\levelletter} - 1 \}^{\spatialdimension}}$ for $\levelletter = \minlevel, \dots, \maxlevel$ and $\vectorial{d}_{\levelletter}^{\populationindex} = (d_{\levelletter, \vectorial{\indexletter}})_{(\levelletter, \vectorial{\indexletter}) \in \setofindices_{\levelletter}}$ for $\levelletter = \minlevel + 1, \dots, \maxlevel$.

\subsection{Tree structure and grading}

\begin{definition}[Tree]
    A set of indices $\Lambda \subset \setofindices$ is said to represent a tree if the following holds:
    \begin{enumerate}
        \item $\setofindices_{\minlevel} \subset \Lambda$ 
        \item If $(\levelletter, \vectorial{\indexletter}) \in \Lambda$ and $\levelletter > \minlevel$, then also $2^{\spatialdimension} - 2$ siblings of $(\levelletter, \indexletter)$ belong to $\Lambda$.
        \item If $(\levelletter + 1, 2\vectorial{\indexletter} + \vectorial{\delta}) \in \Lambda$ for some $\vectorial{\delta} \in \Sigma$, then either its parent $(\levelletter, \vectorial{\indexletter}) \in \Lambda$ or its has $2^d - 1$ siblings belonging to $\Lambda$.
    \end{enumerate}
\end{definition}

Given a tree $\Lambda \subset \setofindices$ according to the previous definition, we indicate $R(\Lambda)$ (complete tree) the set of elements in $\Lambda$ (details cells) completed by their siblings missed due to the construction of $\setofindices_{\levelletter}$ (non-detail cells).
We also introduce the set of leaves $\leaves{\Lambda} \subset \Lambda$, which are the elements of $\Lambda$ without sons.
Adding the non-detail cells to $\leaves{\Lambda}$ we obtain the complete leaves $\physicalleaves{\Lambda}$.
As observed by \cite{cohen2003} $(C_{\levelletter, \vectorial{\indexletter}})_{(\levelletter, \vectorial{\indexletter}) \in \physicalleaves{\Lambda}}$ forms a hybrid partition of the domain $\domain$.

        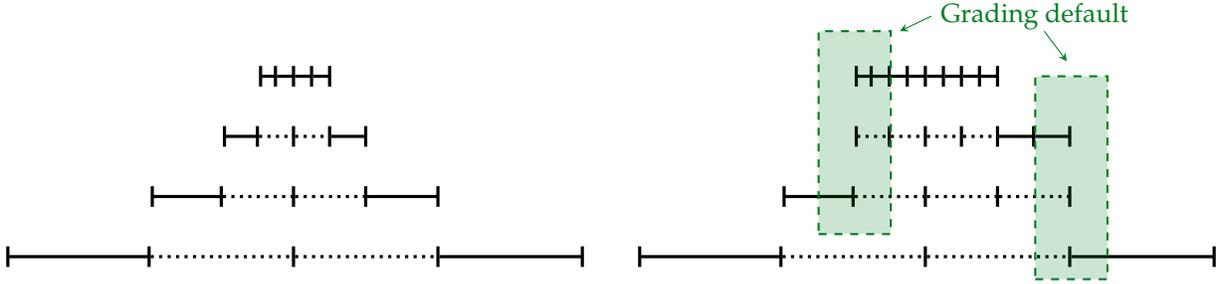
\begin{figure}[h]
        \begin{center}
        \begin{tikzpicture}[x=0.06cm, y=0.1cm]
    \newcount\levdis
    \levdis = 8
    
    \newcount\dist
    \dist = 140

    \draw[very thick, |-|]  (0, 0) -- (32, 0); 
    \draw[very thick,  -|, dotted] (32, 0) -- (64, 0); 
    \draw[very thick,  -|, dotted] (64, 0) -- (96, 0); 
    \draw[very thick,  -|] (96, 0) -- (128, 0); 

    \draw[very thick, |-|]  (32, \levdis) -- (48, \levdis); 
    \draw[very thick, -|, dotted]  (48, \levdis) -- (64, \levdis); 
    \draw[very thick, -|, dotted]  (64, \levdis) -- (80, \levdis); 
    \draw[very thick, -|]  (80, \levdis) -- (96, \levdis); 

    \draw[very thick, |-|]  (48, 2*\levdis) -- (56, 2*\levdis); 
    \draw[very thick, -|, dotted]  (56, 2*\levdis) -- (64, 2*\levdis); 
    \draw[very thick, -|, dotted]  (64, 2*\levdis) -- (72, 2*\levdis); 
    \draw[very thick, -|]  (72, 2*\levdis) -- (80, 2*\levdis); 

    \draw[very thick, |-|]  (56, 3*\levdis) -- (60, 3*\levdis); 
    \draw[very thick, -|]  (60, 3*\levdis) -- (64, 3*\levdis); 
    \draw[very thick, -|]  (64, 3*\levdis) -- (68, 3*\levdis); 
    \draw[very thick, -|]  (68, 3*\levdis) -- (72, 3*\levdis); 
    
    \draw[very thick, |-|]  (\dist, 0) -- (\dist + 32, 0); 
    \draw[very thick,  -|, dotted] (\dist + 32, 0) -- (\dist + 64, 0); 
    \draw[very thick,  -|, dotted] (\dist + 64, 0) -- (\dist + 96, 0); 
    \draw[very thick,  -|] (\dist + 96, 0) -- (\dist + 128, 0); 
    
    \draw[very thick, |-|]  (\dist + 32, \levdis) -- (\dist + 48, \levdis); 
    \draw[very thick, -|, dotted]  (\dist + 48, \levdis) -- (\dist + 64, \levdis); 
    \draw[very thick, -|, dotted]  (\dist + 64, \levdis) -- (\dist + 80, \levdis); 
    \draw[very thick, -|, dotted]  (\dist + 80, \levdis) -- (\dist + 96, \levdis); 

    \draw[very thick, |-|, dotted]  (\dist + 48, 2*\levdis) -- (\dist + 56, 2*\levdis); 
    \draw[very thick, -|, dotted]  (\dist + 56, 2*\levdis) -- (\dist + 64, 2*\levdis); 
    \draw[very thick, -|, dotted]  (\dist + 64, 2*\levdis) -- (\dist + 72, 2*\levdis); 
    \draw[very thick, -|, dotted]  (\dist + 72, 2*\levdis) -- (\dist + 80, 2*\levdis); 
    \draw[very thick, -|]  (\dist + 80, 2*\levdis) -- (\dist + 88, 2*\levdis); 
    \draw[very thick, -|]  (\dist + 88, 2*\levdis) -- (\dist + 96, 2*\levdis); 

    \draw[very thick, |-|]  (\dist + 48, 3*\levdis) -- (\dist + 52, 3*\levdis); 
    \draw[very thick, -|]  (\dist + 52, 3*\levdis) -- (\dist + 56, 3*\levdis); 
    \draw[very thick, -|]  (\dist + 56, 3*\levdis) -- (\dist + 60, 3*\levdis); 
    \draw[very thick, -|]  (\dist + 60, 3*\levdis) -- (\dist + 64, 3*\levdis); 
    \draw[very thick, -|]  (\dist + 64, 3*\levdis) -- (\dist + 68, 3*\levdis); 
    \draw[very thick, -|]  (\dist + 68, 3*\levdis) -- (\dist + 72, 3*\levdis); 
    \draw[very thick, -|]  (\dist + 72, 3*\levdis) -- (\dist + 76, 3*\levdis); 
    \draw[very thick, -|]  (\dist + 76, 3*\levdis) -- (\dist + 80, 3*\levdis); 
    
    \definecolor{bluej}{rgb}{0., 0.31764705882, 0.65490196078}
    \definecolor{greenjp1}{rgb}{0., 0.49803921568, 0.12549019607}

    \draw[greenjp1] (\dist + 65, 32) node [right] {Grading default};
    \draw[greenjp1, -stealth]  (\dist + 65, 32) --  (\dist + 58, 30);
    \draw[greenjp1, -stealth]  (\dist + 90, 30) --  (\dist + 95, 26);
    
    \draw[greenjp1, dashed, thick, fill = greenjp1, fill opacity = 0.2] (\dist+40,3) -- (\dist+16+40,3) -- (\dist+16+40,30) -- (\dist+40,30) -- (\dist+40,3);
    \draw[greenjp1, dashed, thick, fill = greenjp1, fill opacity = 0.2] (\dist+88,-3) -- (\dist+16+88,-3) -- (\dist+16+88,24) -- (\dist+88,24) -- (\dist+88,-3);

\end{tikzpicture}
        \end{center}\caption[]{\label{fig:gradedtree} Graded tree (left) \emph{versus} non-graded tree (right) in the case $d = 1$ for $\gamma = 1$. The complete leaves $\physicalleaves{\Lambda}$ are marked with full lines whereas the remaining cells $\physicaltree{\Lambda} \smallsetminus \physicalleaves{\Lambda}$ of the complete tree are rendered with dotted lines.}
        \end{figure}

Furthermore, Cohen \emph{et al.} \cite{cohen2003} (Proposition 2.3 in their work) have also shown that one can perform the multiresolution transform in an optimal way if the resulting tree structure is graded (see Figure \ref{fig:gradedtree}) with respect to the stencil of the prediction operator.
Still, observe that the lack of grading does not prevent one from performing the multiresolution analysis, because the notion of detail and the decay estimates are still available.
\begin{definition}[Grading]
    Let $\Lambda \subset \setofindices$ be a tree. Then it is said to be graded with respect to the prediction operator $\operatorial{P}_{\vartriangle}$ if for every $(\levelletter, \vectorial{\indexletter}) \in \physicaltree{\Lambda}$ with $\levelletter > \minlevel$ its prediction stencil also belongs to $\physicaltree{\Lambda}$.
\end{definition}
In the sequel, given a tree structure, the operation yielding the corresponding minimal graded tree shall be indicated by $\mathcal{G}$.
Therefore, grading is an important feature because it ensures an efficient implementation of the isomorphic transformation between averages on the complete leaves $\physicalleaves{\Lambda}$, namely $(f_{\levelletter, \vectorial{\indexletter}}^{\populationindex})_{(\levelletter, \indexletter) \in \physicalleaves{\Lambda}}$ and coarse means $(f^{\populationindex}_{\minlevel, \vectorial{\indexletter}})_{(\levelletter, \vectorial{\indexletter}) \in \setofindices_{\minlevel}}$ plus $(d^{\populationindex}_{\levelletter, \vectorial{\indexletter}})_{(\levelletter, \vectorial{\indexletter}) \in \setofindices_{\levelletter}}$ details on the tree for $\levelletter = \minlevel + 1, \dots, \maxlevel$.

\subsection{Mesh coarsening}

We have seen that knowing a function discretized on the whole finest level $\maxlevel$ is equivalent to know the means on the entire coarsest level $\minlevel$ and all the details -- upon eliminating their redundancy by \eqref{eq:DetailRedundancy} -- on each level $\levelletter = \minlevel + 1, \dots, \maxlevel$.
The passage from one representation to the other is performed by the so-called ``fast wavelet transform''.
Apart from this equivalence, the decomposition in terms of details is superior as far as we want to probe the local regularity of the functions as expressed by \eqref{eq:DetailsDecayEstimate}.
This can be exploited to coarsen the computational mesh in areas where the solution is strongly regular, still being sure that we can reconstruct information with accuracy within a certain given tolerance $0 < \thresholdletter \ll 1$.
This is done, theoretically, by setting to zero the details which are below a certain value.
From the practical point of view, one really eliminates the corresponding cells from the data structure at the end of the process.

Let $\Lambda \subset \setofindices$ be a graded tree and  $(f_{\levelletter, \vectorial{\indexletter}}^{\populationindex})_{(\levelletter, \indexletter) \in \physicalleaves{\Lambda}}$ the datum defined on its leaves.
In order to yield a tree structure with the following operations, we must treat detail cells having the same parent at once. This is done by considering the same detail for all of them\footnote{Other authors \cite{duarte2011} consider a $\ell^2$ metric under the form $\rangle d_{\levelletter + 1, 2\vectorial{\indexletter} + \vectorial{\delta}}^{\populationindex}  \langle_{2} \definitionequality 2^{-\spatialdimension/2} \sqrt{\sum_{\vectorial{\pi} \in \Sigma} |d_{\levelletter + 1, 2\vectorial{\indexletter} + \vectorial{\pi}}^{\populationindex} |^2}$ for $\delta \in \Sigma$.}
\begin{linenomath}\begin{equation*}
    \rangle d_{\levelletter + 1, 2\vectorial{\indexletter} + \vectorial{\delta}}^{\populationindex}  \langle_{\infty} \definitionequality \max_{\vectorial{\pi} \in \Sigma} |d_{\levelletter + 1, 2\vectorial{\indexletter} + \vectorial{\pi}}^{\populationindex} |, \qquad \delta \in \Sigma,
\end{equation*}\end{linenomath}
for $\levelletter = \minlevel, \dots, \maxlevel - 1$ and $\vectorial{\indexletter} \in \{0, \dots, N_{\levelletter}  - 1 \}^{\spatialdimension}$.
Thus we consider the coarsening (or threshold) operator $\mathcal{T}_{\epsilon}$ given by
\begin{linenomath}\begin{equation}\label{eq:Thresholding}
    \mathcal{T}_{\epsilon}(\Lambda) \definitionequality \setofindices_{\minlevel} \cup \adaptivecurlybrackets{(\levelletter, \vectorial{\indexletter}) \in \Lambda \smallsetminus \setofindices_{\minlevel} ~ : ~ \max_{\populationindex = 0, \dots, \velocitiesnumber - 1} \rangle d_{\levelletter, \vectorial{\indexletter}}^{\populationindex}  \langle_{\infty} ~ \geq \thresholdletter_{\levelletter}},
\end{equation}\end{linenomath}
where $\thresholdletter_{\minlevel + 1}, \dots, \thresholdletter_{\maxlevel}$ is a sequence of non-negative thresholds.
We also define the reconstruction operator $\reconstructed{~~}{}{}$, where $\reconstructed{\vectorial{f}}{\populationindex}{\maxlevel} = (\reconstructed{f}{\populationindex}{\maxlevel, \vectorial{\indexletter}})_{\vectorial{\indexletter} \in \{0, \dots, N_{\maxlevel} - 1 \}^{\spatialdimension}}$ representing the recursive application of the prediction operator $\operatorial{P}_{\vartriangle}$ until reaching available data stored on the complete leaves $\physicalleaves{\Lambda}$.
Following the proof given by Cohen \emph{et al.} \cite{cohen2003}, we prove that setting $\subscript{\thresholdletter}{\levelletter} \definitionequality \powertwo{\spatialdimension (\levelletter - \maxlevel)} \thresholdletter$ and considering the operator $\operatorial{T}_{\Lambda}$ setting to zero the details on cells which do not belong to $\Lambda$, then
\begin{linenomath}\begin{equation}\label{eq:ErrorControl}
    \norm{\reconstructed{\vectorial{f}}{\populationindex}{\maxlevel} - \multiresolutiontransform^{-1} \operatorial{T}_{\mathcal{G} \circ \mathcal{T}_{\epsilon}(\Lambda)} \multiresolutiontransform \reconstructed{\vectorial{f}}{\populationindex}{\maxlevel}}{\ell^p} \leq C_{\text{MR}} \epsilon,
\end{equation}
\end{linenomath}
with $p \in [1, \infty]$, where $\constantmultiresolution > 0$ depends on $\gamma$ and $p$, but not on the number of levels $\levelsnumber$.
In this section, we tried to succinctly present the basic elements without entering into technicalities.
The interested reader is referred to our previous contribution \cite{bellotti2021} and to \cite{cohen2003, muller2012, duarte2011}. 

\section{Adaptive \lb multiresolution scheme}\label{sec:AdaptiveLatticeBoltzmann}

In this central Section of the paper, we put all the previous ingredients together.
In particular, we shall answer the question on how to utilize the multiresolution procedure to handle time-dependent problems and on the way of adapting a given \lb scheme under the formalism of Section \ref{sec:LatticeBoltzmannOnVolumes}.

\subsection{Construction of the time adaptive mesh}\label{sec:TreeEnlargement}

At each discrete time $\timevariable$, we possess a solution defined on the complete leaves $\physicalleaves{\Lambda(t)}$ of the tree $\Lambda(t) \subset \setofindices$.
To compute the solution at the next time step $\timevariable + \Delta \timevariable$, we need to ensure that the computational lattice is refined enough so that an upper-bound similar to \eqref{eq:ErrorControl} still holds for the new solution.
Of course, due to the fact that at the moment of constructing the mesh, the new solution is still unknown, we have to devise a heuristics to slightly enlarge the tree $\Lambda(\timevariable)$ with the information known at time $t$.
The way of operating is resumed as follows

\begin{linenomath}\begin{equation*}
    \Lambda(\timevariable) \quad \xrightarrow[]{\mathcal{T}_{\thresholdletter}} \quad \mathcal{T}_{\thresholdletter} \adaptiveroundbrackets{\Lambda(\timevariable)} \quad  \xrightarrow[]{\mathcal{H}_{\thresholdletter}} \quad \mathcal{H}_{\thresholdletter} \circ \mathcal{T}_{\thresholdletter} \adaptiveroundbrackets{\Lambda(\timevariable)} \quad \xrightarrow[]{\mathcal{G}} \quad  \mathcal{G} \circ \mathcal{H}_{\thresholdletter} \circ \mathcal{T}_{\thresholdletter} \adaptiveroundbrackets{\Lambda(\timevariable)} =: \Lambda(\timevariable + \Delta \timevariable),
\end{equation*}\end{linenomath}
where $\mathcal{T}_{\thresholdletter}$ is the threshold operator, which eliminates superfluous cells according to \eqref{eq:Thresholding}, $\mathcal{H}_{\thresholdletter}$ is the enlarging operator that is defined in the sequel of the section and finally $\mathcal{G}$ is the grading operator introduced in the previous Section.
Again, we stress that in this phase details are computed with the old solution at time $\timevariable$.

\subsubsection{Propagation of information at finite speed}
On one hand, since we expect propagation of information at finite speed \emph{via} the stream phase \eqref{eq:StreamFinest} of the \lb method, we want to ensure that this flux of information is correctly captured by the computational mesh.
       Inspired by \cite{harten1994}, we thus request that $\mathcal{H}_{\thresholdletter}$ does the following
        \begin{linenomath}\begin{equation}
            \text{If} \quad \adaptiveroundbrackets{\levelletter, \vectorial{\indexletter}} \in \physicaltree{\mathcal{T}_{\thresholdletter} \adaptiveroundbrackets{\Lambda(\timevariable)}}, \quad \text{then} \quad \quad \adaptiveroundbrackets{\levelletter, \vectorial{\indexletter} - \superscript{\vectorial{\eta}}{\populationindex}} \in \physicaltree{\mathcal{H}_{\thresholdletter} \circ \mathcal{T}_{\thresholdletter} \adaptiveroundbrackets{\Lambda(\timevariable)}}, \qquad \populationindex = 0, \dots, \velocitiesnumber - 1. \label{eq:AdditionNeigh}
        \end{equation}\end{linenomath}
The previous formula \eqref{eq:AdditionNeigh} stipulates that if a cell at a certain level of refinement $\levelletter$ is kept in the structure, then we also keep some of its neighbors at the same level $\levelletter$. The number of kept neighbors is determined by the largest shift associated with the discrete velocities of the scheme.
We observe that this procedure is inherent to the hyperbolic equations where it has originally been developed \cite{harten1994, harten1995} and \cite{cohen2003} in the context of Finite Volume schemes.
For parabolic systems, where the propagation is done at infinite speed, this procedure cannot guarantee that we are able to capture all the phenomena \cite{nguessan2019}. 
Nevertheless, the infinite velocity is intrinsic to the continuous equations and one can still claim that the previous rule works pretty well in practice because the numerical scheme behaves ``hyperbolically''.

\subsubsection{Solution blowup}

On the other hand, as the equations we plan to solve are non-linear, and this non-linearity reflects on the collision phase \eqref{eq:CollisionOnFinest} of the \lb method, the regularity of the solution could vary (and decrease) in time. Therefore, especially with hyperbolic systems, the solutions can develop shocks which should be correctly captured by the mesh $\Lambda(t + \Delta t)$. Thus, we enforce that
        \begin{linenomath}\begin{gather*}
            \text{If} \quad \adaptiveroundbrackets{\levelletter, \vectorial{\indexletter}} \in \physicaltree{\mathcal{T}_{\thresholdletter} \adaptiveroundbrackets{\Lambda(\timevariable)}}, \quad \minlevel < \levelletter < \maxlevel \quad \text{and} \quad \max_{\populationindex = 0, \dots, \velocitiesnumber - 1} \rangle \subscript{\superscript{d}{\populationindex}}{\levelletter, \vectorial{\indexletter}}(t) \langle_{\infty} ~ \geq \powertwo{\regularityguess + d} \subscript{\thresholdletter}{\levelletter}, \\
            \text{then} \quad \adaptiveroundbrackets{\levelletter + 1, 2\vectorial{\indexletter} + \vectorial{\delta}} \in \physicaltree{\mathcal{H}_{\thresholdletter} \circ \mathcal{T}_{\thresholdletter} \adaptiveroundbrackets{\Lambda(\timevariable)}}, \qquad \vectorial{\delta} \in \Sigma,
        \end{gather*}\end{linenomath}
        where $\regularityguess \in \adaptivesquarebrackets{0, \predictionorder}$ has to be selected at the beginning of the simulation. Remark that this refinement criterion acts on $2^{\spatialdimension}$ siblings at once by refining all of them once the metric on the details is large enough. 
        This refinement criterion means that if the current cell is kept because its detail is not small enough to coarsen it, but moreover the detail is quite large, then we have to refine such cell. This allows, see \cite{bellotti2021} for more details, to identify areas of the mesh where the solution is undergoing a decrease of smoothness.
        The rationale is based on estimations like \eqref{eq:DetailsDecayEstimate} on details at time $t$ to estimate those (unavailable) at time $t + \Delta t$ and the interested reader can look at \cite{bellotti2021} for details.
        Since we do not know the number of bounded derivatives $\nu = \nu(t)$ of the solution  \emph{a priori}, we have to consider the following parameter $\regularityguess \definitionequality \min{(\nu, \predictionorder)}$ as an adjustable one.
        For example, if we know that the solution is going to develop some shock, it could be a wise choice to select $\regularityguess = 0$.

By proceeding at enlarging the computational mesh in this way by $\mathcal{H}_{\thresholdletter}$, we assume that it is suitable to represent the solution at the new time $t + \Delta t$ within a reasonable tolerance given by $\thresholdletter$.
This assumption is often called Harten's heuristics \cite{harten1994} in the world of Finite Volume.
After having applied the operators $\mathcal{T}_{\epsilon}$ and $\mathcal{H}_{\thresholdletter}$, the grid still has to be graded. Thus, applying $\mathcal{G}$, we obtain the graded mesh $\Lambda(t + \Delta t)$ to be used for the new computations.

Once the new mesh $\Lambda(t+\Delta t)$ is created, the old solution is adapted from $\physicalleaves{\Lambda(t)}$ to $\physicalleaves{\Lambda(t + \Delta t)}$ using the projection operator $\operatorial{P}_{\triangledown}$ when cells are merged by $\mathcal{T}_{\thresholdletter}$ and the prediction operator $\operatorial{P}_{\vartriangle}$ when they are created {either by} $\mathcal{H}_{\epsilon}$ or $\mathcal{G}$.
        
\subsection{Adapting the \lb method}

Once we have a new mesh $\Lambda(t + \Delta t)$ and the solution at time $t$ defined on $\physicalleaves{\Lambda(t+\Delta t)}$ at our disposal, we have to adapt the reference \lb method given by \eqref{eq:CollisionOnFinest} and \eqref{eq:StreamFinest} in order to be deployed only on the complete leaves of such a tree.
Given a cell $C_{j, \bm{k}}$ at level $j$, we consider the set $\mathcal{B}_{j, \bm{k}}$ of virtual cells at the finest level of refinement $\overline{J}$ covering $C_{j, \bm{k}}$ (recall that the grids are nested) and defined by
\begin{linenomath}\begin{equation}\label{eq:CellSplitting}
        \mathcal{B}_{\levelletter, \vectorial{\indexletter}} \definitionequality \adaptivecurlybrackets{\vectorial{\indexletter} \powertwo{\maxlevel - j}  + \vectorial{\delta} \quad : \quad \vectorial{\delta} \in \power{\adaptivecurlybrackets{0, \dots, \powertwo{\maxlevel - j} - 1}}{\spatialdimension}}.
\end{equation}\end{linenomath}

\subsubsection{Collision}\label{sec:Collision}

As seen before, the collision \eqref{eq:CollisionOnFinest} is local to each cell at the finest level $\maxlevel$, thus its generalization to non-uniform meshes is straightforward.
Indeed, consider any complete leaf $\subscript{\cellletter}{\levelletter, \vectorial{\indexletter}}$ with $(\levelletter, \vectorial{\indexletter}) \in \physicalleaves{\Lambda(\timevariable + \Delta \timevariable)}$, then we update with

    \begin{linenomath}\begin{align}
        \subscript{\vectorial{m}}{\levelletter, \vectorial{\indexletter}}(t) &= \operatorial{\changeofbasisletter} \subscript{\vectorial{f}}{\levelletter, \vectorial{\indexletter}} (t), \nonumber \\
        \subscript{\vectorial{f}}{\levelletter, \vectorial{\indexletter}}^{\collided} \adaptiveroundbrackets{\timevariable} &= \power{\operatorial{\changeofbasisletter}}{-1}  \adaptiveroundbrackets{(\operatorial{I} - \operatorial{S})\subscript{\vectorial{m}}{\levelletter, \vectorial{\indexletter}}(t) + \operatorial{S} \superscript{\vectorial{m}}{\text{eq}} \adaptiveroundbrackets{\subscript{\superscript{m}{0}}{\levelletter, \vectorial{\indexletter}}(t), \dots}}.\label{eq:AdaptiveCollision}
    \end{align}\end{linenomath}
As we fully discussed in \cite{bellotti2021}, one should be aware that this can potentially add an error to the solution because we are implicitly assuming that
\begin{equation}\label{eq:sourceoferror}
    \superscript{\vectorial{m}}{\text{eq}} \adaptiveroundbrackets{\subscript{\superscript{m}{0}}{\levelletter, \vectorial{\indexletter}}(t), \dots}  = \ratio{1}{\powertwo{\spatialdimension (\maxlevel - \levelletter)}} \sum_{\overline{\vectorial{k}} \in \subscript{\mathcal{B}}{\levelletter, \vectorial{\indexletter}}} \superscript{\vectorial{m}}{\text{eq}} \adaptiveroundbrackets{\reconstructed{m}{0}{\maxlevel, \overline{\vectorial{\indexletter}}}(t), \dots} ,
\end{equation}
which is true only if the equilibrium functions are linear. Otherwise, we can only hope to  have the equality plus an error of the order of $\thresholdletter$.
The previous formula means that colliding using the data directly available on the cell $C_{\levelletter, \vectorial{\indexletter}}$ is equivalent to doing it with the reconstructed data on the cells at the finest level $\mathcal{B}_{\levelletter, \vectorial{\indexletter}}$ covering $C_{\levelletter, \vectorial{\indexletter}}$ followed by an averaging step.
Performing the collision at the finest level $\maxlevel$ by using the reconstruction operator as on the right-hand side of \eqref{eq:AdaptiveCollision} and \eqref{eq:sourceoferror} (what \cite{hovhannisyan2010} call ``original adaptive scheme'' for the computation of source terms for Finite Volume scheme) clearly guarantees excellent results.
However, this would rely on the recursive nature of adaptive  multiresolution and would yield an explosion of the complexity of the algorithm in most cases where a high compression rate can be reached.
%
This holds even when memoization techniques are employed to reduce the number of evaluations due to the recursive structure of the reconstruction operator.
We verified both with 1D and 2D tests that for the problems we analyzed in \cite{bellotti2021} and in this paper, the use of collision operator on the mere complete leaves \eqref{eq:AdaptiveCollision} has a marginal impact on the accuracy of the adaptive method.

To sum up, performing the collision only on the complete leaves of the tree reduces the computational time associated with the collision phase because one performs a number of change of variables through $\operatorial{\changeofbasisletter}$ and evaluation of the non-linear equilibria equal to the cardinality of $\physicalleaves{t + \Delta t}$ instead of $2^{\spatialdimension \maxlevel}$. 

\subsubsection{Advection}

    \begin{figure}[h]
            \begin{center}
            \begin{tikzpicture}[x=0.11cm, y=0.11cm]

                \definecolor{bluej}{rgb}{0., 0.31764705882, 0.65490196078}
    \definecolor{greenjp1}{rgb}{0., 0.49803921568, 0.12549019607}
    
        \newcount\trans
    \trans = 26
    
\draw[-stealth, thick, densely dotted] (0, -5) -- (8, -5) node [below] {$x$};
\draw[-stealth, thick, densely dotted] (0, -5) -- (0, 8 - 5) node [left] {$y$};

\draw[-stealth, very thick] (0, -5) -- (4, 2 - 5) node [above, right] {$\superscript{\vectorial{\eta}}{\populationindex} = \adaptiveroundbrackets{2, 1}$};

\draw[very thick] (\trans, 0) -- (\trans + 16, 0) -- (\trans + 16, 16)--(\trans, 16)  -- cycle;

\draw[] (\trans + 8, -5) node [below] {$\adaptiveroundbrackets{\levelletter, \vectorial{\indexletter}}$};

\draw[very thick] (24 + \trans, 0)      -- (24 +  \trans + 16, 0) -- (24 +  \trans + 16, 16)--(24 +  \trans, 16)  -- cycle;
\draw[very thick] (24 + \trans + 2, 0)  -- (24 + \trans + 2, 16);
\draw[very thick] (24 + \trans + 4, 0)  -- (24 + \trans + 4, 16);
\draw[very thick] (24 + \trans + 6, 0)  -- (24 + \trans + 6, 16);
\draw[very thick] (24 + \trans + 8, 0)  -- (24 + \trans + 8, 16);
\draw[very thick] (24 + \trans + 10, 0) -- (24 + \trans + 10, 16);
\draw[very thick] (24 + \trans + 12, 0) -- (24 + \trans + 12, 16);
\draw[very thick] (24 + \trans + 14, 0) -- (24 + \trans + 14, 16);
\draw[very thick] (24 + \trans, 2)      -- (24 + \trans + 16, 2);
\draw[very thick] (24 + \trans, 4)      -- (24 + \trans + 16, 4);
\draw[very thick] (24 + \trans, 6)      -- (24 + \trans + 16, 6);
\draw[very thick] (24 + \trans, 8)      -- (24 + \trans + 16, 8);
\draw[very thick] (24 + \trans, 10)     -- (24 + \trans + 16, 10);
\draw[very thick] (24 + \trans, 12)     -- (24 + \trans + 16, 12);
\draw[very thick] (24 + \trans, 14)     -- (24 + \trans + 16, 14);

\draw[] (\trans + 8 + 24, -5) node [below] {$\mathcal{B}_{\levelletter, \vectorial{\indexletter}}$};

\draw[very thick, densely dotted, gray!50] (48 + 4 + \trans, 0)      -- (48 + 4 +  \trans + 16, 0) -- (48 + 4 +  \trans + 16, 16)--(48 + 4 +  \trans, 16)  -- cycle;
\draw[very thick, densely dotted, gray!50] (48 + 4 + \trans + 2, 0)  -- (48 + 4 + \trans + 2, 16);
\draw[very thick, densely dotted, gray!50] (48 + 4 + \trans + 4, 0)  -- (48 + 4 + \trans + 4, 16);
\draw[very thick, densely dotted, gray!50] (48 + 4 + \trans + 6, 0)  -- (48 + 4 + \trans + 6, 16);
\draw[very thick, densely dotted, gray!50] (48 + 4 + \trans + 8, 0)  -- (48 + 4 + \trans + 8, 16);
\draw[very thick, densely dotted, gray!50] (48 + 4 + \trans + 10, 0) -- (48 + 4 + \trans + 10, 16);
\draw[very thick, densely dotted, gray!50] (48 + 4 + \trans + 12, 0) -- (48 + 4 + \trans + 12, 16);
\draw[very thick, densely dotted, gray!50] (48 + 4 + \trans + 14, 0) -- (48 + 4 + \trans + 14, 16);
\draw[very thick, densely dotted, gray!50] (48 + 4 + \trans, 2)      -- (48 + 4 + \trans + 16, 2);
\draw[very thick, densely dotted, gray!50] (48 + 4 + \trans, 4)      -- (48 + 4 + \trans + 16, 4);
\draw[very thick, densely dotted, gray!50] (48 + 4 + \trans, 6)      -- (48 + 4 + \trans + 16, 6);
\draw[very thick, densely dotted, gray!50] (48 + 4 + \trans, 8)      -- (48 + 4 + \trans + 16, 8);
\draw[very thick, densely dotted, gray!50] (48 + 4 + \trans, 10)     -- (48 + 4 + \trans + 16, 10);
\draw[very thick, densely dotted, gray!50] (48 + 4 + \trans, 12)     -- (48 + 4 + \trans + 16, 12);
\draw[very thick, densely dotted, gray!50] (48 + 4 + \trans, 14)     -- (48 + 4 + \trans + 16, 14);

\draw[very thick] (48 -4  + 4 + \trans, 0 - 2)      -- (48 - 4  + 4 +  \trans + 16, 0 - 2 ) -- (48 - 4 + 4 +  \trans + 16, 16 - 2)--(48 - 4 + 4 +  \trans, 16 - 2)  -- cycle;
\draw[very thick] (48 - 4 + 4 + \trans + 2,  0  - 2)     -- (48  + 4- 4 + \trans + 2,  16 - 2);
\draw[very thick] (48 - 4 + 4 + \trans + 4,  0  - 2)     -- (48  + 4- 4 + \trans + 4,  16 - 2);
\draw[very thick] (48 - 4 + 4 + \trans + 6,  0  - 2)     -- (48  + 4- 4 + \trans + 6,  16 - 2);
\draw[very thick] (48 - 4 + 4 + \trans + 8,  0  - 2)     -- (48  + 4- 4 + \trans + 8,  16 - 2);
\draw[very thick] (48 - 4 + 4 + \trans + 10, 0  - 2)     -- (48  + 4- 4 + \trans + 10, 16 - 2);
\draw[very thick] (48 - 4 + 4 + \trans + 12, 0  - 2)     -- (48  + 4- 4 + \trans + 12, 16 - 2);
\draw[very thick] (48 - 4 + 4 + \trans + 14, 0  - 2)     -- (48  + 4- 4 + \trans + 14, 16 - 2);
\draw[very thick] (48 - 4 + 4 + \trans,      2  - 2)     -- (48  + 4- 4 + \trans + 16, 2  - 2);
\draw[very thick] (48 - 4 + 4 + \trans,      4  - 2)     -- (48  + 4- 4 + \trans + 16, 4  - 2);
\draw[very thick] (48 - 4 + 4 + \trans,      6  - 2)     -- (48  + 4- 4 + \trans + 16, 6  - 2);
\draw[very thick] (48 - 4 + 4 + \trans,      8  - 2)     -- (48  + 4- 4 + \trans + 16, 8  - 2);
\draw[very thick] (48 - 4 + 4 + \trans,      10 - 2)     -- (48  + 4- 4 + \trans + 16, 10 - 2);
\draw[very thick] (48 - 4 + 4 + \trans,      12 - 2)     -- (48  + 4- 4 + \trans + 16, 12 - 2);
\draw[very thick] (48 - 4 + 4 + \trans,      14 - 2)     -- (48  + 4- 4 + \trans + 16, 14 - 2);

\draw[] (\trans + 8 + 24 + 24 + 4, -5) node [below] {$\mathcal{B}_{\levelletter, \vectorial{\indexletter}} - \superscript{\vectorial{\eta}}{\populationindex}$};

\draw[very thick, bluej] (4 + 24 + 48 -4  + 4 + \trans, 0 - 2)      -- (4 + 48 + 24 - 4  + 4 +  \trans + 16, 0 - 2);
\draw[very thick, bluej] (4 + 24 + 48 -4  + 4 + \trans, 0 - 2) -- (4 + 48 + 24 - 4 + 4 +  \trans, 16 - 2);
\draw[very thick, bluej] (4 + 24 + 48 - 4 + 4 + \trans + 2,  0  - 2)     -- (4 + 48 + 24  + 4- 4 + \trans + 2,  16 - 2);
\draw[very thick, bluej] (4 + 24 + 48 - 4 + 4 + \trans + 4,  0  - 2)     -- (4 + 48 + 24  + 4- 4 + \trans + 4,  16 - 2);
\draw[very thick, bluej] (4 + 24 + 48 - 4 + 4 + \trans + 6,  0  - 2)     -- (4 + 48 + 24  + 4- 4 + \trans + 6,  2 - 2);
\draw[very thick, bluej] (4 + 24 + 48 - 4 + 4 + \trans + 8,  0  - 2)     -- (4 + 48 + 24  + 4- 4 + \trans + 8,  2 - 2);
\draw[very thick, bluej] (4 + 24 + 48 - 4 + 4 + \trans + 10, 0  - 2)     -- (4 + 48 + 24  + 4- 4 + \trans + 10, 2 - 2);
\draw[very thick, bluej] (4 + 24 + 48 - 4 + 4 + \trans + 12, 0  - 2)     -- (4 + 48 + 24  + 4- 4 + \trans + 12, 2 - 2);
\draw[very thick, bluej] (4 + 24 + 48 - 4 + 4 + \trans + 14, 0  - 2)     -- (4 + 48 + 24  + 4- 4 + \trans + 14, 2 - 2);
\draw[very thick, bluej] (4 + 24 + 48 - 4 + 4 + \trans + 16, 0  - 2)     -- (4 + 48 + 24  + 4- 4 + \trans + 16, 2 - 2);
\draw[very thick, bluej] (4 + 24 + 48 - 4 + 4 + \trans,      2  - 2)     -- (4 + 48 + 24  + 4- 4 + \trans + 16, 2  - 2);
\draw[very thick, bluej] (4 + 24 + 48 - 4 + 4 + \trans,      4  - 2)     -- (4 + 48 + 24  + 4- 4 + \trans + 4, 4  - 2);
\draw[very thick, bluej] (4 + 24 + 48 - 4 + 4 + \trans,      6  - 2)     -- (4 + 48 + 24  + 4- 4 + \trans + 4, 6  - 2);
\draw[very thick, bluej] (4 + 24 + 48 - 4 + 4 + \trans,      8  - 2)     -- (4 + 48 + 24  + 4- 4 + \trans + 4, 8  - 2);
\draw[very thick, bluej] (4 + 24 + 48 - 4 + 4 + \trans,      10 - 2)     -- (4 + 48 + 24  + 4- 4 + \trans + 4, 10 - 2);
\draw[very thick, bluej] (4 + 24 + 48 - 4 + 4 + \trans,      12 - 2)     -- (4 + 48 + 24  + 4- 4 + \trans + 4, 12 - 2);
\draw[very thick, bluej] (4 + 24 + 48 - 4 + 4 + \trans,      14 - 2)     -- (4 + 48 + 24  + 4- 4 + \trans + 4, 14 - 2);
\draw[very thick, bluej] (4 + 24 + 48 - 4 + 4 + \trans,      16 - 2)     -- (4 + 48 + 24  + 4- 4 + \trans + 4, 16 - 2);

\draw[very thick, greenjp1] (4 + 24 + 48 + 4 +  \trans + 16, 0) -- (4 + 24 + 48 + 4 +  \trans + 16, 16);
\draw[very thick, greenjp1] (4 + 24 + 48 + 4 +  \trans, 16) -- (4 + 24 + 48 + 4 +  \trans + 16, 16);
\draw[very thick, greenjp1] (4 + 24 + 48 + 4 + \trans + 0, 14)  -- (4 + 24 + 48 + 4 + \trans + 0, 16);
\draw[very thick, greenjp1] (4 + 24 + 48 + 4 + \trans + 2, 14)  -- (4 + 24 + 48 + 4 + \trans + 2, 16);
\draw[very thick, greenjp1] (4 + 24 + 48 + 4 + \trans + 4, 14)  -- (4 + 24 + 48 + 4 + \trans + 4, 16);
\draw[very thick, greenjp1] (4 + 24 + 48 + 4 + \trans + 6, 14)  -- (4 + 24 + 48 + 4 + \trans + 6, 16);
\draw[very thick, greenjp1] (4 + 24 + 48 + 4 + \trans + 8, 14)  -- (4 + 24 + 48 + 4 + \trans + 8, 16);
\draw[very thick, greenjp1] (4 + 24 + 48 + 4 + \trans + 10,14) -- (4 + 24 + 48 + 4 + \trans + 10, 16);
\draw[very thick, greenjp1] (4 + 24 + 48 + 4 + \trans + 12,0) -- (4 + 24 + 48 + 4 + \trans + 12, 16);
\draw[very thick, greenjp1] (4 + 24 + 48 + 4 + \trans + 14,0) -- (4 + 24 + 48 + 4 + \trans + 14, 16);
\draw[very thick, greenjp1] (4 + 24 + 48 + 4 + \trans + 12, 0)      -- (4 + 24 + 48 + 4 + \trans + 16, 0);
\draw[very thick, greenjp1] (4 + 24 + 48 + 4 + \trans + 12, 2)      -- (4 + 24 + 48 + 4 + \trans + 16, 2);
\draw[very thick, greenjp1] (4 + 24 + 48 + 4 + \trans + 12, 4)      -- (4 + 24 + 48 + 4 + \trans + 16, 4);
\draw[very thick, greenjp1] (4 + 24 + 48 + 4 + \trans + 12, 6)      -- (4 + 24 + 48 + 4 + \trans + 16, 6);
\draw[very thick, greenjp1] (4 + 24 + 48 + 4 + \trans + 12, 8)      -- (4 + 24 + 48 + 4 + \trans + 16, 8);
\draw[very thick, greenjp1] (4 + 24 + 48 + 4 + \trans + 12, 10)     -- (4 + 24 + 48 + 4 + \trans + 16, 10);
\draw[very thick, greenjp1] (4 + 24 + 48 + 4 + \trans + 12, 12)     -- (4 + 24 + 48 + 4 + \trans + 16, 12);
\draw[very thick, greenjp1] (4 + 24 + 48 + 4 + \trans, 14)     -- (4 + 24 + 48 + 4 + \trans + 16, 14);

\draw[bluej] (\trans + 8 + 24 + 24 + 24 + 8, -5) node [below] {$\adaptiveroundbrackets{\mathcal{B}_{\levelletter, \vectorial{\indexletter}} - \superscript{\vectorial{\eta}}{\populationindex}} \smallsetminus \mathcal{B}_{\levelletter, \vectorial{\indexletter}} $};

\draw[greenjp1] (\trans + 8 + 24 + 24 + 24 + 8, 24) node [below] {$ \mathcal{B}_{\levelletter, \vectorial{\indexletter}} \smallsetminus \adaptiveroundbrackets{\mathcal{B}_{\levelletter, \vectorial{\indexletter}} - \superscript{\vectorial{\eta}}{\populationindex}} $};

\end{tikzpicture}
\end{center}
\caption[]{\label{fig:adaptedstream} Example of the sets needed for the adaptive stream phase for the $d = 2$ case. We consider a leaf $\adaptiveroundbrackets{\levelletter, \vectorial{\indexletter}}$ which in this example is at level $\levelletter = \maxlevel - 3$ and a logical velocity $\superscript{\vectorial{\eta}}{\populationindex} = \adaptiveroundbrackets{2, 1}$ (just for illustrative purposes).}
\end{figure}

Again, consider a complete leaf $\subscript{\cellletter}{\levelletter, \vectorial{\indexletter}}$ with $(\levelletter, \vectorial{\indexletter
}) \in \physicalleaves{\Lambda(\timevariable + \Delta \timevariable)}$ and assume to work on the advection of the population $\populationindex = 0, \dots, \velocitiesnumber - 1$.
We again recover the indices $\mathcal{B}_{\levelletter, \vectorial{\indexletter}}$ of the virtual cells at the finest level of resolution $\maxlevel$ covering $C_{\levelletter, \vectorial{\indexletter}}$, see \eqref{eq:CellSplitting}.
    Let $\vectorial{\velocitynormalizedletter} \in \power{\relatives}{\spatialdimension}$, then $\mathcal{B}_{\levelletter, \vectorial{\indexletter}}  - \vectorial{\velocitynormalizedletter}$ indicates the element-wise subtraction of $\vectorial{\velocitynormalizedletter}$.
    Once we have this, the adaptive stream phase is given by
    \begin{equation}\label{eq:AdaptiveAdvectionPhase}
        \subscript{\superscript{f}{\populationindex}}{\levelletter, \vectorial{\indexletter}}\adaptiveroundbrackets{\timevariable + \Delta \timevariable} = \subscript{\superscript{f}{\populationindex, \collided}}{\levelletter, \vectorial{\indexletter}}\adaptiveroundbrackets{\timevariable } + \powertwo{\spatialdimension \adaptiveroundbrackets{\levelletter - \maxlevel}} \adaptiveroundbrackets{
        \sum_{\overline{\vectorial{\indexletter}}  \in \subscript{\superscript{\mathcal{E}}{\populationindex}}{\levelletter, \vectorial{\indexletter}}   } \reconstructed{f}{\populationindex, \collided}{\maxlevel, \overline{\vectorial{\indexletter}}}\adaptiveroundbrackets{\timevariable} - \sum_{\overline{\vectorial{\indexletter}} \in \subscript{\superscript{\mathcal{A}}{\populationindex}}{\levelletter, \vectorial{\indexletter}}} \reconstructed{f}{\populationindex, \collided}{\maxlevel, \overline{\vectorial{\indexletter}}}  \adaptiveroundbrackets{\timevariable}},
    \end{equation}
    with (Figure \ref{fig:adaptedstream})
    \begin{linenomath}\begin{equation*}
        \subscript{\superscript{\mathcal{E}}{\populationindex}}{\levelletter, \vectorial{\indexletter}} \definitionequality  \adaptiveroundbrackets{\mathcal{B}_{\levelletter, \vectorial{\indexletter}} - \superscript{\vectorial{\eta}}{\populationindex}} \smallsetminus \mathcal{B}_{\levelletter, \vectorial{\indexletter}}, \qquad \subscript{\superscript{\mathcal{A}}{\populationindex}}{\levelletter, \vectorial{\indexletter}} \definitionequality \mathcal{B}_{\levelletter, \vectorial{\indexletter}} \smallsetminus \adaptiveroundbrackets{\mathcal{B}_{\levelletter, \vectorial{\indexletter}} - \superscript{\vectorial{\eta}}{\populationindex}}.
    \end{equation*}\end{linenomath}
The updated average on cell $\vectorial{\indexletter}$ at level $\levelletter$ is obtained by the post-collisional average on the same cell completed by the reconstruction of the post-collisional averages on some shifted virtual cells at the finest level of refinement $\maxlevel$, according to the velocity $\populationindex$ we are considering.
Here the first summed term corresponds to an incoming pseudo-flux, whereas the second one is the outgoing one. We mimic the way of proceeding by \cite{cohen2003} performing the so-called ``exact flux reconstruction'' for Finite Volume schemes via the double hat operator.
Conversely, the AMR approach from Fakhari \emph{et al.} \cite{fakhari2014, fakhari2016} relies on a Lax-Wendroff scheme with direct evaluation on some ghost (or halo) cells at the same level of resolution $\levelletter$. This procedure is surely cheaper than \eqref{eq:AdaptiveAdvectionPhase} but cannot ensure any control on the error because the mesh is generated by a heuristic criterion which is not correlated to the way of locally reconstructing the solution \emph{via} a cascade of interpolations.

It is important to observe that once we consider a cell $C_{\levelletter, \vectorial{\indexletter}}$ adjacent to the boundary of the domain $\partial \Omega$, some virtual cell indexed by $\mathcal{E}_{\levelletter, \vectorial{\indexletter}}^{\populationindex}$ could lie outside the domain. This calls for the enforcement of some boundary condition as explained in the following Section \ref{sec:BoundaryConditions}.

It can be easily shown that our approach, devised through a different way of reasoning, is indeed a CTU (Corner Transport Upwind method) Finite Volume discretization by Colella \cite{colella1990multidimensional} (also see Leveque \cite{leveque1992}) of the advection phase, with reconstruction of a piece-wise constant representation of the solution at the finest level made possible by the multiresolution.
At finest level, the advection phase is given by the translation of data in the opposite direction of the velocity field.
Let us illustrate the construction of \eqref{eq:AdaptiveAdvectionPhase}: following \eqref{eq:StreamFinest}
    \begin{linenomath}\begin{equation*}
        \subscript{\superscript{f}{\populationindex}}{\maxlevel, \overline{\vectorial{\indexletter}}} \adaptiveroundbrackets{\timevariable + \Delta \timevariable} = \subscript{\superscript{f}{\populationindex, \collided}}{\maxlevel, \overline{\vectorial{\indexletter}} - \superscript{\vectorial{\eta}}{\populationindex}} \adaptiveroundbrackets{\timevariable }, \qquad \overline{\vectorial{\indexletter}} \in \mathcal{B}_{\levelletter, \vectorial{\indexletter}}.
    \end{equation*}\end{linenomath}
    Unfortunately, the right-hand side of these equations is not always well defined because the cell could fail to be a complete leaf of the tree.
    Nevertheless, we are allowed to reconstruct data on these cells using multiresolution with a precision controlled by $\thresholdletter$
    \begin{linenomath}\begin{equation*}
        \subscript{\superscript{f}{\populationindex}}{\maxlevel, \overline{\vectorial{\indexletter}}} \adaptiveroundbrackets{\timevariable + \Delta \timevariable} \simeq \reconstructed{f}{\populationindex, \collided}{\maxlevel, \overline{\vectorial{\indexletter}} - \superscript{\vectorial{\eta}}{\populationindex}} \adaptiveroundbrackets{\timevariable }, \qquad \overline{\vectorial{\indexletter}} \in \mathcal{B}_{\levelletter, \vectorial{\indexletter}}.
    \end{equation*}\end{linenomath}
    We come back to level $\levelletter$ by $\maxlevel - j$ applications of the projection operator $\operatorial{P}_{\triangledown}$ and consistency, yielding
    \begin{linenomath}\begin{align*}
        \subscript{\superscript{f}{\populationindex}}{\levelletter, \vectorial{\indexletter}} \adaptiveroundbrackets{\timevariable + \Delta \timevariable} &\simeq \powertwo{\spatialdimension(\levelletter - \maxlevel)} \sum_{\overline{\vectorial{\indexletter}} \in \mathcal{B}_{\levelletter, \vectorial{\indexletter}}} \reconstructed{f}{\populationindex, \collided}{\maxlevel, \overline{\vectorial{\indexletter}} - \superscript{\vectorial{\eta}}{\populationindex}} \adaptiveroundbrackets{\timevariable} \\
        &= \powertwo{\spatialdimension(\levelletter - \maxlevel)} \adaptiveroundbrackets{ \sum_{\overline{\vectorial{\indexletter}} \in \mathcal{B}_{\levelletter, \vectorial{\indexletter}}} \reconstructed{f}{\populationindex, \collided}{\maxlevel, \overline{\vectorial{\indexletter}}}\adaptiveroundbrackets{\timevariable } + \sum_{\overline{\vectorial{\indexletter}} \in \mathcal{E}^{\populationindex}_{\levelletter, \vectorial{\indexletter}}} \reconstructed{f}{\populationindex, \collided}{\maxlevel, \overline{\vectorial{\indexletter}}} \adaptiveroundbrackets{\timevariable}  - \sum_{\overline{\vectorial{\indexletter}} \in \mathcal{A}^{\populationindex}_{\levelletter, \vectorial{\indexletter}}} \reconstructed{f}{\populationindex, \collided}{\maxlevel, \overline{\vectorial{\indexletter}}}  \adaptiveroundbrackets{\timevariable}} \\
        &=\subscript{\superscript{f}{\populationindex, \collided}}{\levelletter, \vectorial{\indexletter}} \adaptiveroundbrackets{\timevariable } + \powertwo{\spatialdimension(\levelletter - \maxlevel)} \adaptiveroundbrackets{\sum_{\overline{\vectorial{\indexletter}} \in \mathcal{E}^{\populationindex}_{\levelletter, \vectorial{\indexletter}}} \reconstructed{f}{\populationindex, \collided}{\maxlevel, \overline{\vectorial{\indexletter}}}  \adaptiveroundbrackets{\timevariable }  - \sum_{\overline{\vectorial{\indexletter}} \in \mathcal{A}^{\populationindex}_{\levelletter, \vectorial{\indexletter}}} \reconstructed{f}{\populationindex,\collided}{\maxlevel, \overline{\vectorial{\indexletter}}}  \adaptiveroundbrackets{\timevariable} }.
    \end{align*}\end{linenomath}
   
If the implementation of the reconstruction operator is optimized, one can expect a gain in computational time when using \eqref{eq:AdaptiveAdvectionPhase} on the adaptive grid because, as one can see from Figure \ref{fig:adaptedstream}, the sets $ \mathcal{E}^{\populationindex}_{\levelletter, \vectorial{\indexletter}}$ and $ \mathcal{A}^{\populationindex}_{\levelletter, \vectorial{\indexletter}}$ only involve few cells at the finest level $\maxlevel$ close to the edges of each leaf.

\subsection{Boundary conditions}\label{sec:BoundaryConditions}

We are left to describe our implementation of suitable boundary conditions for the problem.
We consider a cell $C_{\levelletter, \vectorial{\indexletter}}$ touching the boundary $\partial \domain$, so that $|\overline{C_{\levelletter, \vectorial{\indexletter}}} \cap \partial \domain|_{\spatialdimension-1} > 0$ and a velocity field $\superscript{\vectorial{\eta}}{\populationindex}$ with information coming from outside the domain, so $\superscript{\vectorial{\eta}}{\populationindex} \cdot \vectorial{n} < 0$, where $\vectorial{n}$ is the normal vector to $\partial \domain$.
For illustrative purposes, we consider three types of boundary conditions which -- for the reference scheme -- correspond to
\begin{equation}\label{eq:refBoundaryConditions}
    f^{\populationindex}(t + \Delta t) = 
    \begin{cases}
        f^{\populationindex, \collided}(t), \qquad &\text{(Copy)}, \\
        f^{\populationindex^{\dag}, \collided}(t), \qquad &\text{(Bounce Back)}, \\
        -f^{\populationindex^{\dag}, \collided}(t), \qquad &\text{(Anti Bounce Back)},
    \end{cases}
\end{equation}
where we recall that $\vectorial{\eta}^{\populationindex^{\dag}} = - \vectorial{\eta}^{\populationindex}$ yielding the opposite dimensionless velocity. This given the name ``bounce back'' to this kind of conditions.
One can handle non-homogeneous conditions or more intricate ones in the same way.
Observe that, in the advection phase \eqref{eq:AdaptiveAdvectionPhase}, we cannot close the system due to the lack of some pieces of information, indicated without the reconstruction sign
    \begin{linenomath}\begin{align*}
        \subscript{\superscript{f}{\populationindex}}{\levelletter, \vectorial{\indexletter}} \adaptiveroundbrackets{\timevariable + \Delta \timevariable} =\subscript{\superscript{f}{\populationindex, \collided}}{\levelletter, \vectorial{\indexletter}} \adaptiveroundbrackets{\timevariable } &+ \powertwo{\spatialdimension(\levelletter - \maxlevel)} \adaptiveroundbrackets{\sum_{\overline{\vectorial{\indexletter}} \in \mathcal{E}^{\populationindex}_{\levelletter, \vectorial{\indexletter}} \cap \{ 0, \dots, N_{\maxlevel} - 1 \}^{\spatialdimension}} \reconstructed{f}{\populationindex, \collided}{\maxlevel, \overline{\vectorial{\indexletter}}} \adaptiveroundbrackets{\timevariable }  - \sum_{\overline{\vectorial{\indexletter}} \in \mathcal{A}^{\populationindex}_{\levelletter, \vectorial{\indexletter}}} \reconstructed{f}{\populationindex, \collided}{\maxlevel, \overline{\vectorial{\indexletter}}} \adaptiveroundbrackets{\timevariable} } \\
        &+ \powertwo{\spatialdimension(\levelletter - \maxlevel)} \sum_{\overline{\vectorial{\indexletter}} \in \mathcal{E}^{\populationindex}_{\levelletter, \vectorial{\indexletter}} \smallsetminus \{ 0, \dots, N_{\maxlevel} - 1 \}^{\spatialdimension}} \subscript{\superscript{f}{\populationindex, \collided}}{\maxlevel, \overline{\vectorial{\indexletter}}} \adaptiveroundbrackets{\timevariable } .
    \end{align*}\end{linenomath}
Having $\mathcal{A}^{\populationindex}_{\levelletter, \vectorial{\indexletter}} \subset \mathcal{B}^{\populationindex}_{\levelletter, \vectorial{\indexletter}} \subset \{ 0, \dots, N_{\maxlevel} - 1 \}^{\spatialdimension}$, all the quantities appearing in the outgoing pseudo-fluxes are known since they come from inside the domain.
Using \eqref{eq:refBoundaryConditions} and presenting the ``bounce back'' case only for the sake of presentation, this corresponds to
    \begin{linenomath}\begin{align}
        \subscript{\superscript{f}{\populationindex}}{\levelletter, \vectorial{\indexletter}} \adaptiveroundbrackets{\timevariable + \Delta \timevariable} =\subscript{\superscript{f}{\populationindex, \collided}}{\levelletter, \vectorial{\indexletter}} \adaptiveroundbrackets{\timevariable} &+ \powertwo{\spatialdimension(\levelletter - \maxlevel)} \adaptiveroundbrackets{\sum_{\overline{\vectorial{\indexletter}} \in \mathcal{E}^{\populationindex}_{\levelletter, \vectorial{\indexletter}} \cap \{ 0, \dots, N_{\maxlevel} - 1 \}^{\spatialdimension}} \reconstructed{f}{\populationindex, \collided}{\maxlevel, \overline{\vectorial{\indexletter}}}  \adaptiveroundbrackets{\timevariable}  - \sum_{\overline{\vectorial{\indexletter}} \in \mathcal{A}^{\populationindex}_{\levelletter, \vectorial{\indexletter}}} \reconstructed{f}{\populationindex, \collided}{\maxlevel, \overline{\vectorial{\indexletter}}}  \adaptiveroundbrackets{\timevariable } } \nonumber \\
        &+ \powertwo{\spatialdimension(\levelletter - \maxlevel)}
        \sum_{\overline{\vectorial{\indexletter}} \in \adaptiveroundbrackets{ \mathcal{E}^{\populationindex}_{\levelletter, \vectorial{\indexletter}} \smallsetminus \{ 0, \dots, N_{\maxlevel} - 1 \}^{\spatialdimension}} + \superscript{\vectorial{\eta}}{\populationindex}}\subscript{\superscript{f}{\populationindex^{\dag}, \collided}}{\maxlevel, \overline{\vectorial{\indexletter}}} \adaptiveroundbrackets{\timevariable }, \label{eq:AdvectionOnBoundary}
    \end{align}\end{linenomath}
where all the sums pertain to quantities inside $\domain$, so the quantities on the second line of \eqref{eq:AdvectionOnBoundary} can be evaluated either using the usual reconstruction operator or by direct evaluation, a cheaper alternative. 
We implement the latter method, which consists in taking the value directly available on the corresponding leaf.

\section{Implementation}\label{sec:ImplementationAndOptimization}

As already emphasized in the course of the paper, the structure of the multiresolution is intrinsically recursive, leading to important slow-downs when implementing the reconstruction operator in this fashion to obtain the solution on the finest mesh.
This is especially evident for large problems with $d = 2, 3$ (see \cite{cohen2003, muller2012}).
In this Section, we provide an interesting idea on how to speed-up the computation of the reconstruction operator for the sole purpose of the \lb scheme in order to guarantee the proper accurate evaluation of the advection without relying on the full reconstruction at the finest level. Besides, we describe the overall structure of the numerical solver.

\subsection{A non-recursive implementation of the reconstruction operator}

To make computations feasible for large problems, we follow the idea of \cite{cohen2003}, claiming that in the univariate case, the recursive application of a linear prediction operator can be condensed into the computation of the powers of a given matrix at the beginning of the simulation, based on the assumption that the fluxes of the Finite Volume method involve only adjacent cells. For $\spatialdimension = 1$ and $\predictionstencildepth = 1$, this reads
\begin{equation}\label{eq:PredictionWithMatrix}
    \adaptiveroundbrackets{\begin{matrix*}[l]
                                \reconstructed{f}{\populationindex, {\collided}}{\maxlevel, \powertwo{\maxlevel - \levelletter}\indexletter - 2} \\
                                \reconstructed{f}{\populationindex, {\collided}}{\maxlevel, \powertwo{\maxlevel - \levelletter}\indexletter - 1} \\
                                \reconstructed{f}{\populationindex, {\collided}}{\maxlevel, \powertwo{\maxlevel - \levelletter}\indexletter} \\
                                \reconstructed{f}{\populationindex, {\collided}}{\maxlevel, \powertwo{\maxlevel - \levelletter}\indexletter + 1}
                           \end{matrix*}} = 
    \adaptiveroundbrackets{\begin{matrix}
                                \ratioobl{1}{8} & 1 & -\ratioobl{1}{8} & 0 \\
                                -\ratioobl{1}{8} & 1 & \ratioobl{1}{8} & 0 \\
                                0 & \ratioobl{1}{8} & 1 & -\ratioobl{1}{8} \\
                                0 & -\ratioobl{1}{8} & 1 & \ratioobl{1}{8} 
                           \end{matrix}}^{\maxlevel - \levelletter}
    \adaptiveroundbrackets{\begin{matrix*}[l]
                                f^{\populationindex, {\collided}}_{\levelletter, \indexletter - 2} \\
                                f^{\populationindex, {\collided}}_{\levelletter, \indexletter - 1} \\
                                f^{\populationindex, {\collided}}_{\levelletter, \indexletter} \\
                                f^{\populationindex, {\collided}}_{\levelletter, \indexletter + 1}
                           \end{matrix*}}.
\end{equation}
Clearly, this method works far from the boundary $\partial \domain$ and in areas where the local refinement level of the leaves is constant. Since for $\spatialdimension > 1$ we have constructed the prediction operator by tensor product, the matrices involved in this framework shall just be the Kronecker product of that for $\spatialdimension = 1$ in \eqref{eq:PredictionWithMatrix} $\spatialdimension$-times with itself.

        \begin{figure}
        \begin{center}
            \def\svgwidth{0.92\textwidth}
            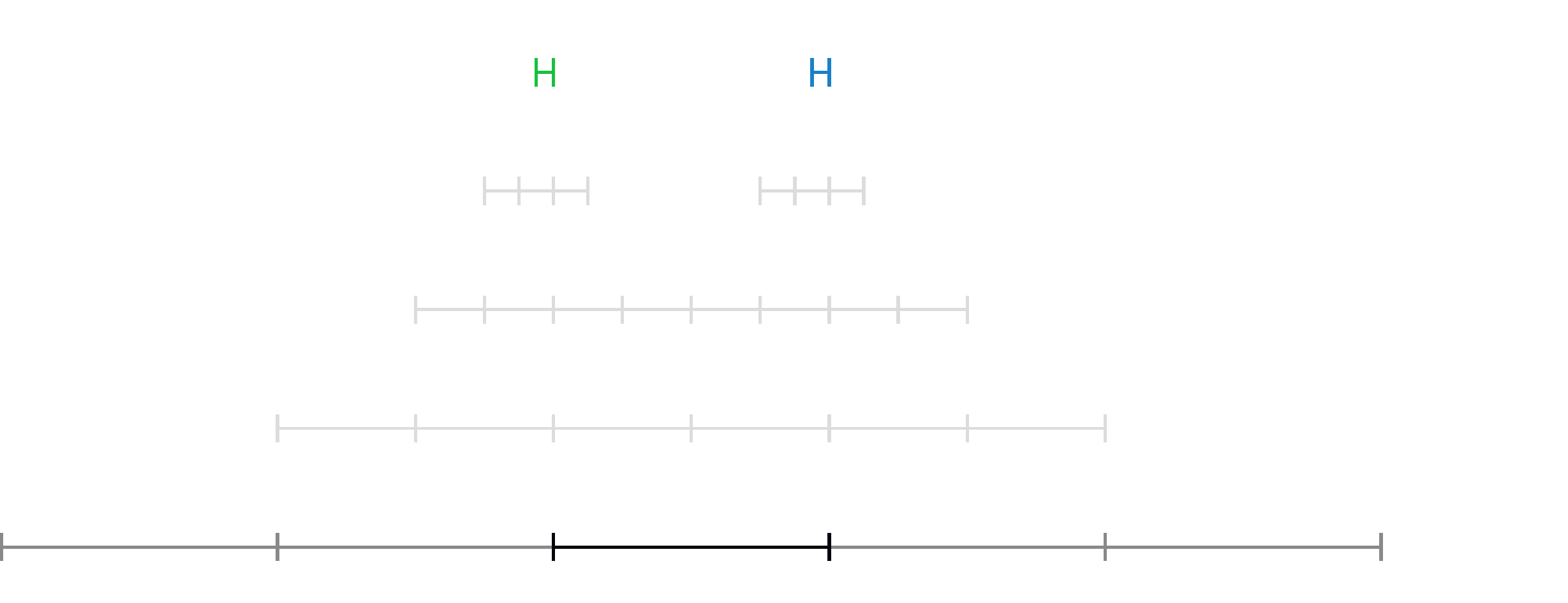
        \end{center}
        \caption{\label{fig:flattening}Example of non-recursive implementation of the reconstruction operator for $\spatialdimension = 1$ and $\predictionstencildepth = 1$ for the velocity $\latticevelocity$. The cell on which one updates the solution is four level far from the finest level and indexed by $\indexletter$. The green cell on top corresponds to the cell at the finest level $\maxlevel$ giving the incoming pseudo-flux for the cell $\indexletter$ at level $\levelletter$, whereas that in blue corresponds to that of the outgoing one. The prediction operator is progressively applied spanning the intermediate (non existing) cells in pale grey inside the green (resp. blue) funnel, until reaching cells at the same level $\levelletter$ in black and grey.}
        \end{figure}

We start to describe the procedure for larger spatial dimensions $\spatialdimension$ by selecting a complete leaf $\subscript{\cellletter}{\levelletter, \vectorial{\indexletter}}$ and assuming that it is surrounded by enough leaves of the same level.
At the beginning of the numerical simulation, we can once for all compute by recursion (we also did it analytically), for any level $\levelletter$ and velocity $\populationindex$, the set of logical shifts $\superscript{\subscript{\Xi}{\levelletter}}{\populationindex} \subset \power{\relatives}{\spatialdimension}$ and real weights $\subscript{(\superscript{\subscript{w}{\levelletter, \vectorial{\delta}}}{\populationindex})}{\vectorial{\delta} \in \superscript{\subscript{\Xi}{\levelletter}}{\populationindex}}$ so that
\begin{equation}\label{eq:PredictionWithMatrixNew}
        \sum_{\overline{\vectorial{\indexletter}} \in \mathcal{E}_{\levelletter, \vectorial{\indexletter}}^{\populationindex}} \reconstructed{f}{\populationindex, {\collided}}{\maxlevel, \overline{\vectorial{\indexletter}}} \adaptiveroundbrackets{\timevariable }  = \sum_{\vectorial{\delta} \in \superscript{\subscript{\Xi}{\levelletter}}{\populationindex}} \superscript{\subscript{w}{\levelletter, \vectorial{\delta}}}{\populationindex} ~ \subscript{\superscript{f}{\populationindex, {\collided}}}{\levelletter, \vectorial{\indexletter} + \vectorial{\delta}} \adaptiveroundbrackets{\timevariable }, 
\end{equation}
and the same for $\mathcal{A}_{\levelletter, \vectorial{\indexletter}}^{\populationindex}$, transforming the recursive left-hand side of \eqref{eq:PredictionWithMatrixNew} into a right-hand side made up of linear combinations of data known on the leaves with previously computed weights, in analogy with \eqref{eq:PredictionWithMatrix}. 
An illustration of such a process in a one-dimensional setting is given in Figure \ref{fig:flattening}.

On the other hand, if the surrounding leaves at the same level are not enough (we could, for example, fall on some ghost or halo cell\footnote{We added $2\predictionstencildepth$ of them in each Cartesian and diagonal direction to ensure the feasibility of any operation concerning the prediction operator.}), we are not sure that the value we retrieve is accurate enough according to the multiresolution analysis.
Let us study this more in detail, considering $ \vectorial{\delta} \in \superscript{\subscript{\Xi}{\levelletter}}{\populationindex}$. There are two possibilities:
\begin{itemize}
    \item $\subscript{\cellletter}{\levelletter, \vectorial{\indexletter} + \vectorial{\delta}}$ is a complete leaf, thus everything is fine, because multiresolution allows us to employ this value in a recursion formula to reconstruct at finest level without adding any detail.
    \item $\subscript{\cellletter}{\levelletter, \vectorial{\indexletter} + \vectorial{\delta}}$ is not a complete leaf. There are two situations (and not more thanks to the construction of the mesh), one of which requires particular care:
    \begin{itemize}
        \item It intersects a leaf at the coarser level $\levelletter - 1$. The value we retrieve is fine because multiresolution guarantees that quantities computed by applying the prediction operator without adding details are fine if their respective cell is situated on top of a leaf.
        \item It intersects a leaf at the finer level $\levelletter + 1$. This is the critical situation, because the retrieved value is not accurate enough to be employed in reconstructions, according to the multiresolution analysis and we have to add the proper detail information in order to preserve our target of error control.
    \end{itemize}
\end{itemize}

Comparing for $\spatialdimension = 2$, $\predictionstencildepth = 1$ and a D2Q9 scheme to the Lax-Wendroff advection phase by \cite{fakhari2014, fakhari2015, fakhari2016}, which has to recover three values of the solution at the considered level $\levelletter$ in each velocity direction, multiply each of them by a weight and then summing them, it can be shown that our method using \eqref{eq:PredictionWithMatrixNew} needs to recover at most 25 values -- whatever the velocity satisfying $\max(|\eta^{\populationindex}_x|, |\eta^{\populationindex}_y|) \leq 2$, thus for a large family of schemes -- multiply each of them by a weight and then summing them.
Even if this should be quantified precisely in terms of computational cost for a given problem, a given implementation and a given architecture,
the complexity of the algorithms are sensibly at the same level.
The gain comes once considering that we can deal with a very large class of schemes achieving a control on the error which is not obtained by the AMR procedures and that current investigations using formal expansions indicate that our adaptive method reproduces the behavior of the reference method on the finest grid in terms of both error control and equivalent equations up to order three, instead of two for the Lax-Wendroff approach.

\subsection{Overall solver}

To summarize, we present the overall procedure in a synthetic fashion.
In the actual algorithm, we swap the collision and the advection phase as in Algorithm \ref{alg:OverallSolver}.
This does not change the result since we start with initial data at the equilibrium.

\begin{algorithm}[H]
\SetAlgoLined
 Initial time : $t = 0$\;
 Full finest mesh : $\Lambda(0) = \{ (\maxlevel, \vectorial{\indexletter}) \quad : \quad \vectorial{\indexletter} \in \{0, \dots, N_{\maxlevel}-1 \}^{\spatialdimension} \}$\;
 Initial datum at equilibrium : $\subscript{\vectorial{f}}{\maxlevel, \vectorial{\indexletter}}(0) = \operatorial{M}^{-1} \subscript{\superscript{\vectorial{m}}{\text{eq}}}{\maxlevel, \vectorial{\indexletter}}(0)$\;
 \While{$t \leq T$}{
  Mesh adaptation : $\Lambda(t+\Delta t) = \mathcal{G} \circ \subscript{\mathcal{H}}{\thresholdletter} \circ \subscript{\mathcal{T}}{\thresholdletter} (\Lambda(t))$ computing the details using $\subscript{\operatorial{P}}{\triangledown} $ and $ \subscript{\operatorial{P}}{\vartriangle}$\;
  The data structure on $\physicalleaves{\Lambda(t)}$ is updated on $\physicalleaves{\Lambda(t+\Delta t)}$\;
  \For{$(\levelletter, \vectorial{\indexletter}) \in \physicalleaves{\Lambda(t+\Delta t)}$}{
    \eIf{$(\levelletter, \vectorial{\indexletter})$ internal to $\domain$}{
        Advect using \eqref{eq:AdaptiveAdvectionPhase} and optimizing using \eqref{eq:PredictionWithMatrixNew} when possible \;
    }{
        Advect using \eqref{eq:AdvectionOnBoundary} and optimizing using \eqref{eq:PredictionWithMatrixNew} when possible\;    
    }
    Collide using \eqref{eq:AdaptiveCollision}\;
    }
 }
 \caption{\label{alg:OverallSolver}Overall solver.}
\end{algorithm}

\begin{remark}
 In this work, we consider problems of moderate size where starting from the full finest mest $\Lambda(0) = \{(\overline{J}, \bm{k}) ~ : ~ \bm{k} \in \{0, \dots, N_{\overline{J} - 1} \}^d \}$ is not an issue. However, as emphasized by Cohen \emph{et al.} [8] for Finite Volume, the initial datum could not to fit into memory for large realistic problems, especially when $q$, the number of discrete velocities, is large. We warn the reader that in this case, two approaches are customary.
 The first one is to start from a heuristically adapted mesh $\Lambda(0) \neq \{(\overline{J}, \bm{k}) ~ : ~ \bm{k} \in \{0, \dots, N_{\overline{J} - 1} \}^d \}$ which is coarsened in areas the initial solution is very smooth.
 The second one is to initialize the problem performing the multiresolution analysis on a less refined grid and then to refine every remaining cell a certain number of times and reinitializing the solution on the new cells, see \cite{duarte2012new, duarte2011}.
\end{remark}

The code is sequentially implemented in \texttt{C++} using a code called \texttt{SAMURAI}\footnote{The code with the test cases as well as its documentation can be found at \url{https://github.com/hpc-maths/samurai}.} (\textbf{S}tructured \textbf{A}daptive mesh and \textbf{MU}lti-\textbf{R}esolution based on \textbf{A}lgebra of \textbf{I}ntervals) which is currently under development and than can handle general problems involving dynamically refined Cartesian meshes, stemming both from MR and AMR, in an efficient way.
The data structure \texttt{SAMURAI} relies on intervals of contiguous cells along each Cartesian direction gathered in different sets.
The availability of inter-sets operations allows one to perform all the operations involved in multiresolution and to easily write the corresponding adaptive numerical scheme (lattice Boltzmann, Finite Volume, \dots).
Our implementation relies on the library \texttt{xtensor}\footnote{\url{https://xtensor.readthedocs.io}} for the data storage and in order to take advantage of the lazy evaluations of expressions.

\section{Validation and numerical experiments}\label{sec:VerificationAndExperiments}

To validate our approach, we investigate two types of system for $d = 2$, namely the full Euler system -- taken as the prototype of non-linear hyperbolic system -- on one hand and the incompressible Navier-Stokes equations on the other hand, which represent the archetypal non-linear parabolic system.  
Finally, for $d = 3$, we simulate the advection equation.

For the first test, the purpose is twofold: show that we can ensure a precise error control \cite{bellotti2021} using the procedure described in Section \ref{sec:AdaptiveLatticeBoltzmann}, which is tailored to fit hyperbolic systems; show that in problems involving steep fronts, we can achieve interesting compression rates and thus potentially a reduction of the computation time.
For this test case, we employ a quite unsophisticated vectorial four-velocities scheme (D2Q4).
Using vectorial five-velocities schemes (D2Q5) or a nine-velocities scheme (D2Q9) is straightforward.

The second test aims at showing that our procedure can be coupled with the most well-known \lb scheme -- namely the nine-velocities scheme (D2Q9) -- to tackle a problem which is not typically addressed by multiresolution, because of the smoothness of the solutions induced by the important role of the diffusion.
We compare the adaptive scheme with the reference scheme using integral quantities such as the drag coefficient and the Strouhal number.

For the third test, we employ a six-velocities scheme (D3Q6) and we aim at providing a first example of three-dimensional computation.

All the tests are performed using $\predictionstencildepth = 1$.

\subsection{Systems of hyperbolic conservation laws}\label{sec:EulerTest}

The non-isothermal Euler system is probably the most emblematic example of system of hyperbolic conservation laws.
This kind of problem naturally develops weak solutions with shocks and therefore constitutes the ideal setting to test the efficiency and the reliability of multiresolution.

\subsubsection{Problem}

In the context of $\spatialdimension = 2$, we aim at finding an approximation of the weak entropic solution of the system
\begin{equation}\label{eq:EulerSystem}
    \begin{cases}
        \partialderreduced{t} \vectorial{u} + \partialderreduced{x} \Phi_{x} (\vectorial{u}) + \partialderreduced{y} \Phi_{y}(\vectorial{u}) = 0, \qquad t \in [0, T], \quad &\vectorial{x} \in \reals^2, \\
        \vectorial{u}(t = 0, \vectorial{x}) = \vectorial{u}_0 (\vectorial{x}), \qquad &\vectorial{x} \in \reals^2,
    \end{cases}
\end{equation}
where $\vectorial{u}: [0, T] \times \reals^{2} \to \reals^{m}$, for $m \geq 1$, that we assume to be hyperbolic \cite{toro2013}.
In this paper, we consider some of the test cases introduced by Lax and Liu \cite{lax1998} for the complete Euler system ($m = 4$), given by
\begin{linenomath}\begin{align*}
    \vectorial{u} = \adaptiveroundbrackets{\rho, \rho u, \rho v, E}^{{T}}, \quad \Phi_x(\vectorial{u}) &= \adaptiveroundbrackets{\rho u, \rho u^2 + p, \rho uv, u(E + p)}^{{T}}, \quad 
    \Phi_y(\vectorial{u}) &= \adaptiveroundbrackets{\rho v, \rho uv, \rho v^2 + p, v(E + p)}^{{T}}, 
\end{align*}\end{linenomath}
with pressure law given by $E = \ratioobl{p}{(\gamma_{\text{gas}} - 1)} + \ratioobl{\rho (u^2 + v^2)}{2}$, selecting $\gamma_{\text{gas}} = 1.4$. The computational domain is $\domain = [0, 1]^2$ and we enforce the copy boundary conditions for each variable. The initial datum, given as a function of the density, the velocities and the pressure, is
\begin{linenomath}\begin{equation*}
    \vectorial{u}_0(x, y) = \begin{cases}
                                (\rho_{\text{UR}}, u_{\text{UR}}, v_{\text{UR}}, p_{\text{UR}})^{{T}}, \qquad &x > \ratioobl{1}{2}, \quad y > \ratioobl{1}{2}, \\
                                (\rho_{\text{UL}}, u_{\text{UL}}, v_{\text{UL}}, p_{\text{UL}})^{{T}}, \qquad &x < \ratioobl{1}{2}, \quad y > \ratioobl{1}{2}, \\
                                (\rho_{\text{LL}}, u_{\text{LL}}, v_{\text{LL}}, p_{\text{LL}})^{{T}}, \qquad &x < \ratioobl{1}{2}, \quad y < \ratioobl{1}{2}, \\
                                (\rho_{\text{LR}}, u_{\text{LR}}, v_{\text{LR}}, p_{\text{LR}})^{{T}}, \qquad &x > \ratioobl{1}{2}, \quad y < \ratioobl{1}{2}.
                            \end{cases}
\end{equation*}\end{linenomath}
As test cases, we consider Configuration 3 and Configuration 12 by \cite{lax1998}.
\

\subsubsection{Numerical scheme}

We employ a batch of $m$ D2Q4 schemes with only one conserved moment, coupled through their equilibri in order to achieve the solution of the whole system, see Remark \ref{rem:VectorialSchemes}.
The base D2Q4 scheme is obtained by selecting $\velocitiesnumber = 4$ for the following choice of discrete velocities
\begin{linenomath}\begin{equation*}
    \superscript{\vectorial{\xi}}{\populationindex} = \latticevelocity \adaptiveroundbrackets{\cos{\adaptiveroundbrackets{\ratio{\pi}{2} \populationindex}}, \sin{\adaptiveroundbrackets{\ratio{\pi}{2} \populationindex}}}^{{T}}, \qquad \populationindex = 0, 1, 2, 3,
\end{equation*}\end{linenomath}
regardless of which scheme $i = 1, \dots, m$ we are looking at.
The change of variable yielding the moments is fixed as well
\begin{linenomath}\begin{equation*}
    \operatorial{M} = \adaptiveroundbrackets{
    \begin{matrix}
        1 & 1 & 1 & 1 \\
        \latticevelocity & 0 & -\latticevelocity & 0 \\
        0 & \latticevelocity & 0 & -\latticevelocity\\
        \latticevelocity^2 & -\latticevelocity^2 & \latticevelocity^2 & -\latticevelocity^2
    \end{matrix}}, \qquad
    \operatorial{S}_i = \adaptiveroundbrackets{
    \begin{matrix}
        0 & 0 & 0 & 0 \\
        0 & s_{i}^q & 0 & 0 \\
        0 & 0 & s_{i}^q & 0 \\
        0 & 0 & 0 & s_{i}^{xy}
    \end{matrix}}, \qquad i = 1, \dots, m,
\end{equation*}\end{linenomath}
with the relaxation matrix depending on the component of the system we are considering with $s_{i}^q, s_{i}^{xy} \in (0, 2)$ for $i = 1, \dots, m$. 
We can show, using the equivalent equations \cite{dubois2009}, that selecting the equilibri as
\begin{linenomath}\begin{equation*}
    m_i^{0, \text{eq}} = u_i, \qquad m_i^{1, \text{eq}} = \Phi_{x, i}(\vectorial{u}), \qquad  m_i^{2, \text{eq}} = \Phi_{y, i}(\vectorial{u}), \qquad  m_i^{3, \text{eq}} = 0,
\end{equation*}\end{linenomath}
the scheme is consistent with system \eqref{eq:EulerSystem} up to order $\Delta \timevariable$. 
We have also tested our strategy with the D2Q4 twisted scheme, see Février \cite{fevrier2014}, where the velocities  are not parallel to the axis and are given by
\begin{linenomath}\begin{equation*}
    \superscript{\vectorial{\xi}}{\populationindex} = \latticevelocity \adaptiveroundbrackets{\cos{\adaptiveroundbrackets{\ratio{\pi}{2} \populationindex + \ratio{\pi}{4}}}, \sin{\adaptiveroundbrackets{\ratio{\pi}{2} \populationindex + \ratio{\pi}{4}}}}^{{T}}, \qquad \populationindex = 0, 1, 2, 3,
\end{equation*}\end{linenomath}
observing good behaviors like for the standard D2Q4.
We do not present such tests.
For the examined configurations, we found that $\latticevelocity = 5$ and
\begin{equation*}
    s_i^q = \begin{cases}
                1.9,\qquad &i = 1 \\
                1.75, \qquad &i = 2, 3, 4.
            \end{cases} \qquad \text{and} \qquad
    s_i^{xy} = 1,
\end{equation*}
provide adequate performances and a reasonable amount of numerical diffusion.

\subsubsection{Points of emphasis}

\begin{itemize}
    \item The first important point is to verify that we control the additional error.
    In our previous contribution \cite{bellotti2021}, we have shown that under certain assumptions on the \lb scheme, we are able to control the additional error introduced by the mesh adaptation
    \begin{linenomath}\begin{equation*}
        E[m^{\populationindex}](t) = \frac{\sum_{\vectorial{\indexletter} \in \{0, \dots, N_{\maxlevel} - 1 \}^{\spatialdimension}} \Delta x \adaptiveabs{\reconstructed{m}{\populationindex}{\maxlevel, \vectorial{\indexletter}}(t) - m^{\text{REF}, \populationindex}_{\maxlevel, \vectorial{\indexletter}}(t)}}{\sum_{\vectorial{\indexletter} \in \{0, \dots, N_{\maxlevel} - 1 \}^{\spatialdimension}} \Delta x \adaptiveabs{m^{\text{REF}, \populationindex}_{\maxlevel, \vectorial{\indexletter}}(t)}},
    \end{equation*}\end{linenomath}
    where the solution denoted by ``REF'' is the one by the reference scheme on a uniform mesh at finest level $\maxlevel$ and $\populationindex$ spans the conserved moments.
    We are looking for an upper bound with $C(T) > 0$ of the kind
    \begin{equation}\label{eq:AdditionalErrorEstimate}
        E[m^{\populationindex}](T) \leq C(T) ~ \thresholdletter,
    \end{equation}
    where $T > 0$ is the final time.
    Observe that as we are computing the difference with the reference scheme, the solution of the latter depends on the maximum level of resolution $\maxlevel$. This should be taken into account when comparing results for different $\maxlevel$. 

    \item The second important axis of analysis is the gain in terms of computational time and memory impact thanks to the joint work of the adaptive scheme with multiresolution.
    The memory occupation rate and the mesh occupation rate at time $\timevariable$, given by
    \begin{linenomath}\begin{equation*}
        {\text{MemOR}(t) = \ratio{\#(\text{total cells in }\Lambda(t + \Delta t))}{\#(\text{total cells in the full mesh})}, \qquad \text{MeshOR}(t) = \ratio{\#(\physicalleaves{\Lambda(t + \Delta t)}}{\powertwo{\maxlevel}},}
    \end{equation*}\end{linenomath}
    are taken as metrics as far as the memory impact is concerned.
    {Observe that occupation rates much smaller than one are good and correspond to high compression rates (one has $\text{CR} = 1 - \text{OR}$).}
    Since the parallelization and the optimization of the computational code is not the central topic of this work  and is still undergoing major improvements which shall be a subject on their own, we do not perform an extensive analysis, which could follow the guidelines of \cite{deiterding2016}, of the gain in term of computational time.
    However, notice that since the adaptive scheme only works on the complete leaves of the tree, $\text{MeshOR}(t)$ can be reasonably considered an indicator of the computational cost of the method, apart from the time spent adapting the mesh.
    This will need to be confirmed in future works. However, at each time step, since the number of collisions to perform is equal to $\#(S(\Lambda(t + \Delta t)))$ and $\#(\mathcal{E}_{\levelletter, \vectorial{\indexletter}}^{\populationindex}), \#(\mathcal{A}_{\levelletter, \vectorial{\indexletter}}^{\populationindex}) \sim 2^{(\spatialdimension-1)(\maxlevel - \levelletter)} \ll 2^{\spatialdimension (\maxlevel - \levelletter)}$ if the level $\levelletter$ is far enough from $\maxlevel$, corresponding to the number of reconstruction to be done for the stream phase, we believe that the number of operations to perform in one time step can be made significantly smaller than for the reference scheme. This is particularly flagrant when handling large three-dimensional problems like in Section \ref{sec:3dAdvection}.
    
\end{itemize}

\subsubsection{Results and discussion}

\begin{figure}
\begin{center}
        \includegraphics[width=1.\textwidth]{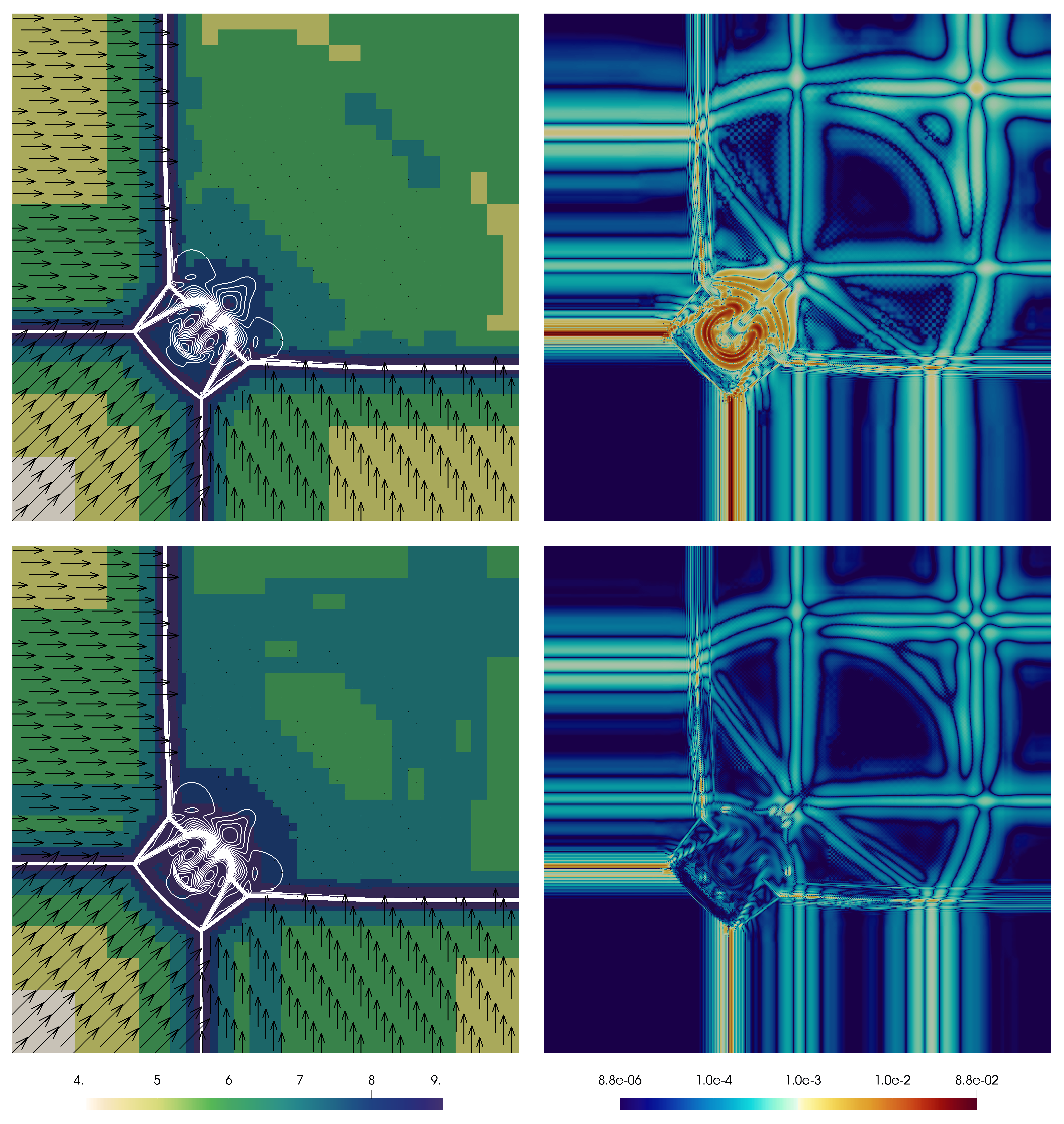}
    \end{center}\caption{\label{fig:Configuration3_error}Configuration 3. Top $\thresholdletter = 0.005$, bottom $\thresholdletter = 0.001$ and $\regularityguess = 0$. On the left, level $\levelletter$ (colored), contours of the density field (white) and velocity field (black). On the right: local relative error of the adaptive method with respect to the reference method. Time $T = 0.3$, $\minlevel = 2$ and $\maxlevel = 9$.}
\end{figure}

\begin{figure}
    \begin{center}
        \includegraphics[width=1.\textwidth]{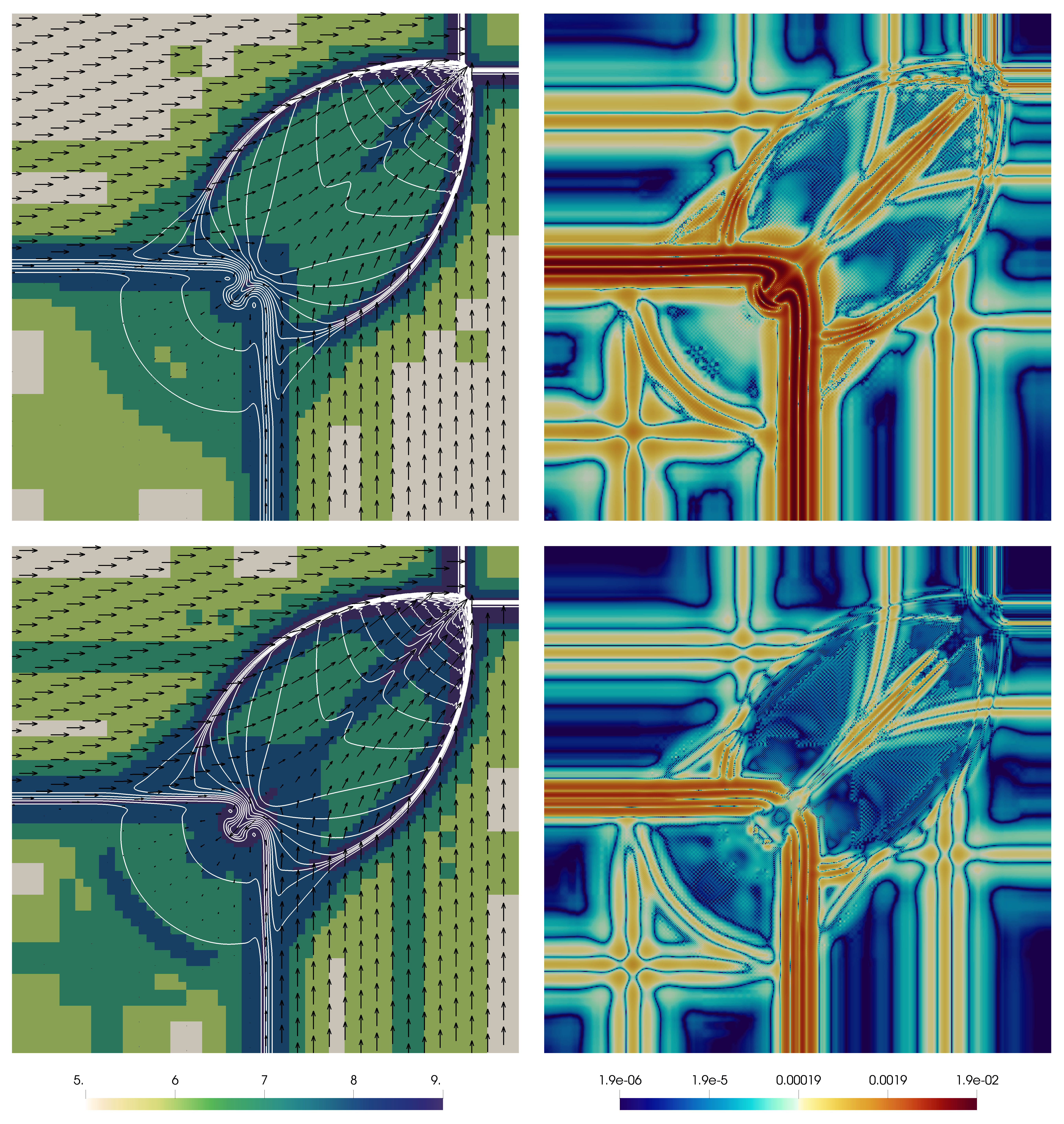}
    \end{center}\caption{\label{fig:Configuration12_error}Configuration 12. Top $\thresholdletter = 0.005$, bottom $\thresholdletter = 0.001$ and $\regularityguess = 0$. On the left, level $\levelletter$ (colored), contours of the density field (white) and velocity field (black). On the right: local relative error of the adaptive method with respect to the reference method. Time $T = 0.25$, $\minlevel = 2$ and $\maxlevel = 9$.}
\end{figure}

\begin{figure}
    \begin{center}
        \includegraphics[width=1.\textwidth]{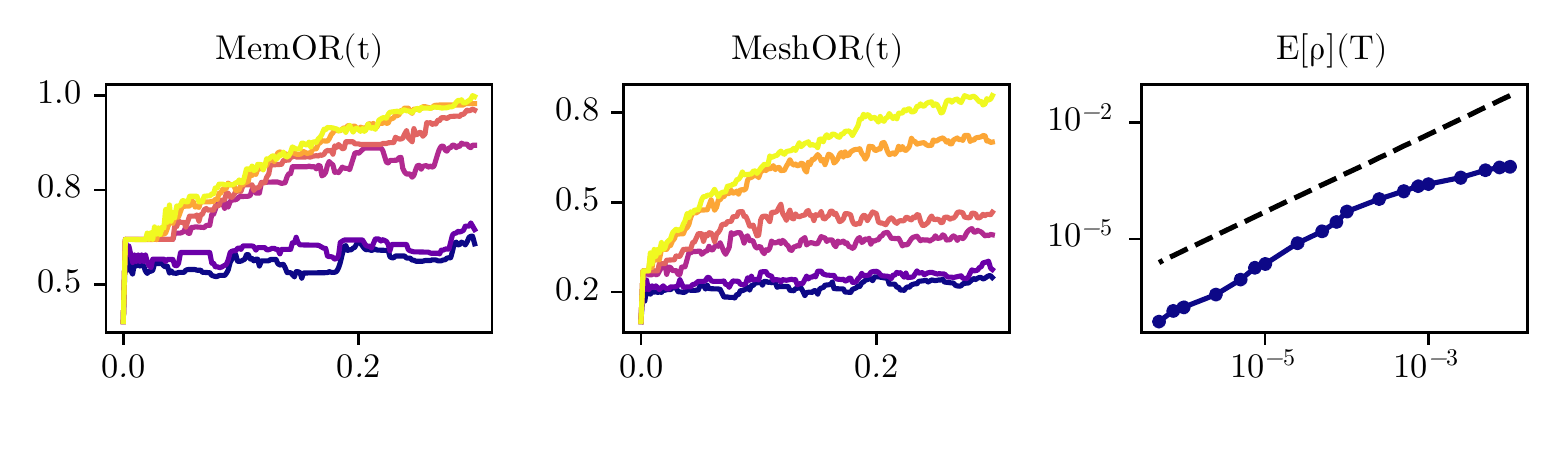}
        \includegraphics[width=1.\textwidth]{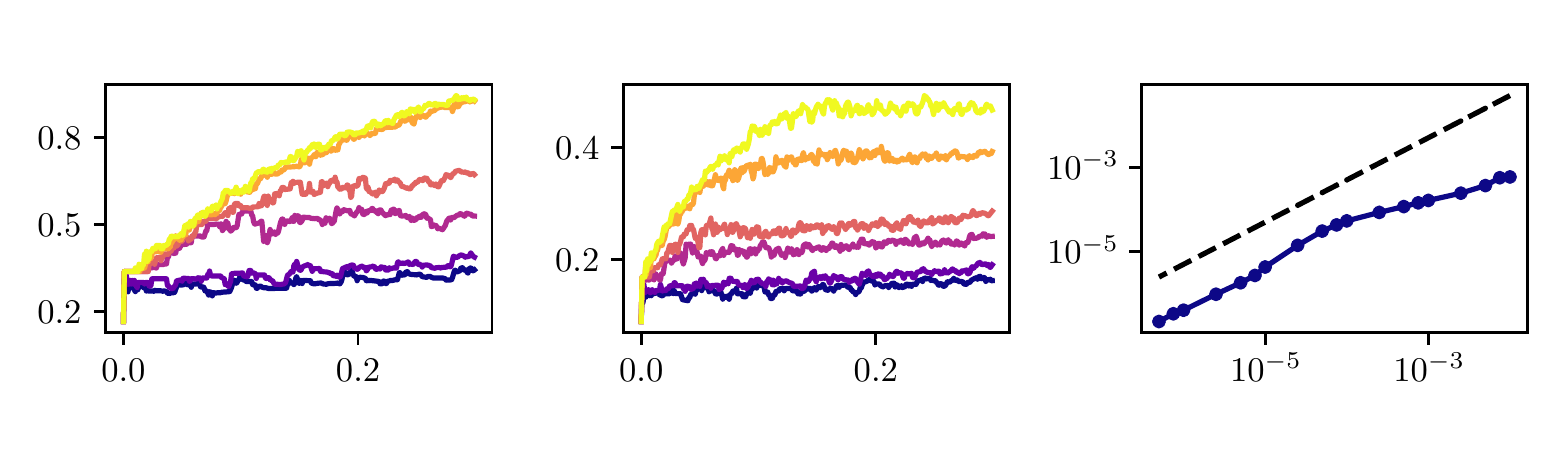}
        \includegraphics[width=1.\textwidth]{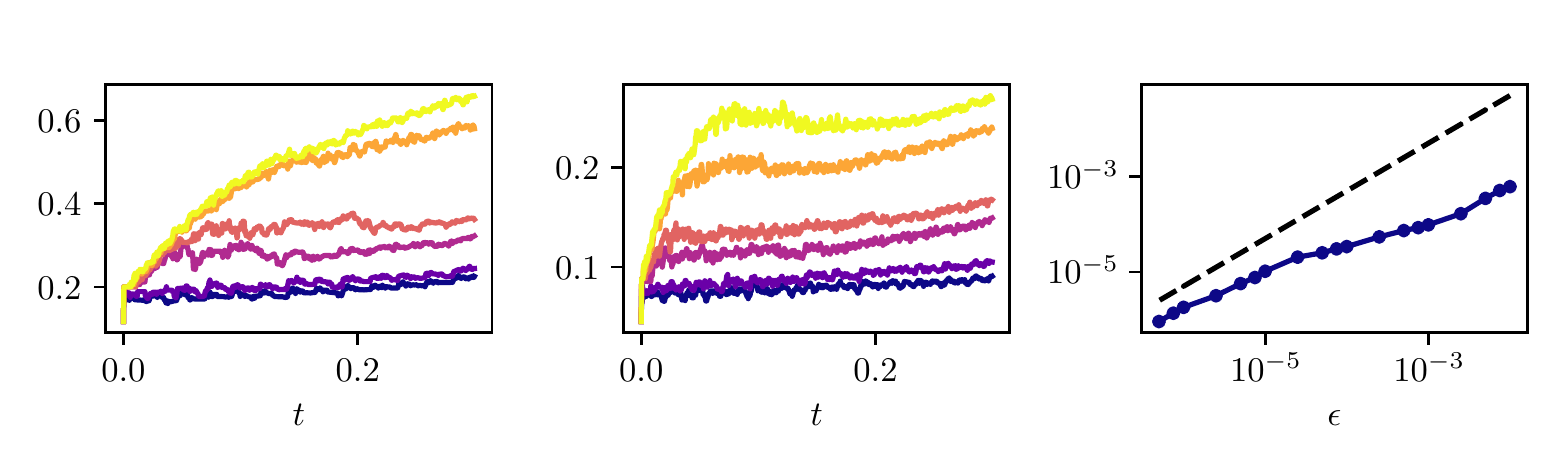}
    \end{center}\caption{\label{fig:Configuration3_eps}Configuration 3 with $\minlevel = 2$ and $\regularityguess = 0$. From top to bottom $\maxlevel = 7, 8$ and $9$.
    From left to right: memory occupation and mesh occupation as functions of time and additional error as function of $\thresholdletter$. For the two plots on the left: \textcolor{legendBlue}{$\blacksquare$} $\thresholdletter = 10^{-2}$, \textcolor{legendViolet}{$\blacksquare$}  $\thresholdletter = 5 \cdot 10^{-3}$, \textcolor{legendMagenta}{$\blacksquare$}  $\thresholdletter = 10^{-3}$, \textcolor{legendRed}{$\blacksquare$}  $\thresholdletter = 5 \times 10^{-4}$, \textcolor{legendOrange}{$\blacksquare$}  $\thresholdletter = 10^{-4}$, \textcolor{legendYellow}{$\blacksquare$}  $\thresholdletter = 5 \times 10^{-5}$. The dashed black line gives the slope $\thresholdletter$.}
\end{figure}

\begin{figure}
    \begin{center}
        \includegraphics[width=1.\textwidth]{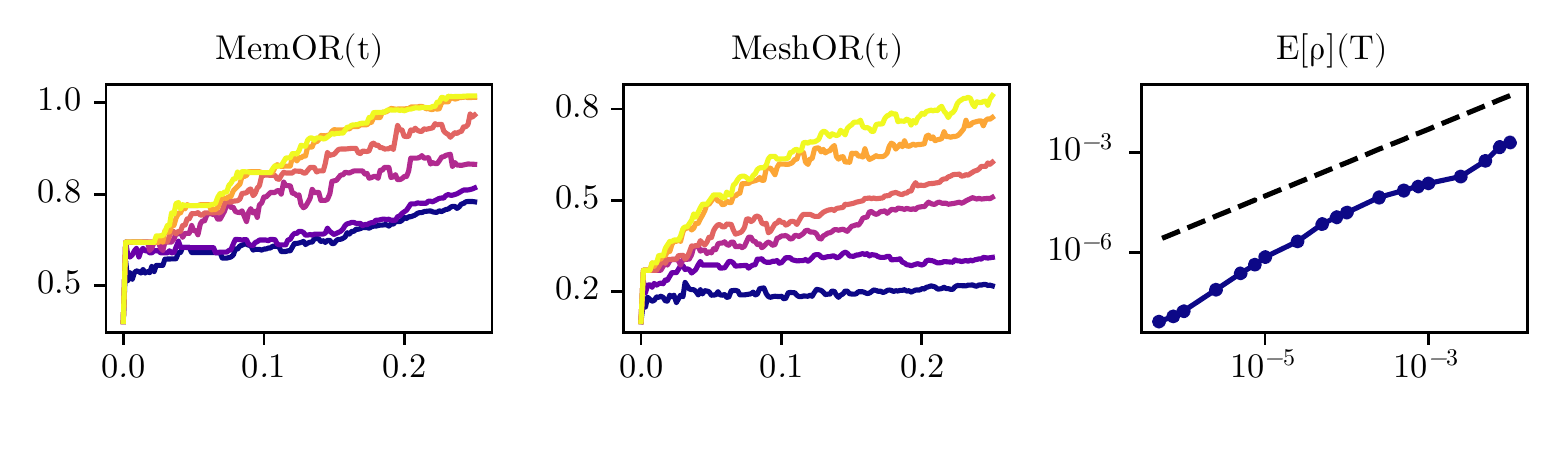}
        \includegraphics[width=1.\textwidth]{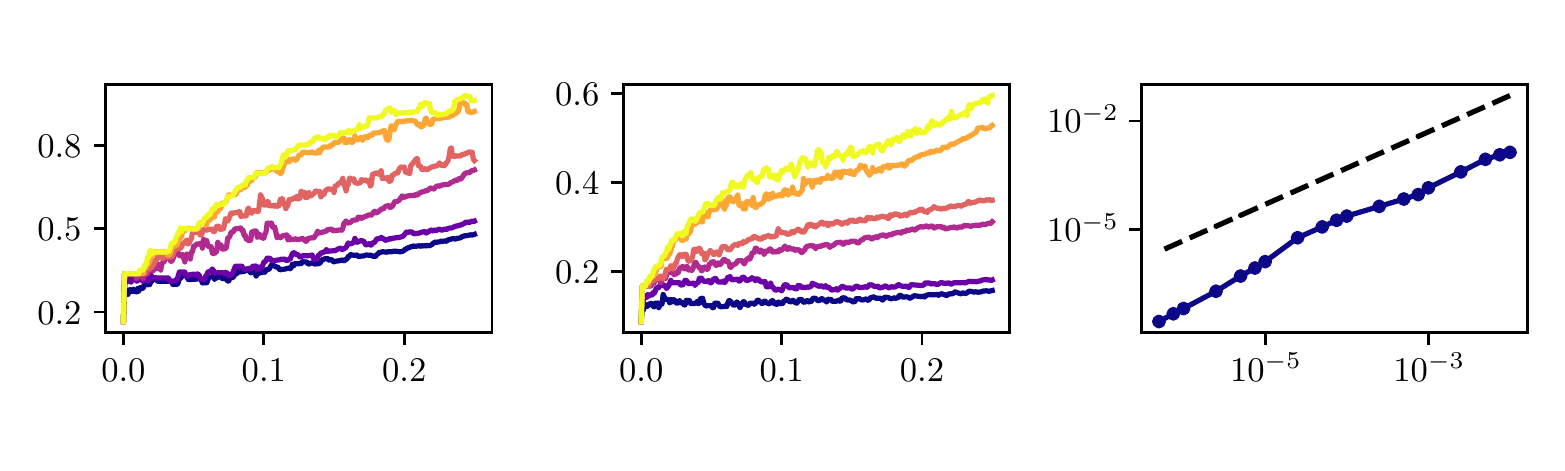}
        \includegraphics[width=1.\textwidth]{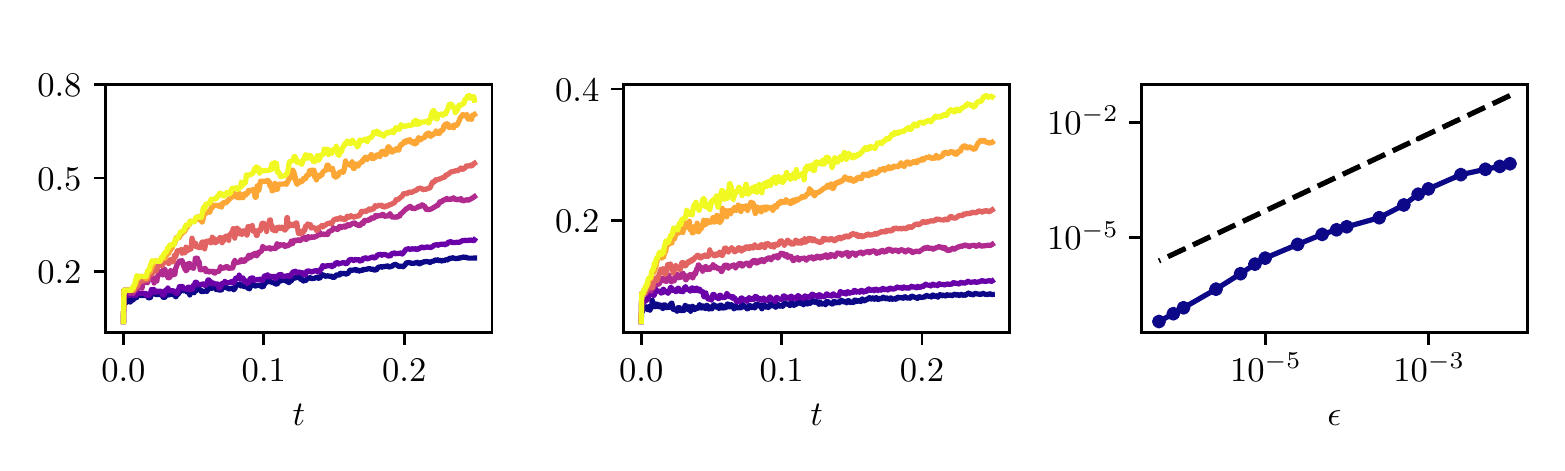}
    \end{center}\caption{\label{fig:Configuration12_eps}Configuration 12 with $\minlevel = 2$ and $\regularityguess = 0$. From top to bottom $\maxlevel = 7, 8$ and $9$.
    From left to right: memory occupation and mesh occupation as functions of time and additional error as function of $\thresholdletter$. For the two plots on the left: \textcolor{legendBlue}{$\blacksquare$} $\thresholdletter = 10^{-2}$, \textcolor{legendViolet}{$\blacksquare$}  $\thresholdletter = 5 \cdot 10^{-3}$, \textcolor{legendMagenta}{$\blacksquare$}  $\thresholdletter = 10^{-3}$, \textcolor{legendRed}{$\blacksquare$}  $\thresholdletter = 5 \times 10^{-4}$, \textcolor{legendOrange}{$\blacksquare$}  $\thresholdletter = 10^{-4}$, \textcolor{legendYellow}{$\blacksquare$}  $\thresholdletter = 5 \times 10^{-5}$. The dashed black line gives the slope $\thresholdletter$.}
\end{figure}

\paragraph{General remarks}

The structure of the solution and the additional error for Configuration 3 are given in Figure \ref{fig:Configuration3_error} at final time $T = 0.3$; those for Configuration 12 are given in Figure \ref{fig:Configuration12_error} at final time $T = 0.25$. 
In the former case, we remark that the four shocks, where all the conserved moments are discontinuous, are well resolved and finely meshed, so that we can observe some hydrodynamic instabilities typical of such systems\footnote{With the limitation of dealing with a low order scheme.} (see \cite{liska2003} for example).
In the latter case, the two shocks propagating towards the upper-right corner are followed by the finest discretization of the mesh, whereas we observe a coarsening of one level (for $\thresholdletter = 5\cdot 10^{-3}$) close to the static contact discontinuities. This phenomenon shall be clarified in a moment with a finer analysis.
Overall, this qualitative analysis allows us to conclude that the adaptive \lb scheme succeeds in reproducing the expected behavior of the solution \cite{lax1998, liska2003} of the Euler system and that the adaptive grid follows the shock structures propagating with finite velocity.

\paragraph{Error control}

According to the results shown in Figure \ref{fig:Configuration3_eps} and \ref{fig:Configuration12_eps}, respectively for Configuration 3 and 12, we verify that the upper bound \eqref{eq:AdditionalErrorEstimate} with $\thresholdletter$ is verified for these test cases.
This is in agreement with the theoretical analysis made in \cite{bellotti2021} and holds for any choice of finest resolution $\maxlevel$.

Besides corroborating the theoretical analysis in a multidimensional setting, this confirms that, as far as we are dealing with hyperbolic problems, the enlargement strategy devised in Section \ref{sec:TreeEnlargement} is capable of ensuring that the adaptive mesh correctly follows the temporal evolution of the solution.

\paragraph{Compression}

Still looking at Figure \ref{fig:Configuration3_eps} and \ref{fig:Configuration12_eps}, we observe that the occupation rates, namely $\text{MemOR}(t)$ and $\text{MeshOR}(t)$ become more interesting as one approaches larger maximum levels of resolution $\maxlevel$. This is in accordance with the intuition that since a fine sampling of the solution is needed only close to the shocks, the number of needed leaves shall grow in $\maxlevel$ more slowly than $2^{\spatialdimension \maxlevel}$, also because shock are $(\spatialdimension-1)$-dimensional entities. 

As far as time is concerned, after an important initial growth guided by the refinement criterion $\mathcal{H}_{\thresholdletter}$, the trend of the occupation rates eventually stabilizes, especially for Configuration 3, where the secondary structures close to the hydrodynamic instabilities induced by the contact discontinuities do not grow too much in size as time advances.
For Configuration 12, the occupation rates grow linearly even close to the final time because of the expanding curvilinear front linking the main shocks and the contact discontinuities.

Finally, the occupation rates are less interesting once we decrease $\thresholdletter$. In these tests, we deliberately used unrealistically small $\thresholdletter$ only aiming at showing convergence.
We conclude that the choice of $\thresholdletter$ should be the result of an arbitration between a desired target error and performances.

\paragraph{Coupling between MR and LBM. Role of $\regularityguess$}

\begin{figure}
    \begin{center}
        \includegraphics[width=0.8\textwidth]{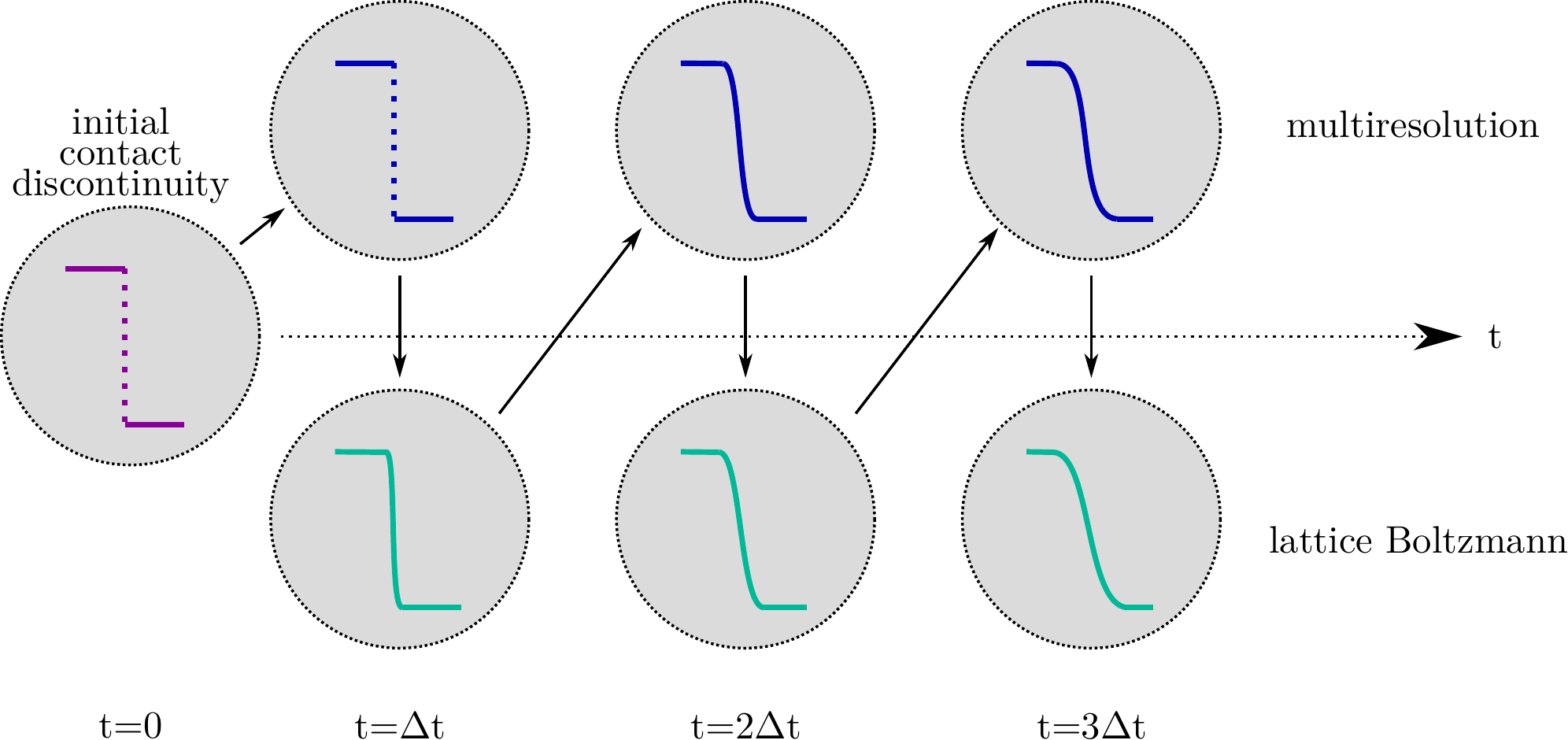}
    \end{center}\caption{\label{fig:brouillon_contact_discontinuity}Coupling between the \lb scheme and the multiresolution to amplify the errors on the contact discontinuities.}
\end{figure}

\begin{figure}
    \begin{center}
        \includegraphics[width=0.85\textwidth]{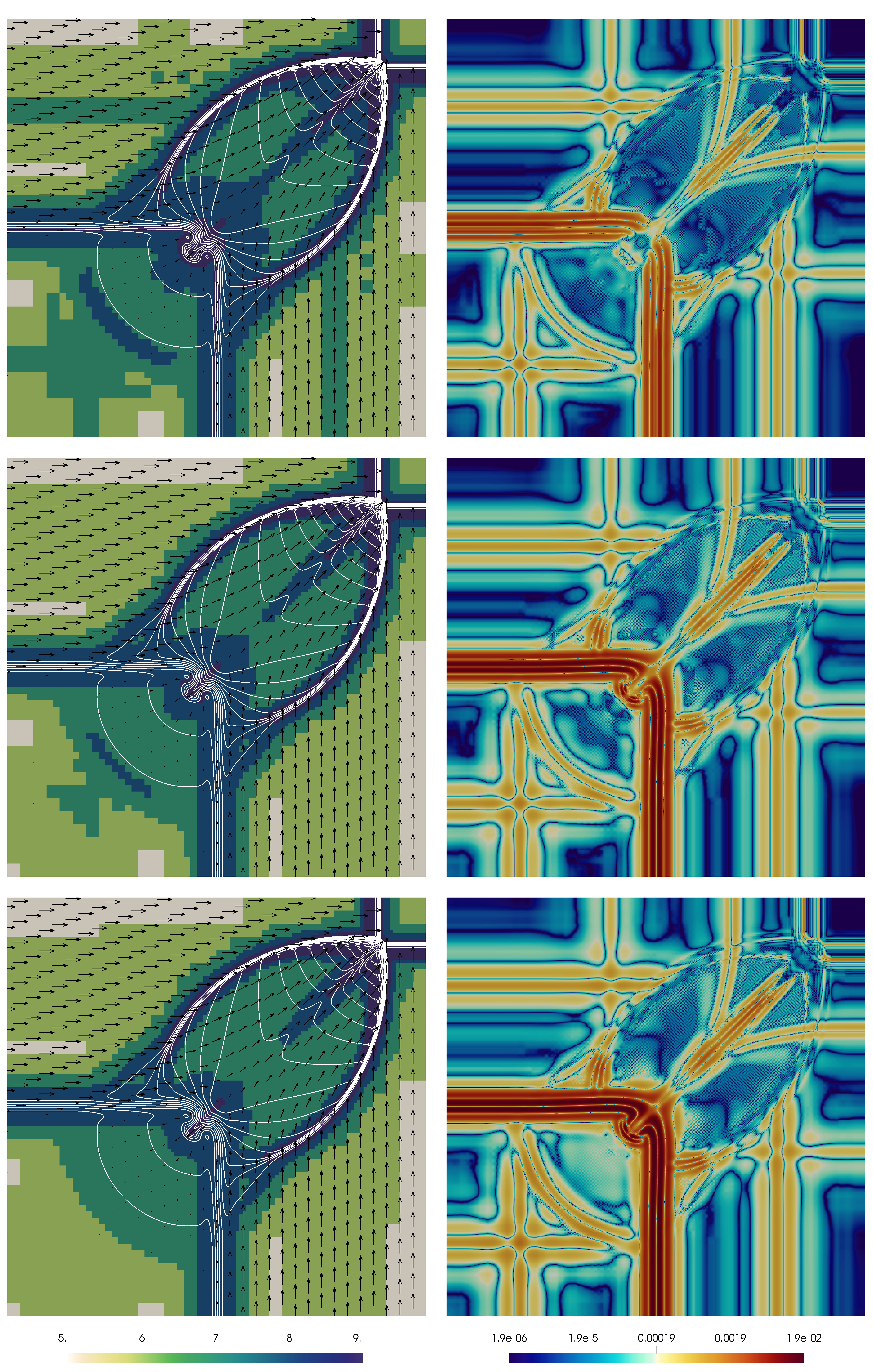}
    \end{center}\caption{\label{fig:Configuration12_varying_mu}Configuration 12. For all the simulations, time $T = 0.25$, $\minlevel = 2$ and $\maxlevel = 9$, $\epsilon = 0.001$. On the left, level $\levelletter$ (colored), contours of the density field (white) and velocity field (black). On the right: local relative error of the adaptive method with respect to the reference method. From top to bottom $\regularityguess = 0, 1$ and $2$.}
\end{figure}

As it has been hinted while discussing Figure \ref{fig:Configuration12_error}, we remark that the larger errors for Configuration 12 are situated close to the static contact discontinuities.
This also holds, but less spectacularly, for Configuration 3 in the area where contact discontinuities are present causing the hydrodynamic instabilities to appear. 
This is an interesting coupling phenomenon between the poor behavior of the reference scheme on the contact discontinuities -- which is inherent to this class of vectorial schemes (see Graille \cite{graille2014}) -- and the multiresolution. Indeed, the reference scheme smears the contact discontinuities (decreases the magnitude of the details in these areas) from the very beginning and then the multiresolution adaptation coarsens the mesh causing a local accumulation of error in time. This pattern is schematized in Figure \ref{fig:brouillon_contact_discontinuity}.

We also verified that the fact of performing the collision on the complete leaves without reconstruction, see Section \ref{sec:Collision}, has a negligible impact on this particular phenomenon, even if the equilibrium functions are strongly non-linear.

As Figure \ref{fig:Configuration12_error} and \ref{fig:Configuration3_error} show, this problem is clearly alleviated by decreasing the quality factor $\thresholdletter$, from $0.005$ to $0.001$. 
Still, it is more interesting to study how the regularity guess $\regularityguess$ for the solution influences the behavior of the hybrid scheme with respect to this issue.
To study this for Configuration 12, we fixed $\epsilon = 0.001$ and varied $\regularityguess = 0, 1$ and $2$, as shown in Figure \ref{fig:Configuration12_varying_mu}. We observe that this parameter, involved in the refinement process $\mathcal{H}_{\thresholdletter}$, does not affect the structures which are already well refined after the coarsening process $\mathcal{T}_{\thresholdletter}$, namely the shocks.
Close to these structures, if we assume that \eqref{eq:DetailsDecayEstimate} is sharp, the details do not decrease with $\levelletter$ and thus precision is ensured if $\thresholdletter$ is reasonably small.
On the other hand, we observe that a wise choice of $\regularityguess$ is effective in diminishing the coupling effect between the multiresolution and the \lb scheme in smearing contact discontinuities (consider that the color-scale is logarithmic), without having to drastically reduce $\thresholdletter$ causing a degradation of the performances of the algorithm.
Since $\regularityguess$ represents the number of bounded derivatives of the expected solution of the problem, the advice for solutions developing shocks and contact discontinuities (thus only $L^{\infty}$ solutions) is to set $\regularityguess = 0$. This has proved to allow for a reduction of the artificial smearing of the contact discontinuities by the numerical method.

\FloatBarrier

\subsection{Incompressible Navier-Stokes system}\label{sec:NavierStokesTest}

The simulation of the incompressible Navier-Stokes equations is probably the most well-known application of the \lb schemes.
Though the study of this problem with \lb schemes is limited to low Mach and Reynolds numbers, thus limiting the interest of multiresolution, we still perform this test to showcase the applicability of our strategy to this class of problems. 

\subsubsection{Problem}

We consider the problem of a flow around an obstacle occupying the open set $\Theta \subset \mathbb{R}^2$, with flow supposed to be incompressible and considering a Newtonian fluid
\begin{equation}
    \begin{cases}
        \nabla \cdot \vectorial{u} = 0, \qquad &t \geq 0, \quad \vectorial{x} \in \reals^2 \smallsetminus \Theta, \\
        \rho_0 \adaptiveroundbrackets{\partialderreduced{t} \vectorial{u} + \vectorial{u} \cdot \nabla \vectorial{u}} = -\nabla p + \nabla \cdot \adaptiveroundbrackets{2\mu_{\text{visc}} \ratio{\nabla \vectorial{u} + \nabla \vectorial{u}^T}{2}}, \qquad &t \geq 0, \quad \vectorial{x} \in \reals^2 \smallsetminus \Theta, \\
        \vectorial{u} = \vectorial{u}_0 = (u_0, 0), \qquad &t = 0, \quad \vectorial{x} \in \reals^2 \smallsetminus \Theta, \\
        \vectorial{u} = 0, \qquad &t \geq 0, \quad \vectorial{x} \in \partial \Theta.
    \end{cases}
\end{equation}
The leading dimensionless quantity in this problem, determining the flow regime, is the Reynolds number
\begin{linenomath}\begin{equation*}
    \text{Re} = \ratio{\rho_0 u_0 L_{\Theta}}{\mu_{\text{visc}}},
\end{equation*}\end{linenomath}
given the characteristic length of the obstacle $L_{\Theta} \sim (|\Theta|_2)^{1/2}$.
It is well-known that for $\text{Re} > 90$, the flow passes from being fully laminar to a periodic regime where vortices periodically shed. This phenomenon is known as von Kármán vortex street.

\subsubsection{Numerical scheme}

The numerical scheme we test is the popular D2Q9 scheme by Lallemand and Luo \cite{lallemand2000}\footnote{Many other schemes are proposed in the literature and we tested some of them with similar results.}, with $q = 9$ velocities given by
\begin{linenomath}\begin{equation*}
    \superscript{\vectorial{\xi}}{\populationindex} = 
        \begin{cases}
            \adaptiveroundbrackets{0, 0}^T , \qquad &\populationindex = 0,\\
            \latticevelocity \adaptiveroundbrackets{\cos{\adaptiveroundbrackets{\ratio{\pi}{2} (\populationindex - 1)}}, \sin{\adaptiveroundbrackets{\ratio{\pi}{2} (\populationindex - 1)}}}^T, \qquad &\populationindex = 1,2,3,4,\\
            \latticevelocity \adaptiveroundbrackets{\cos{\adaptiveroundbrackets{\ratio{\pi}{2} (\populationindex - 5) + \ratio{\pi}{4}}}, \sin{\adaptiveroundbrackets{\ratio{\pi}{2} (\populationindex - 5) + \ratio{\pi}{4}}}}^T, \qquad &\populationindex = 5,6,7,8,
        \end{cases}
\end{equation*}\end{linenomath}
with relaxation $\operatorial{S} = \text{diag} \adaptiveroundbrackets{0, 0, 0, s_1, s_1, s_1, s_1, s_2, s_2}$ and change of basis given by
\begin{linenomath}\begin{equation*}
    \operatorial{\changeofbasisletter} = \adaptiveroundbrackets{
    \begin{matrix}
        1 &  1  &  1   &  1  &  1  &    1   &   1   &   1   &   1   \\
        0 & \lambda  &  0   & -\lambda &  0  &   \lambda   &  -\lambda  &  -\lambda  &  \lambda   \\
        0 &  0  &  \lambda  &  0  & -\lambda &   \lambda   &  \lambda   &  -\lambda  &  -\lambda  \\
        -4\lambda^2 & -\lambda^2  & -\lambda^2   & -\lambda^2  & -\lambda^2  &  2\lambda^2  & 2\lambda^2  & 2\lambda^2  & 2\lambda^2  \\
        0 &  -2\lambda^3  &  0   &  2\lambda^3  &  0  &   \lambda^3   & -\lambda^3   & -\lambda^3   &  \lambda^3   \\
        0 &  0  &  -2\lambda^3   &  0  &  2\lambda^3  &   \lambda^3   &  \lambda^3   & -\lambda^3   & -\lambda^3   \\
        4\lambda^4 &  -2\lambda^4  &  -2\lambda^4   &  -2\lambda^4 &  -2\lambda^4  &   \lambda^4   &  \lambda^4   &  \lambda^4   &  \lambda^4   \\
        0 & \lambda^2  & -\lambda^2  & \lambda^2  & -\lambda^2 &    0   &   0   &   0   &   0   \\
        0 &  0  &  0   &  0  &  0  &   \lambda^2   & -\lambda^2   &  \lambda^2   & -\lambda^2  
    \end{matrix}}.
\end{equation*}\end{linenomath}
The moments at the equilibrium are, indicating $c_s = \ratioobl{\latticevelocity}{\sqrt{3}}$, $\subscript{\velocity}{x} = \ratioobl{\subscript{q}{x}}{\density}$, $\subscript{\velocity}{y} = \ratioobl{\subscript{q}{y}}{\density}$ and $\adaptiveabs{\vectorial{\velocity}}^2 = \power{\subscript{\velocity}{x}}{2} + \power{\subscript{\velocity}{y}}{2}$
\begin{linenomath}\begin{equation*}
    \superscript{\vectorial{m}}{\text{eq}} = \adaptiveroundbrackets{\density, \subscript{q}{x}, \subscript{q}{y}, 3(-2 c_s^2 \rho + |\vectorial{q}|/\rho), -3 c_s^2 q_x, -3 c_s^2 q_y, 9 (c_s^4 - c_s^2 |\vectorial{q}|/\rho), (q_x^2 - q_y^2)/\rho, q_x q_y/\rho) }^T.
\end{equation*}\end{linenomath}
Assuming to be in the low-Mach setting, namely $|\vectorial{u}_0| / c_s \ll 1$, it can be shown \cite{fevrier2014} that we obtain a quasi-incompressible regime ($\rho \simeq \rho_0$) and that the equivalent equations are, neglecting parasitic terms proportional to the cube of the velocity
\begin{linenomath}\begin{equation*}
    \begin{cases}
        \partialderreduced{\timevariable} \density + \nabla \cdot \adaptiveroundbrackets{\density \vectorial{\velocity}} = O(\Delta t^2), \\
        \partialderreduced{\timevariable} \adaptiveroundbrackets{\density \vectorial{\velocity}} + \nabla \cdot \adaptiveroundbrackets{\density \vectorial{\velocity} \otimes \vectorial{\velocity}} = - \nabla (c_s^2 \rho ) + c_s^2 \Delta t \nabla \cdot \adaptiveroundbrackets{2 \rho {\sigma_2} \adaptiveroundbrackets{\ratio{\nabla \vectorial{\velocity} + \superscript{\nabla \vectorial{\velocity}}{T}}{2}} +  \rho {\sigma_1} \adaptiveroundbrackets{\nabla \cdot \vectorial{\velocity}} \operatorial{I}} + O(\Delta t^2),
    \end{cases}
\end{equation*}\end{linenomath}
with $\sigma_k \definitionequality \adaptiveroundbrackets{\ratio{1}{s_k} - \ratio{1}{2}}$.
Thus, to enforce the diffusivity $\mu_{\text{visc}}$, one has to set $s_2$
\begin{linenomath}\begin{equation*}
    s_2 = \adaptiveroundbrackets{\ratio{1}{2} + \ratio{\mu_{\text{visc}}}{c_s^2 \Delta t \rho_0}}^{-1}.
\end{equation*}\end{linenomath}

Special care must be devoted to the representation of the obstacle $\Theta$.
Usually, one employs bounce back boundary condition to enforce zero velocity on $\partial \Theta$ and there exist a vast specialized literature on this matter.
However, in this work, we do not want to adapt the  domain and its discretization to fit $\partial \Theta$ and so we proceed in the following way: we mesh the entire space $\domain = [0, 2] \times [0, 1]$ without considering the obstacle and we apply the adaptive scheme as always.
At the end of each time step, for the leaves $\subscript{C}{\levelletter, \vectorial{\indexletter}}$ intersecting $\Theta$, we estimate $|\subscript{C}{\levelletter, \vectorial{\indexletter}} \cap \Theta|_d$ and we set
\begin{linenomath}\begin{equation*}
    \subscript{\vectorial{f}}{\levelletter, \vectorial{\indexletter}}(t + \Delta t) = \ratio{|\subscript{C}{\levelletter, \vectorial{\indexletter}} \cap \Theta|_d}{|\subscript{C}{\levelletter, \vectorial{\indexletter}}|_d} \operatorial{M}^{-1} \adaptiveroundbrackets{\superscript{\vectorial{m}}{\text{eq}}|_{\substack{\rho = \rho_0 \\ u_x = 0 \\ u_y = 0}}} + \adaptiveroundbrackets{1 - \ratio{|\subscript{C}{\levelletter, \vectorial{\indexletter}} \cap \Theta|_d}{|\subscript{C}{\levelletter, \vectorial{\indexletter}}|_d}} \subscript{\vectorial{f}}{\levelletter, \vectorial{\indexletter}}(t + \Delta t),
\end{equation*}\end{linenomath}
following the direction of Mohamad and Succi \cite{mohamad2009}.
For the external boundaries, we impose a bounce back boundary condition with the given velocity $\vectorial{u}_0$ on the left, top and bottom boundary and a copy boundary condition on the outlet.

\subsubsection{Points of emphasis}

Due to the important role of the viscosity in the whole domain and not only where the mesh is more refined, we cannot compare the solution at each time frame because the von Kármán instability is going to develop differently (see comments).
Thus, we analyze some integral quantities \cite{eitel2013} like the drag coefficient $C_D$, the lift coefficient $C_L$ and the Strouhal number, given by
\begin{linenomath}\begin{equation*}
    C_D = \frac{2 F_x}{\rho_0 u_0^2 L_{\Theta}}, \qquad C_L = \frac{2 F_y}{\rho_0 u_0^2 L_{\Theta}}, \qquad \text{St} = \ratio{L_{\Theta} \omega}{u_0},
\end{equation*}\end{linenomath}
where $\vectorial{F} = (F_x, F_y)$ is the force acting on the obstacle and $\omega$ is the shedding frequency of the vortices.
In particular, the shedding frequency is computed using the fast Fourier transform on the available lift coefficient.
We are interested in comparing these dimensionless quantities between the reference method on a uniform mesh and the adaptive method on the evolving adaptive mesh, as well as the occupation rates  introduced in the previous Section.

\subsubsection{Results and discussion}

\begin{figure}
    \begin{center}
        \includegraphics[width=1.\textwidth]{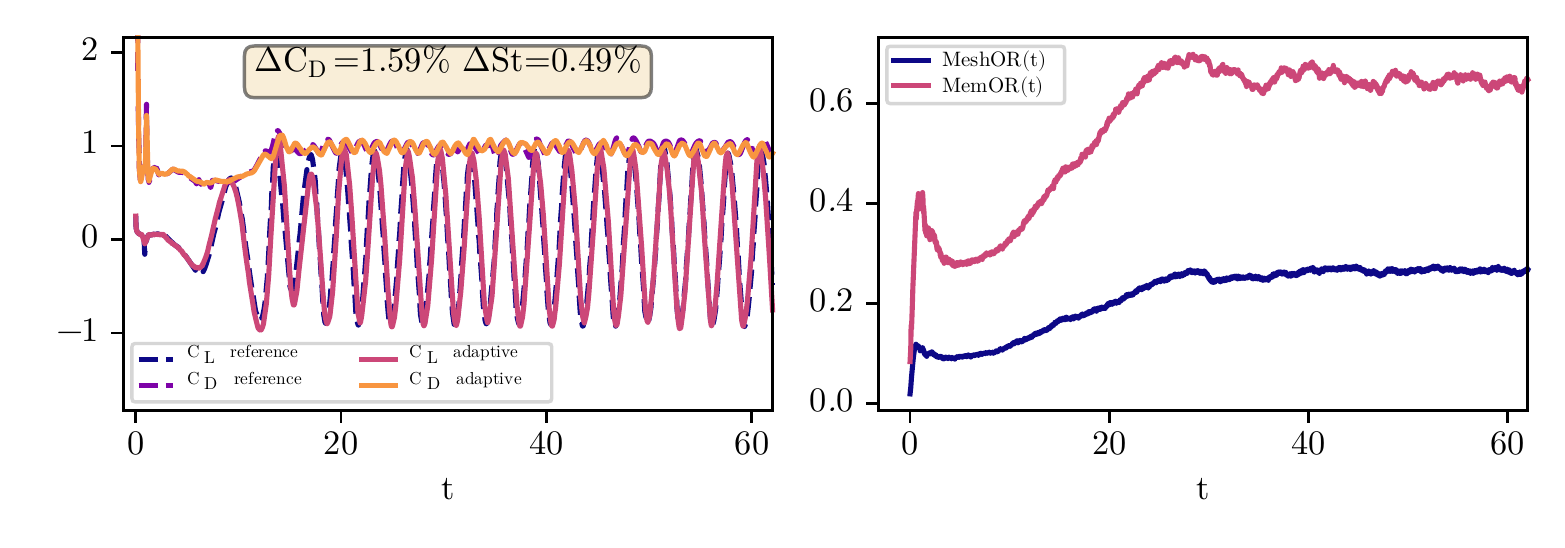}
        \includegraphics[width=1.\textwidth]{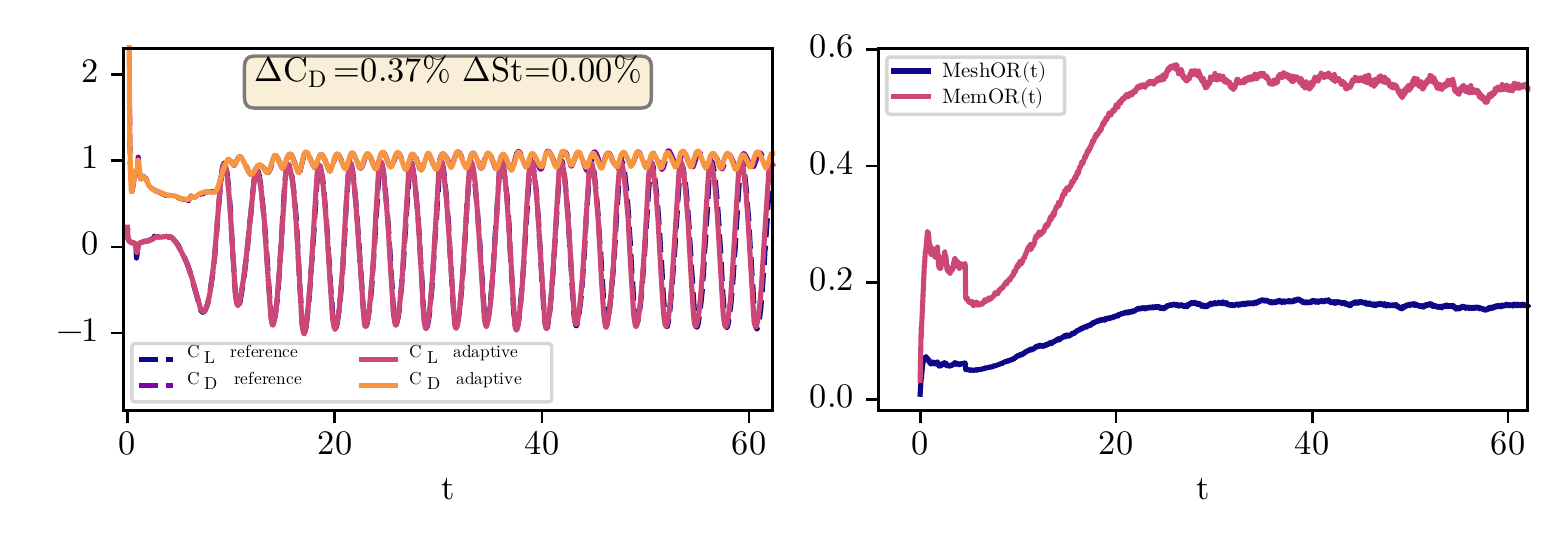}
        \includegraphics[width=1.\textwidth]{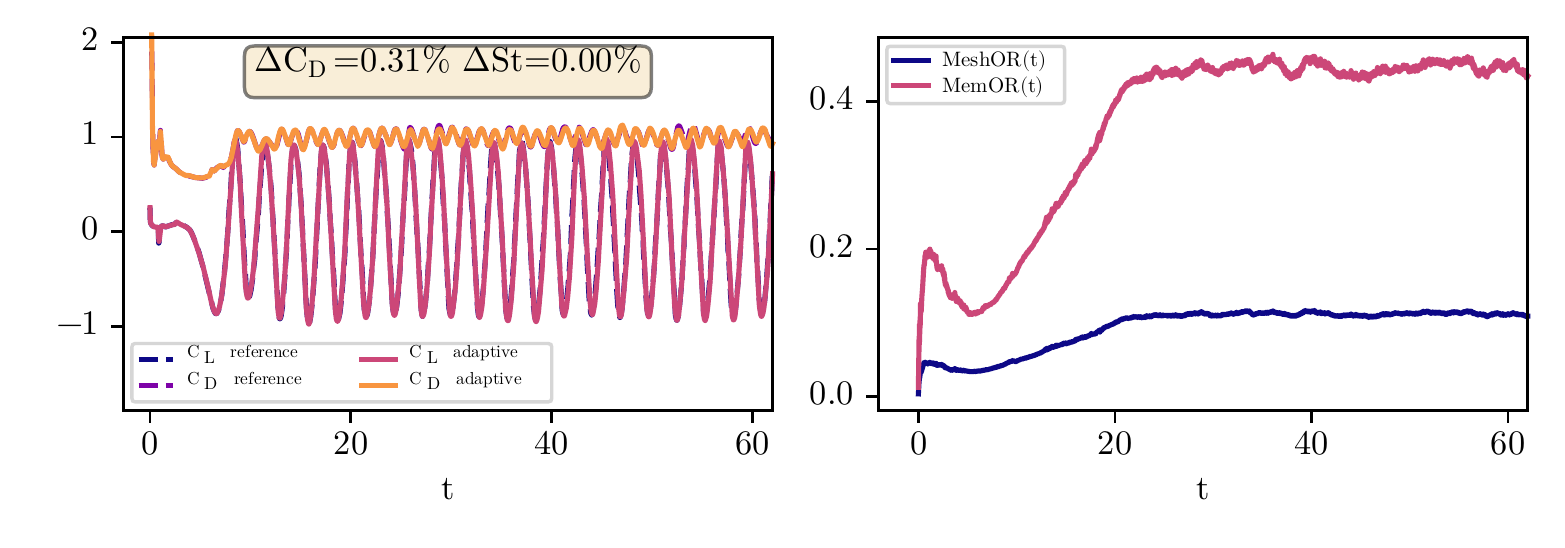}
    \end{center}\caption{\label{fig:1200_integrals_all}On the left, time behavior of the drag and lift coefficient for the reference and the adaptive scheme and on the right, mesh and memory occupation rates, for $\text{Re} = 1200$, $\minlevel = 2$ and $\regularityguess = 1$.
    From top to bottom: $(\maxlevel, \thresholdletter) = (7, 0.00075)$, $(8, 0.000375)$ and $(9, 0.000175) $. The drag coefficient has been normalized with its average for the reference scheme and the lift coefficient with the maximum for the reference scheme.}
\end{figure}

We test the method for $\text{Re} = 1200$ (we also tested it for other moderate Reynolds numbers, obtaining analogous results) with $\lambda = 1$, fixing $\minlevel = 2$ and using $\maxlevel = 7, 8$ and $9$ with $\thresholdletter = 0.00075, 0.000375$ and $0.000175$ respectively and selecting $\regularityguess = 1$.
The reason why we have chosen to decrease $\thresholdletter$ as $\maxlevel$ increases shall be clear in a moment.

On Figure \ref{fig:1200_integrals_all} we observe an excellent agreement on the value of the integral quantities between the adaptive method and the reference method. The discrepancies are obviously reduced as $\maxlevel$ increases because of the variation of the respective threshold parameter $\thresholdletter$.
This comes from the fact that the most important contribution to these quantities comes from the area around the obstacle $\Theta$, where the flow regime can be considered to be, to some extent, highly inertial (or hyperbolic). For this kind of regime, previous works \cite{cohen2003} and \cite{bellotti2021} and the previous Section have shown that the Harten heuristics is respected with our choice of $\mathcal{T}_{\thresholdletter}$ and $\mathcal{H}_{\thresholdletter}$. 
The occupation rates are interesting and become better and better with $\maxlevel$ as we also observed for the solution of the Euler system, despite the fact that we reduced the threshold parameter $\thresholdletter$ as we increased $\maxlevel$. 
The initial growth of these rates followed by a decrease is due to some initial acoustic waves quickly propagating radially in the domain which are eventually damped by the external boundary conditions and are inherent to the \lb method and its way of treating boundary conditions.
These waves do not affect computations on a longer time scale.
The values of the occupation rates stabilize after the complete onset of the instability followed by the periodic regime, that is for $t > 20$, because the number of vortices present in the domain is roughly the same for every time.

\begin{figure}
    \begin{center}
        \includegraphics[width=1.\textwidth]{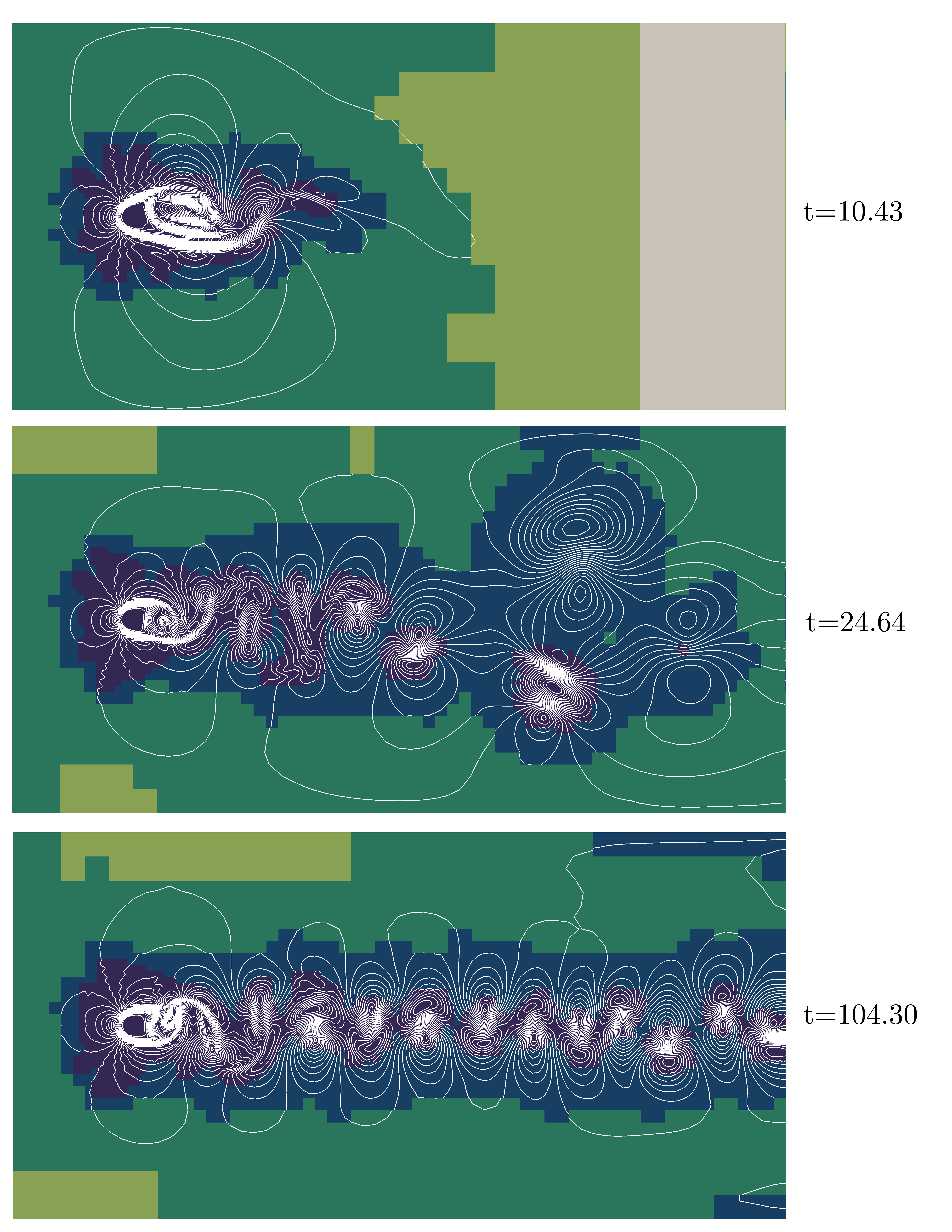}
    \end{center}\caption{\label{fig:1200_fs_2_7}Snapshots of the solution of the adaptive scheme for $\text{Re} = 1200$, $\minlevel = 2$, $\maxlevel = 7$, $\regularityguess = 1$ and $\thresholdletter = 0.00075$. The colors represent the levels of the mesh and the white contours are that of the velocity modulus.}
\end{figure}

\begin{figure}
    \begin{center}
        \includegraphics[width=1.\textwidth]{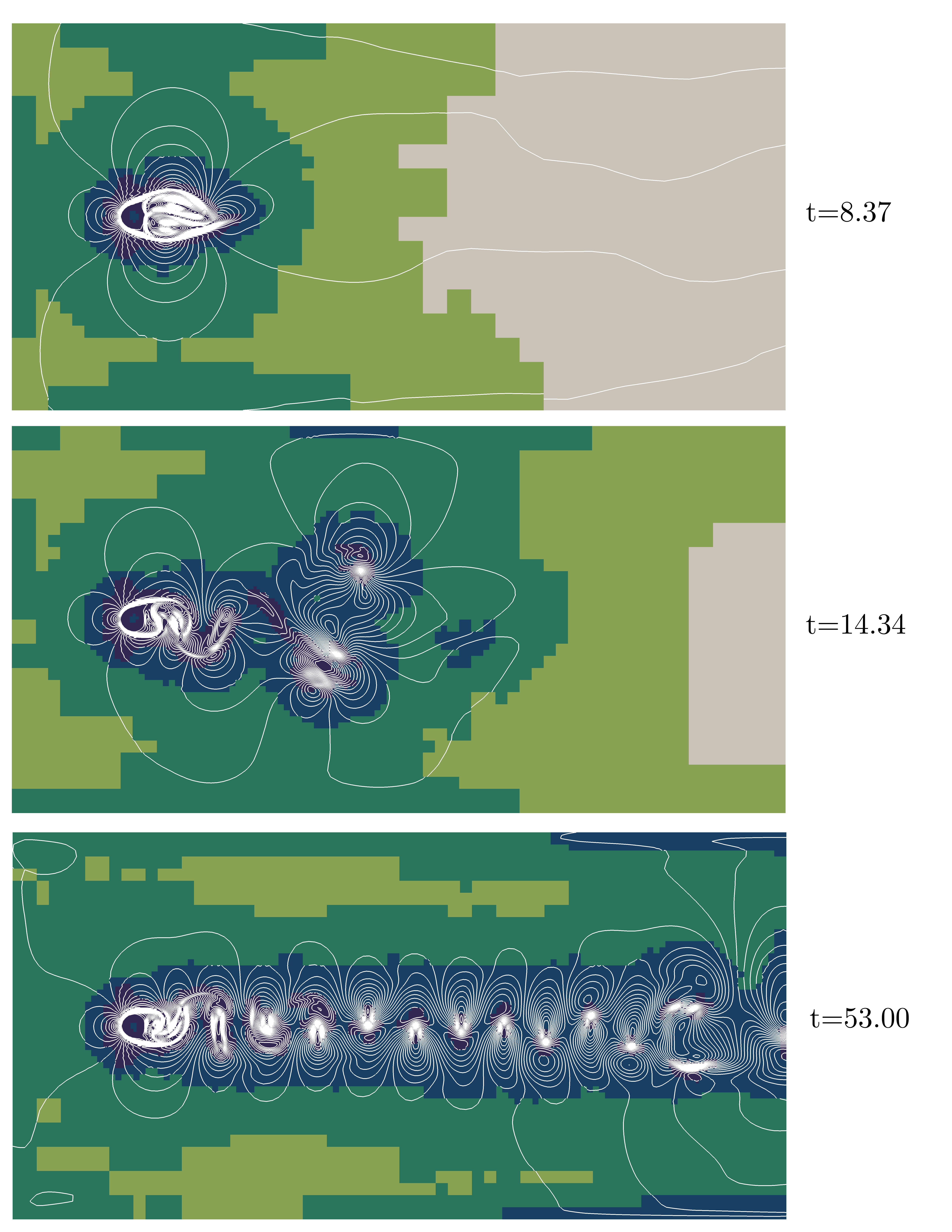}
    \end{center}\caption{\label{fig:1200_fs_2_8}Snapshots of the solution of the adaptive scheme for $\text{Re} = 1200$, $\minlevel = 2$, $\maxlevel = 8$, $\regularityguess = 1$ and $\thresholdletter = 0.000375$. The colors represent the levels of the mesh and the white contours are that of the velocity modulus.}
\end{figure}

\begin{figure}
    \begin{center}
        \includegraphics[width=1.\textwidth]{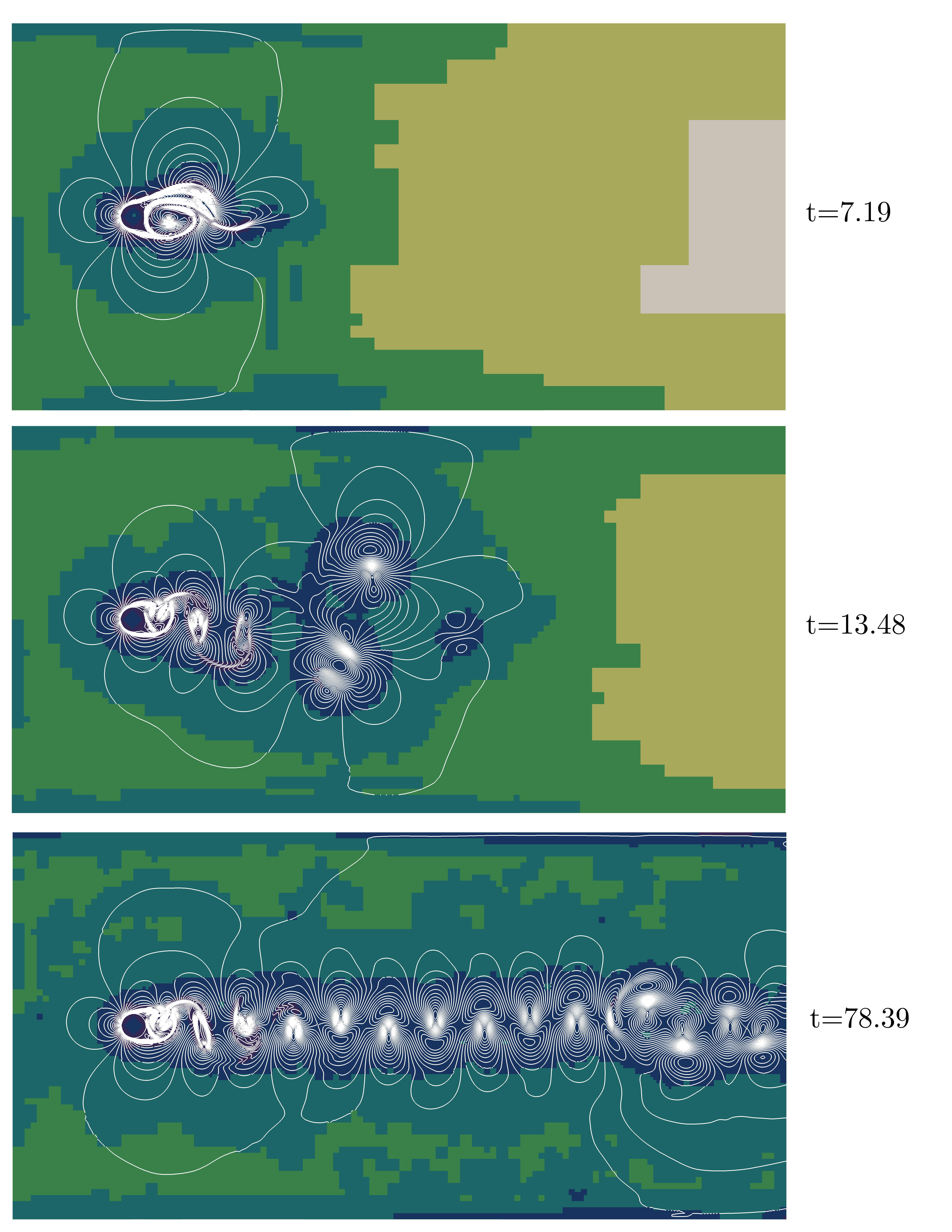}
    \end{center}\caption{\label{fig:1200_fs_2_9}Snapshots of the solution of the adaptive scheme for $\text{Re} = 1200$, $\minlevel = 2$, $\maxlevel = 9$, $\regularityguess = 1$ and $\thresholdletter = 0.0001875$. The colors represent the levels of the mesh and the white contours are that of the velocity modulus.}
\end{figure}

On the opposite side, far from the obstacle and close to the outlet, the regime can be considered to be mostly diffusive (or parabolic) and here the Harten heuristics could be violated. 
The purpose behind the division of the parameter $\thresholdletter$ by two at each time we increased the maximum level $\maxlevel$ was to try to follow the vortices until they reach the outlet with the finest possible resolution. As shown in Figures \ref{fig:1200_fs_2_7}, \ref{fig:1200_fs_2_8} and \ref{fig:1200_fs_2_9}, this attempt of following the vortices until the outlet with the finest resolution $\overline{J}$ was only partially successful.
More investigations are needed to clarify the role of the ``copy'' boundary conditions coupled with the multiresolution procedure.
This is not surprising by looking at \eqref{eq:DetailsDecayEstimate} and by concluding that the exact solution of the Navier-Stokes equations for this Reynolds number must have more than only one bounded derivative (\emph{i.e.} more than just bounded vorticity), because the details scale of a factor larger than two at each change of level. This is related to the nature of the solution and we refer to \cite{nguessan2020} for a detailed study of the use of multiresolution in order to resolve the incompressible Navier-Stokes equations. 
It is also worthwile mentioning that the D2Q9 scheme for the Navier-Stokes system under acoustic scaling $\Delta t \sim \Delta x$ is not converging for $\Delta x \to 0$: therefore one should be really careful once comparing the results for different $\maxlevel$, as we did.

To conclude, this test shows that our method is effective when applied to parabolic problems solved with the \lb method.
This reflects on the excellent results concerning the integral quantities.
Current investigations being conducted by our team using formal expansions \cite{dubois2009} show that our adaptive method does not modify (for $\predictionstencildepth \geq 1$) the viscosity of the flow, thus confirming once more that it is suitable to simulate the Navier-Stokes equations.
Furthermore, the method also preserves higher order terms and thus is less likely to modify the stability properties of the reference scheme.

\FloatBarrier

\subsection{Three dimensional advection equation}\label{sec:3dAdvection}

\subsubsection{Problem}

The only test we show for $\spatialdimension = 3$ aims at solving the transport equation with constant velocity $\bm{V} = (V_x, V_y, V_z)^T$ given by
\begin{equation*}
    \begin{cases}
     \partial_t u + \partial_x (V_x u) + \partial_y (V_y u) + \partial_z (V_z u) = 0, \qquad t \in [0, T], \qquad &\vectorial{x} \in \mathbb{R}^3, \\
     u(t = 0, \vectorial{x}) = \chi_{\sqrt{x^2 + y^2 + z^2} \leq 0.15}(\vectorial{x}), \qquad &\vectorial{x} \in \mathbb{R}^{3}.
    \end{cases}
\end{equation*}

\subsubsection{Numerical scheme}
We employ a D3Q6 scheme, thus $\velocitiesnumber = 6$, with one conserved moment with choice of velocities
\begin{equation*}
    \vectorial{\xi}^{\populationindex} = 
    \begin{cases}
     \latticevelocity \left  (\cos \left ( \pi \populationindex \right ), 0, 0 \right )^T, \qquad &\populationindex = 0, 1, \\
     \latticevelocity \left  (0, \cos \left ( \pi \populationindex \right ), 0 \right )^T, \qquad &\populationindex = 2, 3, \\
     \latticevelocity \left  (0, 0, \cos \left (\pi \populationindex \right ) \right )^T, \qquad &\populationindex = 4, 5, \\
    \end{cases}
\end{equation*}
with relaxation matrix $\operatorial{S} = \text{diag}(0, s_1, s_1, s_1, s_2, s_2)$ and change of basis given by
\begin{equation*}
    \operatorial{M} = 
    \begin{pmatrix}
     1 & 1 & 1 & 1 & 1 & 1\\
     \latticevelocity & -\latticevelocity & 0 & 0 & 0 & 0\\
     0 & 0 & \latticevelocity & -\latticevelocity & 0 & 0\\
     0 & 0 & 0 & 0 & \latticevelocity & -\latticevelocity\\
     \latticevelocity^2 & \latticevelocity^2 & - \latticevelocity^2 & -\latticevelocity^2 & 0 & 0\\
     \latticevelocity^2 & \latticevelocity^2 & 0 & 0 & -\latticevelocity^2 & -\latticevelocity^2
    \end{pmatrix}.
\end{equation*}
The vector of the equilibria is $\vectorial{m}^{\text{eq}} = (u, V_x u, V_y u, V_z u, 0, 0)^T$.
For this problem, the relaxation parameters are $s_1 = 1.4$ and $s_2 = 1$. We take $\latticevelocity = 1$.

\subsubsection{Results and discussion}

\begin{figure}
    \begin{center}
        \includegraphics[width=1.\textwidth]{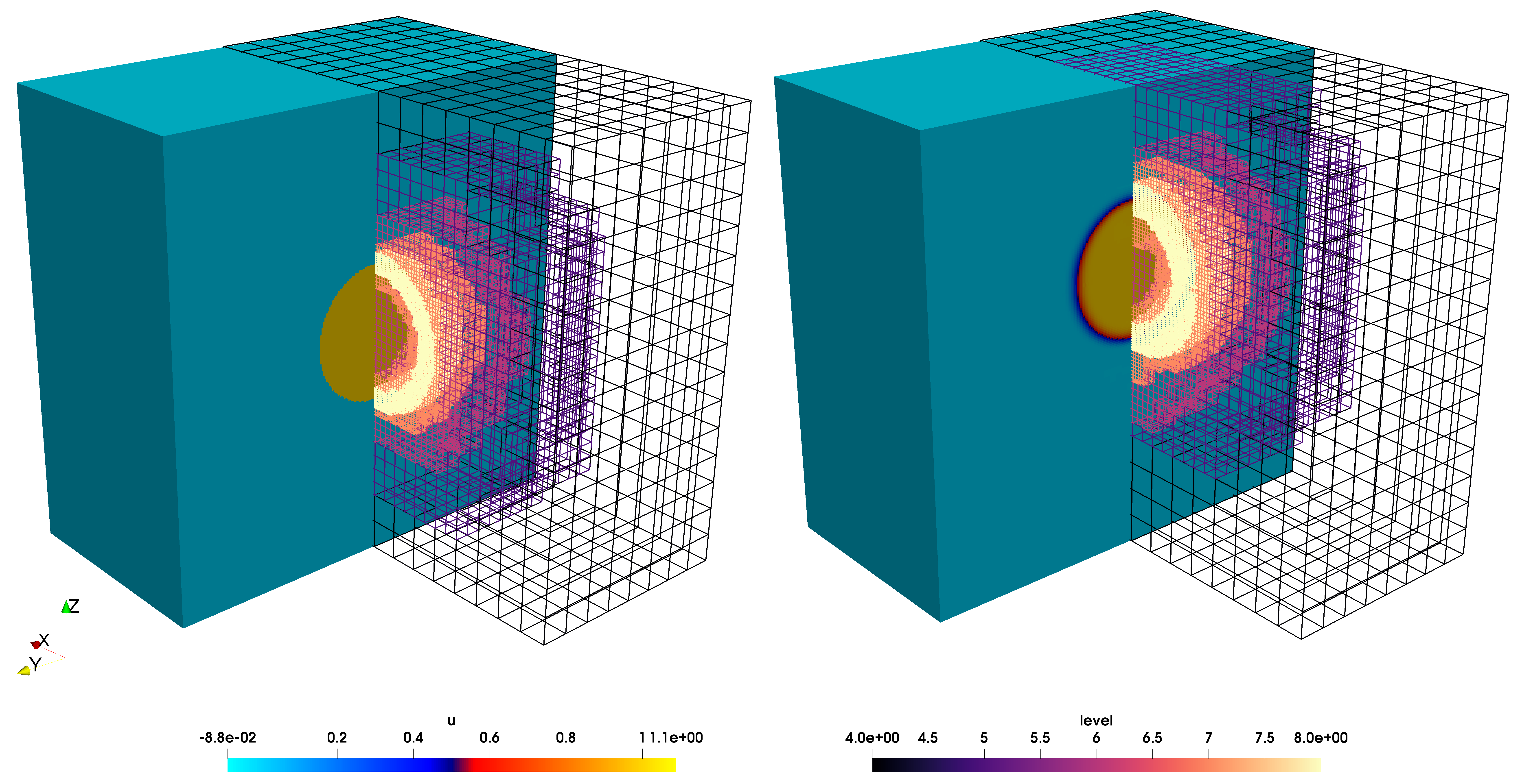}
    \end{center}\caption{\label{fig:3Dadvection}Snapshots at times $t = 0$ (left) and $t = 0.3125$ (right) of the solution of the adaptive scheme for the 3D advection equation. The domain has been cut so that the left half of each snapshots shows the value of the conserved variable $u$ whereas the right half represents the structure of the mesh with the corresponding local level of resolution.}
\end{figure}

\begin{figure}
    \begin{center}
        \includegraphics[width=0.5\textwidth]{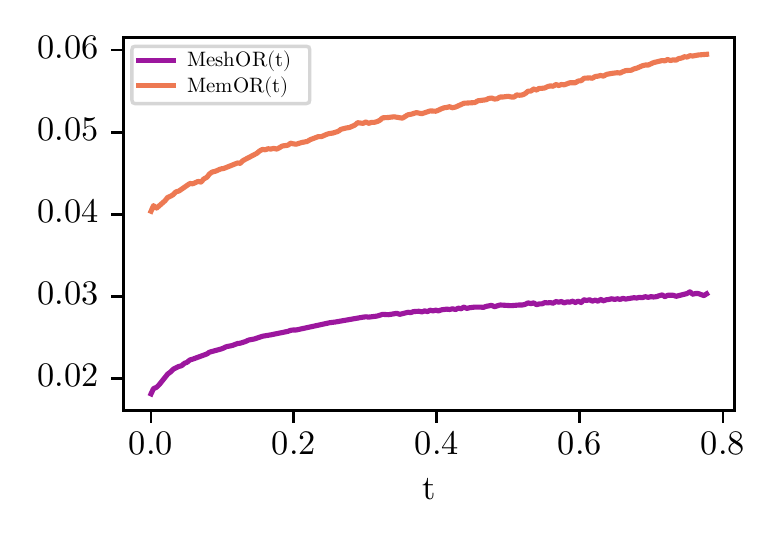}
    \end{center}\caption{\label{fig:3DadvectionCompression}Temporal variation of the occupation rates for the adaptive scheme simulating the 3D advection equation.}
\end{figure}

We carry on the simulation until $T = 0.78125$ with $\minlevel = 1$, $\maxlevel = 8$, $\regularityguess = 2$ and $\thresholdletter = 0.001$. Some snapshots of the solutions and the mesh at different time steps are provided in Figure \ref{fig:3Dadvection}: one can appreciate the significantly hollow mesh produced by multiresolution. The time behavior of the occupation rates is presented in Figure \ref{fig:3DadvectionCompression}, corresponding to an excellent compression rate of around $97$\%.
This shows that our method is capable of coping with this three dimensional problem properly and to achieve really interesting occupation rates which are crucial to conduct large simulations for 3D problems.
Indeed, the full problem would have needed $256^3 \sim 16$ millions cells, which is quite demanding for a sequential code. Thanks to multiresolution, we have made the problem easily treatable with at most $0.03 \times 256^3 \sim 50000$ cells involved (for all sizes). This means that the adaptive mesh occupies 3\% of its original size, which is quite impressive and shows the interest of the method to be used in 3D problems.

\section{Conclusions and future developments}\label{sec:Conclusions}

In this paper, we have introduced a novel way of performing mesh adaptation on multidimensional \lb schemes in order to reduce their memory footprint particularly on problems with steep fronts. 
This approach, based on adaptive multiresolution, modifies the mesh according to the local regularity of the solution and adapts the advection phase to reduce the total number of operations, while doing the collision only on the leaves of the adaptive tree.
Numerical experiments have shown that this technique can be combined with a wide range of existing schemes, both to simulate non-linear hyperbolic and non-linear parabolic systems. 
Therefore, the technique is problem independent and provides an effective cost reduction still being able to control the quality of the solution by means of a single threshold parameter $\thresholdletter$.

Several questions remain open and give room for improvement in future research works. First, one has to reach a real performance enhancement coming from dedicated forms of data storage and the parallelization of the code on various architectures.
Second, it would be interesting to investigate how to handle complex geometries, which are common in real problems, possibly using the numerous high-order boundary conditions present in the literature.
Detailed comparisons with the existing AMR approaches, even if always a difficult task,  could also be conducted with the tools we have developed.
Finally, we are currently investigating the possibility of constructing a more reliable collision phase following the lines by Hovhannisyan and M{\"u}ller \cite{hovhannisyan2010}.

\section*{Acknowledgments}

The authors deeply thank Laurent Séries and Christian Tenaud for fruitful discussions on multiresolution and the referee for having provided valuable inputs and suggestions to improve this work.
Thomas Bellotti is supported by a PhD funding (year 2019) from the Ecole polytechnique.

\bibliographystyle{acm}
\bibliography{bibliography.bib}

\section*{Appendix 1}

In this appendix, we provide the details on how recovering the fully discretized volumetric \lb scheme \eqref{eq:LBMReferenceSchemeNotSplit} starting from the discrete form of the Boltzmann equation \eqref{eq:BoltzmannEquation}.
Integrating this last equation on the characteristic of the velocity $\vectorial{\xi}^{\populationindex}$ gives
\begin{linenomath}\begin{align*}
     \int_{\timevariable}^{\timevariable + \Delta \timevariable} \adaptiveroundbrackets{\ratio{\partial \superscript{f}{\populationindex}}{\partial \timevariable} + \superscript{\vectorial{\xi}}{\populationindex} \cdot \nabla \superscript{f}{\populationindex}} \adaptiveroundbrackets{s, \vectorial{\spacevariable} + s \superscript{\vectorial{\xi}}{\populationindex}} \text{d} s &= - \sum_{l = 0}^{\velocitiesnumber - 1} \subscript{\omega}{hl} \int_{\timevariable}^{\timevariable + \Delta \timevariable} \adaptiveroundbrackets{\superscript{f}{l} - \superscript{f}{l, \text{eq}}}\adaptiveroundbrackets{s, \vectorial{\spacevariable} + s \superscript{\vectorial{\xi}}{\populationindex}} \text{d} s, \\
     &=\superscript{f}{\populationindex}(\timevariable + \Delta \timevariable, \vectorial{\spacevariable} + \Delta \timevariable \superscript{\vectorial{\xi}}{\populationindex}) - \superscript{f}{\populationindex}(\timevariable, \vectorial{\spacevariable}), \\
     &= -\frac{\Delta \timevariable}{2} \sum_{l = 0}^{\velocitiesnumber - 1} \subscript{\omega}{hl}  \adaptiveroundbrackets{\superscript{f}{l} - \superscript{f}{l, \text{eq}}}\adaptiveroundbrackets{\timevariable + \Delta \timevariable, \vectorial{\spacevariable} + \Delta \timevariable \superscript{\vectorial{\xi}}{\populationindex}} \\
     &\phantom{= } -\frac{\Delta \timevariable}{2} \sum_{l = 0}^{\velocitiesnumber - 1} \subscript{\omega}{hl}  \adaptiveroundbrackets{\superscript{f}{l} - \superscript{f}{l, \text{eq}}}\adaptiveroundbrackets{\timevariable, \vectorial{\spacevariable}},
\end{align*}\end{linenomath}
where the collision terms has been approximated by a trapezoidal rule, see Dellar \cite{dellar2013}.
Calling $\superscript{\tilde{f}}{\populationindex} = \superscript{f}{\populationindex} + \frac{\Delta \timevariable}{2} \sum_{l = 0}^{\velocitiesnumber - 1} \subscript{\omega}{hl}  (\superscript{f}{l} - \superscript{f}{l, \text{eq}})$ to keep the method explicit and rearranging gives
\begin{linenomath}\begin{equation*}
  \superscript{\tilde{f}}{\populationindex} \adaptiveroundbrackets{\timevariable + \Delta \timevariable, \vectorial{\spacevariable} + \Delta \timevariable \superscript{\vectorial{\xi}}{\populationindex}} = \superscript{\tilde{f}}{\populationindex} \adaptiveroundbrackets{\timevariable, \vectorial{\spacevariable}} -\Delta \timevariable \sum_{l = 0}^{\velocitiesnumber - 1} \subscript{\omega}{hl}  \adaptiveroundbrackets{\superscript{f}{l} - \superscript{f}{l, \text{eq}}}\adaptiveroundbrackets{\timevariable, \vectorial{\spacevariable}}.
\end{equation*}\end{linenomath}
We are left to dealing with the last term on the right hand side of the previous equation
\begin{linenomath}\begin{equation*}
    \superscript{f}{l} - \superscript{f}{l, \text{eq}} = \superscript{\tilde{f}}{l} - \frac{\Delta \timevariable}{2} \sum_{m = 0}^{\velocitiesnumber - 1} \subscript{\omega}{lm}  (\superscript{f}{m} - \superscript{f}{m, \text{eq}})  - \superscript{f}{l, \text{eq}}, 
\end{equation*}\end{linenomath}
which also reads, gathering the collision frequencies in a matrix $\operatorial{\omega}$ with entries $\omega_{ij}$
\begin{linenomath}\begin{equation*}
    (\operatorial{I} + \Delta \timevariable \operatorial{\omega}/2) (\vectorial{f} - \vectorial{f}^{\text{eq}}) = \tilde{\vectorial{f}} - \vectorial{f}^{\text{eq}}, \qquad \text{thus} \quad \vectorial{f} - \vectorial{f}^{\text{eq}} = (\operatorial{I} + \Delta \timevariable \operatorial{\omega}/2)^{-1}(\tilde{\vectorial{f}} - \vectorial{f}^{\text{eq}})
\end{equation*}\end{linenomath}
Therefore, the semi-discretized numerical scheme reads
\begin{linenomath}\begin{align*}
\tilde{{f}}^{\populationindex} \adaptiveroundbrackets{\timevariable + \Delta \timevariable, \vectorial{\spacevariable} + \Delta \timevariable \superscript{\vectorial{\xi}}{\populationindex}} &= {\tilde{{f}}^{\populationindex} \adaptiveroundbrackets{\timevariable, \vectorial{\spacevariable}} -\Delta \timevariable \operatorial{\omega}  (\vectorial{f} - \vectorial{f}^{\text{eq}})|_{\populationindex}\adaptiveroundbrackets{\timevariable, \vectorial{\spacevariable}},} \\
  &= {\tilde{{f}}^{\populationindex} \adaptiveroundbrackets{\timevariable, \vectorial{\spacevariable}} -\Delta \timevariable \operatorial{\omega}  (\operatorial{I} + \Delta \timevariable \operatorial{\omega}/2)^{-1}(\tilde{\vectorial{f}} - \vectorial{f}^{\text{eq}})|_{\populationindex}\adaptiveroundbrackets{\timevariable, \vectorial{\spacevariable}}.}
\end{align*}\end{linenomath}
Ahead, the tilde sign is understood. For the discretization of the spatial variable, consider a cell $\subscript{\cellletter}{\maxlevel, \vectorial{\indexletter}}$ for any $\vectorial{k} \in \mathbb{Z}$ at the finest level of discretization. Taking the average yields
\begin{linenomath}\begin{align*}
    {
    \ratio{1}{|\subscript{\cellletter}{\maxlevel, \vectorial{\indexletter}}|_{\spatialdimension}} \int_{\subscript{\cellletter}{\maxlevel, \vectorial{\indexletter}}} {{f}}^{\populationindex} \adaptiveroundbrackets{\timevariable + \Delta \timevariable, \vectorial{\spacevariable} + \Delta \timevariable \superscript{\vectorial{\xi}}{\populationindex}} \text{d} \vectorial{\spacevariable}} &=  {\ratio{1}{|\subscript{\cellletter}{\maxlevel, \vectorial{\indexletter}}|_{\spatialdimension}} \int_{\subscript{\cellletter}{\maxlevel, \vectorial{\indexletter}}} {f}^{\populationindex} \adaptiveroundbrackets{\timevariable, \vectorial{\spacevariable}} \text{d} \vectorial{\spacevariable}}  \\
    &{-  \ratio{\Delta \timevariable \operatorial{\omega}  (\operatorial{I} + \Delta \timevariable \operatorial{\omega}/2)^{-1}}{|\subscript{\cellletter}{\maxlevel, \vectorial{\indexletter}}|_{\spatialdimension}} \int_{\subscript{\cellletter}{\maxlevel, \vectorial{\indexletter}}} \adaptiveroundbrackets{\vectorial{f} - \superscript{\vectorial{f}}{\text{eq}}} \adaptiveroundbrackets{\timevariable, \vectorial{\spacevariable}} \text{d} \vectorial{\spacevariable}\hspace{0.2cm} |_{\populationindex}.}
\end{align*}\end{linenomath}
A change of variable in the first integral provides, indicating the averages with a bar
\begin{linenomath}\begin{equation*}
    {\subscript{\overline{{f}}}{\maxlevel, \vectorial{\indexletter} + \superscript{\vectorial{\eta}}{\populationindex}}^{\populationindex} \adaptiveroundbrackets{\timevariable + \Delta \timevariable} = \subscript{\overline{{f}}}{\maxlevel, \vectorial{\indexletter}}^{\populationindex} \adaptiveroundbrackets{\timevariable}  - \Delta \timevariable \operatorial{\omega}  (\operatorial{I} + \Delta \timevariable \operatorial{\omega}/2)^{-1} \adaptiveroundbrackets{\overline{\vectorial{f}}_{\maxlevel, \vectorial{\indexletter}}(\timevariable) - \ratio{1}{|\subscript{\cellletter}{\maxlevel, \vectorial{\indexletter}}|_{\spatialdimension}} \int_{\subscript{\cellletter}{\maxlevel, \vectorial{\indexletter}}} \superscript{\vectorial{f}}{\text{eq}}\adaptiveroundbrackets{\superscript{f}{0}(\timevariable, \vectorial{\spacevariable}), \dots, \superscript{f}{q-1}(\timevariable, \vectorial{\spacevariable})}  \text{d} \vectorial{\spacevariable}} |_{\populationindex}.}
\end{equation*}\end{linenomath}
We commute the integral operator and the equilibrium operator, thus yielding the fully discrete scheme
\begin{linenomath}\begin{equation*}
    {\subscript{\overline{{f}}}{\maxlevel, \vectorial{\indexletter} + \superscript{\vectorial{\eta}}{\populationindex}}^{\populationindex} \adaptiveroundbrackets{\timevariable + \Delta \timevariable} = \subscript{\overline{{f}}}{\maxlevel, \vectorial{\indexletter}}^{\populationindex} \adaptiveroundbrackets{\timevariable}  - \Delta \timevariable \operatorial{\omega}  (\operatorial{I} + \Delta \timevariable \operatorial{\omega}/2)^{-1} \adaptiveroundbrackets{\overline{\vectorial{f}}_{\maxlevel, \vectorial{\indexletter}}(\timevariable) - \superscript{\vectorial{f}}{\text{eq}}\adaptiveroundbrackets{\superscript{\overline{f}}{0}_{\maxlevel, \vectorial{\indexletter}}(\timevariable), \dots, \superscript{\overline{f}}{q-1}_{\maxlevel, \vectorial{\indexletter}}(\timevariable)} } |_{\populationindex}.}
\end{equation*}\end{linenomath}

{
\section*{Appendix 2}}
{
In this appendix, we briefly recall how the prediction coefficient presented in Table \ref{tab:PredictionCoefficients} are recovered, following \cite{cohen2003} and \cite{duarte2011}.
In a one-dimensional setting, for a cell $C_{\levelletter, \indexletter}$, one considers the local reconstruction polynomial of degree $2\predictionstencildepth$ written in the canonical basis
\begin{equation*}
    \pi_{\levelletter, \indexletter}^{\populationindex} (\spacevariable) = \sum_{m = 0}^{2\predictionstencildepth} A_{\levelletter, \indexletter}^{\populationindex, m} \spacevariable^m.
\end{equation*}
The coefficients $(A_{\levelletter, \indexletter}^{\populationindex, m})_{m = 0}^{m = 2\predictionstencildepth} \subset \mathbb{R}$ are obtained enforcing the linear constraints
\begin{equation*}
 \frac{1}{\Delta x_{\levelletter}} \int_{C_{\levelletter, \indexletter + \delta}} \pi_{\levelletter, \indexletter}^{\populationindex} (\spacevariable) \text{d}\spacevariable = f_{\levelletter, \indexletter + \delta}^{\populationindex}, \qquad \delta = -\predictionstencildepth, \dots, \predictionstencildepth.
\end{equation*}
This signifies that the means of the reconstruction polynomial $\pi_{\levelletter, \indexletter}^{\populationindex}$ on the $2\predictionstencildepth$ neighbors of $C_{\levelletter, \indexletter}$, plus the cell itself, must yield the actual mean value on these cells. 
This gives a linear system with matrix $\operatorial{T} \in \text{GL}_{2\predictionstencildepth + 1}(\mathbb{R})$
\begin{equation}\label{eq:systemPrediction}
    \operatorial{T} (A_{\levelletter, \indexletter}^{\populationindex, m})_{m = 0}^{m = 2\predictionstencildepth} = (f_{\levelletter, \indexletter + \delta}^{\populationindex})_{\delta = -\predictionstencildepth}^{\delta = +\predictionstencildepth}.
\end{equation}
Notice that $\operatorial{T}$ does not depend on $\levelletter$, nor on $\indexletter$, nor on $\populationindex$.
Once the system is solved and thus the reconstruction polynomial is known, the prediction is obtained by averaging the local reconstruction polynomial over the siblings $C_{\levelletter + 1, 2\indexletter}$ and $C_{\levelletter + 1, 2\indexletter + 1}$
\begin{equation*}
    \predicted{f}{\populationindex}{\levelletter + 1, 2\indexletter + \delta} = \frac{1}{\Delta x_{\levelletter + 1}} \int_{C_{\levelletter + 1, 2\indexletter + \delta}} \pi_{\levelletter, \indexletter}^{\populationindex} (\spacevariable) \text{d}\spacevariable = \superscript{\subscript{f}{\levelletter, \indexletter}}{\populationindex} + \power{\adaptiveroundbrackets{-1}}{\delta} \subscript{\superscript{\predictionnoncenteredletter}{\predictionstencildepth}}{1} \adaptiveroundbrackets{\indexletter; \superscript{\subscript{\vectorial{f}}{\levelletter}}{\populationindex}}, \qquad \delta = 0, 1,
\end{equation*}
with $\subscript{\superscript{\predictionnoncenteredletter}{\predictionstencildepth}}{1} \adaptiveroundbrackets{\indexletter; \superscript{\subscript{\vectorial{f}}{\levelletter}}{\populationindex}} \definitionequality \sum_{\alpha = 1}^{\predictionstencildepth} c_{\alpha} \adaptiveroundbrackets{\superscript{\subscript{f}{\levelletter, \indexletter + \alpha}}{\populationindex} - \superscript{\subscript{f}{\levelletter, \indexletter - \alpha}}{\populationindex}}$, where one can check that the weights are those given in Table \ref{tab:PredictionCoefficients}. }

{
For $\spatialdimension > 1$ the previous procedure is generalized by tensor product \cite{bihari1997multiresolution}, thus considering an analogous of system \eqref{eq:systemPrediction} where the matrix is the Kronecker product ($\spatialdimension$-times) of $\operatorial{T}$ with itself.
}

\end{document}